\documentclass[12pt]{amsart}
\usepackage{amsbsy}
\usepackage{graphicx,epsfig,subfigure}%,psfrag}%,psfrag}
\usepackage[scanall]{psfrag}
\begin{document}

%%%%%%%%%%%%%%%%%%%%%%%%%%%%%%%%%%%%%%%%%%%%%%%%%%%%%%%%%%%%%%%%%%%%
% Theorem, definition, lemma, proposition, corollary and proof
%%%%%%%%%%%%%%%%%%%%%%%%%%%%%%%%%%%%%%%%%%%%%%%%%%%%%%%%%%%%%%%%%%%%
%%%%%%%%%%%%%%%%%%%%%%%%%%%%%%%%%%%%%%%%%%%%%%%%
\newtheorem{theorem}{Theorem}%[section]
\newtheorem{proposition}{Proposition}%[section]
\newtheorem{lemma}{Lemma}%[section]
\newtheorem{corollary}{Corollary}%[section]
\newtheorem{definition}{Definition}%[section]
\newtheorem{remark}{Remark}%[section]
%%%%%%%%%%%%%%%%%%%%%%%%%%%%%%%%%%%
%%%%%%%%%%%%%%%%%%%%%%%%%%%%%%%%%%%%%%%%%%%%%%                  NEW
%%\newcommand{\be}{\begin{equation}}
%%\newcommand{\ee}{\end{equation}}
%%%%%%%%%%%%%%%%%%%%%%%%%%%%%%%%%%%%%%%%%%%%%%%
%%%%%%%%%%%%%%%%%%%%%%%%%%%%%%%%%%%%%%%%%%%%
\newcommand{\tex}{\textstyle}
%%%%%%%%%%%%%%%%%%%%%%%%%%%%%%%%%%%%%%%%%%%%
%%%%%%%%%%%%%%%%%%%%%%%%%%%%%%%%%%%%%%%%%%%%
\numberwithin{equation}{section} \numberwithin{theorem}{section}
\numberwithin{proposition}{section} \numberwithin{lemma}{section}
\numberwithin{corollary}{section}
\numberwithin{definition}{section} \numberwithin{remark}{section}
%%%%%%%%%%%%%%%%%%%%%%%%%%%%%%%%%%%%%%%%%%%%
%%%%%%%%%%%%%%%%%%%%%%%%%%%%%%%%%%%%%%%%%%%%
\newcommand{\ren}{\mathbb{R}^N}
\newcommand{\re}{\mathbb{R}}
\newcommand{\n}{\nabla}
\newcommand{\iy}{\infty}
\newcommand{\pa}{\partial}
\newcommand{\fp}{\noindent}
\newcommand{\ms}{\medskip\vskip-.1cm}
\newcommand{\mpb}{\medskip}
%%%%%%%%%%%%%%%%%%%%%%%%%%%%%%%%%%%%%%%%%%%%%%%%%
\newcommand{\AAA}{{\bf A}}
\newcommand{\BB}{{\bf B}}
\newcommand{\CC}{{\bf C}}
\newcommand{\DD}{{\bf D}}
\newcommand{\EE}{{\bf E}}
\newcommand{\FF}{{\bf F}}
\newcommand{\GG}{{\bf G}}
\newcommand{\oo}{{\mathbf \omega}}
\newcommand{\Am}{{\bf A}_{2m}}
\newcommand{\CCC}{{\mathbf  C}}
\newcommand{\II}{{\mathrm{Im}}\,}
\newcommand{\RR}{{\mathrm{Re}}\,}
\newcommand{\eee}{{\mathrm  e}}
%%%%%%%%%%%%%%%%%%%%%%%%%%%%%%%%%%%%%%%%%%%%%%%%%%%%%%%%%%%%%%%%%%%%%%% L^2\rho...
\newcommand{\LL}{L^2_\rho(\ren)}
\newcommand{\LLL}{L^2_{\rho^*}(\ren)}
%%%%%%%%%%%%%%%%%%%%%%%%%%%%%%%%%%
%%%%%%%%%%%%%%%%%%%%%%%%%%%%%%%%%%%%%%%%%%%%%%%%%%%%
\renewcommand{\a}{\alpha}
\renewcommand{\b}{\beta}
\newcommand{\g}{\gamma}
\newcommand{\G}{\Gamma}
\renewcommand{\d}{\delta}
\newcommand{\D}{\Delta}
\newcommand{\e}{\varepsilon}
\newcommand{\var}{\varphi}
\newcommand{\lll}{\l}
\renewcommand{\l}{\lambda}
\renewcommand{\o}{\omega}
\renewcommand{\O}{\Omega}
\newcommand{\s}{\sigma}
\renewcommand{\t}{\tau}
\renewcommand{\th}{\theta}
\newcommand{\z}{\zeta}
\newcommand{\wx}{\widetilde x}
\newcommand{\wt}{\widetilde t}
\newcommand{\noi}{\noindent}
 %%%%%%%%%%%%%%%%%%%%%%%%%%%%%%%%%%%%%%%%%%%
\newcommand{\uu}{{\bf u}}
\newcommand{\xx}{{\bf x}}
\newcommand{\yy}{{\bf y}}
\newcommand{\zz}{{\bf z}}
\newcommand{\aaa}{{\bf a}}
\newcommand{\cc}{{\bf c}}
\newcommand{\jj}{{\bf j}}
\newcommand{\ggg}{{\bf g}}
\newcommand{\UU}{{\bf U}}
\newcommand{\YY}{{\bf Y}}
\newcommand{\HH}{{\bf H}}
\newcommand{\GGG}{{\bf G}}
\newcommand{\VV}{{\bf V}}
\newcommand{\ww}{{\bf w}}
\newcommand{\vv}{{\bf v}}
\newcommand{\hh}{{\bf h}}
\newcommand{\di}{{\rm div}\,}
\newcommand{\ii}{{\rm i}\,}
%%%%%%%%%%%%%%%%%%%%%%%%%%%%%%%%%%
%%%%%%%%%%%%%%%%%%%%%%%%%%%%%%%%%%%%%   VAG, NEW
\newcommand{\inA}{\quad \mbox{in} \quad \ren \times \re_+}
\newcommand{\inB}{\quad \mbox{in} \quad}
\newcommand{\inC}{\quad \mbox{in} \quad \re \times \re_+}
\newcommand{\inD}{\quad \mbox{in} \quad \re}
\newcommand{\forA}{\quad \mbox{for} \quad}
\newcommand{\whereA}{,\quad \mbox{where} \quad}
\newcommand{\asA}{\quad \mbox{as} \quad}
\newcommand{\andA}{\quad \mbox{and} \quad}
\newcommand{\withA}{,\quad \mbox{with} \quad}
\newcommand{\orA}{,\quad \mbox{or} \quad}
\newcommand{\atA}{\quad \mbox{at} \quad}
\newcommand{\onA}{\quad \mbox{on} \quad}
\newcommand{\ef}{\eqref}
\newcommand{\ssk}{\smallskip}
\newcommand{\LongA}{\quad \Longrightarrow \quad}
%%%%%%%%%%%%%%%%%%%%%%%%%%%%%%%%
%%%%%%%%%%%%%%%%%%%%%%%%%%%%%%%%%%
\def\com#1{\fbox{\parbox{6in}{\texttt{#1}}}}
%%%%%%%%%%%%%%%%%%%%%%%%%%%%%%%%%%
%%%%%%%%%%%%%%%%%%% From Paper1
\def\N{{\mathbb N}}
\def\A{{\cal A}}
\newcommand{\de}{\,d}
\newcommand{\eps}{\varepsilon}
\newcommand{\be}{\begin{equation}}
\newcommand{\ee}{\end{equation}}
\newcommand{\spt}{{\mbox spt}}
\newcommand{\ind}{{\mbox ind}}
\newcommand{\supp}{{\mbox supp}}
\newcommand{\dip}{\displaystyle}
\newcommand{\prt}{\partial}
\renewcommand{\theequation}{\thesection.\arabic{equation}}
\renewcommand{\baselinestretch}{1.1}
%%%%%%%%%%%%%%%%%%%%%%%%%%%%%%%%%%%%%%%%%%%%%%%
\newcommand{\Dm}{(-\D)^m}

%%%%%%%%%%%%%%%%%%%%%%%%% VICTOR
\title
%%%%%[Positivity]
%%%%%%%%%%%%%%%%%%%%%%%%%
%%%%{\bf How  $\sqrt{\mbox{log\,log}}$ factor occurs in blow-up
{\bf Very singular similarity solutions and Hermitian spectral
theory for semilinear odd-order PDEs}

%%towards Petrovskii-like criterion by blow-up spectral
%%theory}
%%%    $\sqrt{\mbox{log\,log}}$
%%\sqrt{\ln|\ln(T-t)|}}
%% factor occurs}
%% solutions of higher-order
%%%nonlinear Schr\"odinger equations}

%% applications to
%%%large time and blow-up nodal set evolution}

\author {R.S.~Fernandes and V.A.~Galaktionov}

\address{Department of Mathematical Sciences, University of Bath,
 Bath BA2 7AY, UK}
%%% and Keldysh Institute of Applied Mathematics,
%%%% Miusskaya Sq. 4, 125047 Moscow, RUSSIA}
\email{vag@maths.bath.ac.uk}

\address{Department of Mathematical Sciences, University of Bath,
 Bath BA2 7AY, UK}
%%% and Keldysh Institute of Applied Mathematics,
%%%% Miusskaya Sq. 4, 125047 Moscow, RUSSIA}
\email{rsf21@maths.bath.ac.uk}

%%\address{Department of Mathematical Sciences, University of Bath,
%% Bath BA2 7AY, UK}
%%\email{ivk20@maths.bath.ac.uk}

%%\thanks{Research supported by  RTN network
%%HPRN-CT-2002-00274 and CERN-INTAS00-0136}

\keywords{Odd-order linear and semilinear PDEs,
 fundamental solution, Hermitian spectral theory, polynomial eigenfunctions,
 self-similarity,  very singular solutions, bifurcations, branching.
 %%%discrete spectrum.
%% $\sqrt{\ln|\ln(T-t)|}$-factor
%%% global existence, uniform bounds
%% discrete real spectrum, complete set of
%%%eigenfunctions, asymptotic behaviour.
 %% \\
 %%{\bf Submitted to:}
  }

\date{\today}
%%% \quad {\bf Ray/RayI.tex}}
%%%%%%%%%%%%%%%%%%%%%%%%%%%%%%%%%%%%%%%%%%%%%%%%%%%%%%%%

\pagenumbering{arabic}

\begin{abstract}

%%The text of the abstract goes here.

Asymptotic large- and short-time behaviour  of solutions of the
{\em linear dispersion equation}
 $$
 u_t=u_{xxx} \quad \mbox{in}
 \quad \re \times \re_+,
 $$
and its $(2k+1)$th-order extensions are studied. Such a refined
scattering is based on a ``Hermitian"
 spectral theory for a pair $\{\BB,\BB^*\}$ of non self-adjoint rescaled
 operators
  $$
   \mbox{$
  {\bf B}= D_y^3 + \frac 13\, y D_y + \frac 13\, I, \,\,\,
 \mbox{and the adjoint one} \,\,\,
   \BB^*=D_y^3- \frac 13 \, y D_y
   %%% \mbox{with} \quad \sigma({\bf B})=\big\{\l_l=-\frac l3, \,\, l=0,1,2,...\,\big\},
 $}
$$
  with the discrete spectrum $\sigma({\bf
B})=\s(\BB^*)=\big\{\l_l=-\frac l3, \,\, l=0,1,2,...\,\big\}$
 and eigenfunctions $\{\psi_l(y)= \frac {(-1)^l}{\sqrt {l!}} D^l_y {\rm Ai}(y), \, l \ge 0\}$, where ${\rm Ai}(y)$ is Airy's
  classic function.
 %%% and for its adjoint one $\BB^*=D_y^3- \frac 13 \, y D_y$.
 %%% is developed.
 Applications to {\em very singular similarity solutions} (VSSs)
  of the
  semilinear dispersion equation
with absorption,
  $$
   \tex{
  u_{\rm S}(x,t)=t^{- \frac 1{p-1}}f\big(\frac x{t^{1/3}}\big)
  }: \quad
%%   $$
 %%$$
 u_t=u_{xxx}-|u|^{p-1}u \quad \mbox{in}
 \quad \re \times \re_+, \quad
 %% \mbox{where}
 %%\quad
 p>1,
 $$
 and to its higher-order counterparts are presented. The goal is, by using
 various techniques, to show that there exists a countable
 sequence of critical exponents $\{p_l=1+ \frac 3{l+1}, \,
 l=0,1,2,...\}$ such that, at each $p=p_l$, a $p$-branch of VSSs
 bifurcates from the corresponding eigenfunction $\psi_l$ of the linear operator $\BB$ above.

 %%which serves as a basic model for various applications including
 %%the classic KdV area.

\end{abstract}

\maketitle
%\newpage

%\begin{acknowledgements}
%Acknowledgements
%\end{acknowledgements}

%%\tableofcontents

%\numberwithin{equation}{section}

%%%\newpage \pagenumbering{arabic}

%%%%%\section{Introduction}

%%%%%%%%%%%%%%%%%%%%%%%%%%%%%%%%%%%%%%%%%%%%%%%%%%%%%%%%%%%%%%%%%%%%%%
\section{Introduction: semilinear odd-order models, history, and results}
 \label{S1}

\subsection{Basic dispersion models and applications}

As a first basic model,  we will study  higher odd-order partial
differential equations (PDEs) of the form
\begin{equation}
\label{generalpde}
u_t =(-1)^{k+1}D_x^{2k+1} u + \tilde{g}(u)
\quad \text{in} \quad \mathbb{R} \times \mathbb{R_+}, \quad
k=1,2,...\, ,
\end{equation}
with bounded integrable  initial data $u(x,0) = u_0(x)$ in $\re$.
Here $D_x^m =(\partial/\partial x)^m$ denotes
 the $m^{\rm th}$ partial derivative in $x$.
 %%% and $u_t$ denotes the partial derivative of $u(x,t)$ with
 %%respect to the time variable $t$.
 The odd-order semilinear
dispersion equation \eqref{generalpde} can be considered as  a
counterpart of the better known semilinear higher-even-order {\em
parabolic} PDE of {\em reaction-diffusion} type,
\begin{equation}
\label{par1}
 u_t =(-1)^{k+1}D_x^{2k} u + \tilde{g}(u) \quad
\text{in} \quad \mathbb{R} \times \mathbb{R_+}, \quad k=1,2,...\,
.
\end{equation}
 For $k=1$, (\ref{par1}) becomes a standard reaction-diffusion
 equation from combustion theory
\begin{equation}
\label{par12}
 u_t = u_{xx} + \tilde{g}(u) \quad
\text{in} \quad \mathbb{R} \times \mathbb{R_+},
\end{equation}
 to which dozens of well-known monographs are devoted to.
These parabolic equations indeed belongs to an entirely different
type of PDEs and were much better studied in the twentieth
century. However,
 the analogy
between odd and even-order PDEs, such as \eqref{generalpde} and
\eqref{par1}, is rather fruitful and will be used later on.

The function (a nonlinear operator) $ \tilde {g}(u)$ in
(\ref{par1}) usually corresponds to some absorption-reaction type
phenomena and sometimes is assumed to include differential terms,
such as $D_x^{{m}}u$, with ${m} < 2k+1$ (although we do not
consider such cases).
%%Note that
%%any constant coefficients may be placed in front of any of the
%%terms, as these can easily be scaled out to obtain a PDE of the
%%same form.
 It is worth mentioning again that, besides some special and
completely integrable PDEs, general odd-order models such as
\eqref{generalpde} are  less studied in the mathematical
literature, than the parabolic even-order ones \eqref{par1}.

Indeed, the most classical  example of such an odd-order equation
is the {\em KdV equation}:
\begin{equation}
\label{KdV}
 u_t = u_{xxx} + uu_x \quad(\tilde g(u):=u u_x),
\end{equation}
which was  introduced by Boussinesq  in 1872 together with its
soliton solution \cite{Bouss}.
  The KdV equation models long waves in
shallow water and generates a hierarchy of other more complicated
PDEs with linear and nonlinear dispersion (dispersive) mechanisms.
 See further  amazing
historical aspects concerning \ef{KdV} and related integrable PDEs
in \cite[pp.~226-229]{GSVR}.

%%{It has long been recognized
   %%%%%% for a long time
%%  that both the KdV equation and its soliton solution
%% were  derived  earlier by  Boussinesq in 1872 \cite{Bouss2}, so
%% the abbreviation {\em BKdV} can be used for
%%% (\ref{KV.1}).}

Concerning higher-order extensions,
%% those are also well known in
%%different hierarchies of integrable PDEs.
 %%% Talking about odd-order PDEs under consideration,
%% \noi\underline{\em On
%%integrable NDEs}.   Concerning applications of equations such as
%%(\ref{1}),
%%%let us mention that
%%%we begin with the well-known fact
 %%various odd-order PDEs
 %%%are a natural quasilinear counterpart of
  these naturally appear in classic theory of integrable
 PDEs from shallow water applications,  including
 the {\em fifth-order KdV equation},
 $$
 u_t + u_{xxxxx} + 30 \, u^2 u_x + 20 \, u_x u_{xx} + 10 \, u
 u_{xxx}=0.
%%%  \quad \mbox{in} \,\,\,\, \re \times \re,
 $$
Let us also mention such classic examples as
 %%We can enlarge  this list talking about possible
 %%  quasilinear  extensions of the integrable
 {\em Lax's
 seventh-order KdV equation}
  %%%\index{equation!Lax 7th-order KdV}  %%%%%%%%%%%%%%%%%%%%%%%%%%%%%%%%%%%%%%% ???? Where
 $$
 u_t+[35 u^4 + 70(u^2 u_{xx}+ u(u_x)^2)+7(2 u u_{xxxx}+3 (u_{xx})^2 +
 4 u_x u_{xxx})+u_{xxxxxx}]_x=0,
  $$
  and the {\em seventh-order Sawada--Kotara
  equation}
   %%%%%\index{equation!Sawada-Kotara 7th-order}
 $$
 u_t+[63 u^4 + 63(2u^2 u_{xx}+ u(u_x)^2)+21( u u_{xxxx}+ (u_{xx})^2 +
  u_x u_{xxx})+u_{xxxxxx}]_x=0;
  $$
see  \cite[p.~234]{GSVR} for references and other odd-order models
from compacton theory.

\ssk

Of course, PDEs such as \eqref{KdV} and \eqref{generalpde} were
popular in mathematical physics last fifty years at least (see
references below), are have been well-known in mathematical
theory. We refer to a number of papers on local and global
existence, uniqueness, smoothing, and various asymptotic
(scattering-like) properties; see as a guide \cite{Cai97, Cr90,
Cr92, UniqPropDispEq, Hos99, Ion04, Ion06, Lar06, Lev01, Miz06,
Tak06, Tao00}, where some papers on linear Schr\"odinger operators
that exhibit related aspects of smoothing and evolution are also
included. However, for our purposes of studying refined asymptotic
properties of \eqref{generalpde} (a {\em refined scattering
theory}), these strong results turn out to be  not quite enough.
Therefore, first of all,
  we will need to  develop  new  Hermitian spectral theory of linear rescaled odd-order
 operators, which will refine some necessary delicate asymptotic properties of PDEs
  under consideration.

  In particular, as a first unavoidable step, we will need to study in detail
 the long- and short-time (on blow-up micro-scales) behaviour of
 solutions of the simplest {\em linear dispersion equation} for
 $k=1$ (the LDE--3)
  \begin{equation}
   \label{lde3}
    \fbox{$
    {\bf Basic \,\,Model \,\,I:} \quad
   u_t=u_{xxx} \inB \re \times \re_+,
    $}
    \end{equation}
 which, even in this simple standard form,  will create a number of surprises that require
   difficult unusual ``spectral"
 techniques.
 As a result, we then will be able to study similarity solutions
 (``nonlinear eigenfunctions") of the semilinear dispersion
 equation
  \be
   \label{k1N}
    \fbox{$
    {\bf Basic\,\, Model \,\,II:} \quad
   u_t=u_{xxx}-|u|^{p-1}u \inB \re \times \re_+ \whereA p>1.
    $}
    \ee
These two PDEs  \ef{lde3} and  \ef{k1N} are our main simplest
canonical models to study.

%% see
%%%references below.
%%% where $u=u(x,t)$ is the wave surface, with $x$
%%and $t$ denoting position and time.

%and the modified KdV (mKdV) equation
%\begin{equation}
%\label{mKdV} u_t = u_{xxx} + u^2u_x.
%\end{equation}
%Where \eqref{mKdV} is connected to the KdV equation \eqref{KdV} by the {\em Miura transformation}
%\begin{equation}
%\upsilon = u^2 + \sqrt{-3}u_x.
%\end{equation}
%In both, we have first-order perturbations while in \eqref{general
%pde} we take simpler zero-order ones.

%%%%%%%%%%%%%%%%%%%%%%%%%%%%%%%%%%%%%%%%%%%%%%%%%%%%%%%%%%%%
\subsection{First similarity solutions: earlier history
and Prandtl--Blasius results for odd-order boundary layer problem}

For both
 semilinear (\ref{par1}) and linear (\ref{lde3}) odd-order equations, we will
 apply the idea of self-similar solutions.
 In general,
 there are several techniques of
solving PDEs, including using travelling wave and other
group-invariant solutions.  As customary, similarity solutions are
a common way of attempting to understand and sometimes solve some
problems. Similarity solutions are widely used in  PDE theory,
since they simplify the problem by transforming PDEs into much
simpler ODEs, and hence this reduces the number of independent
variables. Rescaled variables can be used in studying PDEs, since
the coordinate system in which the problem is posed does not
affect the formulation of any fundamental physical laws involved.

 The first
general type of explicit solutions was traveling waves in
eighteen's century {\em d'Alembert's formula} for the {\em linear
wave equation}. The method of separation of variables was
developed by Fourier in the study of
heat conduction problems, and was later generalized and  %%%essentially
extended  by  Sturm and Liouville in the 1830s.
 Many famous mathematicians,
such as Euler, Poisson,  Lagrange, Liouville, Sturm, Laplace,
Darboux, B\"acklund, Lie, Jacobi, Boussinesq,  Goursat, and others
%D.F. Egorov
%L. Euler, Lagrange, J. Liouville, Laplace, G. Darboux, S. Lie,
%C.G.L. Jacobi, E. Goursat,
%D.F. Egorov
developed various techniques for obtaining explicit solutions of a
variety of linear and nonlinear models from physics and mechanics.
 Their methods  included a number of  particular transformations,
symmetries, expansions, separation of variables, {etc.} Similarity
solutions appeared in the works by Weierstrass around 1870, and by
Bolzman around 1890.  General principles for finding  solutions of
systems of ODEs and PDEs by symmetry reductions date back to the
famous Lie papers
%% \cite{L81}--\cite{Lie}
 published in the 1880s and 1890s.

More complicated and of key importance
 similarity solutions next appeared in {\em
Prandtl's equation}, occurring in {\em Prandtl's boundary layer
theory}, which was proposed in 1904, \cite{Prandtl}. After that,
 similarity solutions of
  linear and nonlinear  boundary-value problems became more common in the
 literature.
 Namely,
 this
 similarity solution in $\re^2$ was due to Blasius (1908)
\cite{Blas} for the equation of  incompressible fluids,
 \be
 \label{Pran.1}
 \psi_{yyy} + \psi_x \psi_{yy} - \psi_y \psi_{xy} - \psi_{yt} +
 u u_x + U_t=0,
 \ee
where $\psi=\psi(x,y,t)$ is the stream function and $U(x,t)$ is
the given external far-field (at $y=\infty$) velocity
distribution. Remarkably,
 this is an {\em odd-order} PDE.
 %% though different from those
 %%studied here.

%%% see further mathematical  details in
%%%Oleinik--Samokhin \cite{OlSam}.

%%%%%%%%%%%%%%%%%%%%%%%%%%%%%%%%%%%%%%%%%%%%%%%%%%%%%%%%%%%%%%
\subsection{On even-order models: the canonical heat equation}

 For a long period, until 1980s, whilst many asymptotic and singularity
 aspects of higher odd-order models have not been studied in great
detail, for a number of years, various even-order, both lower and
higher order, models have been reasonably well understood.  The
basic linear and most classical case of such a PDE is the {\em
heat equation}.  It is one of the most important linear partial
differential equations. A natural logic of the first part of our
research of the LDE--3 \ef{lde3} can be conveniently connected
with this classic area.

%%%%It is curious that for us it is convenient to begin with it.
%% and has only recently come to be fairly
%%%well understood.

Thus, the one-dimensional {heat equation} is the canonical PDE
\begin{equation}
\label{HE11}
 u_t =u_{xx} \quad \text{in} \quad
\mathbb{R}\times\mathbb{R}_+.
\end{equation}
%where $\Delta$ denotes the Laplacian, $\Delta \equiv
%\frac{\partial^2}{\partial x^2}$.
It has the classic fundamental solution
\begin{equation}
\label{Fund11}
 b(x,t) = \mbox{$\frac{1}{\sqrt{t}}\,F(y), \quad
\mbox{with the similarity variable} \quad  y=\frac{x}{\sqrt{t}}$},
\end{equation}
which takes Dirac's delta as an initial function, i.e., in the
sense of distributions,
%%\begin{equation*}
$b(x,0)=\delta(x)$.
%%\end{equation*}
Substituting \eqref{Fund11} into the heat equation \eqref{HE11}
yields a simple ODE for the rescaled kernel $F$:
 %%which reduces the PDE to ODE
\begin{equation}
\label{Herm1}
 \mbox{${\bf B}F \equiv F^{\prime\prime}
+\frac{1}{2}\,(Fy)' = 0$ \,\,in\,\, $\re$, \quad $\int F=1$}.
\end{equation}
This is  solved explicitly to show that the solution, $F$, is the
Gaussian:
\begin{equation}
\label{heatgauss}
  \mbox{$F(y) =
\frac{1}{2\sqrt{\pi}}\,\mathrm{e}^{-\frac{y^2}{4}}$},
\end{equation}
which  is strictly positive and does not have oscillatory
components to be traced out for linear dispersion equations. Full
spectral theory has been developed, with eigenvalues of the linear
second-order operator ${\bf B}$  and eigenfunctions given by (see
Birman--Solomjak \cite[p.~48]{BS} for full details and spectral
theory in $\ren$)
\begin{equation}
 \label{Herm1S}
\mbox{$ \s({\bf B})= \big\{\l_l=-\frac l2, \, l=0,1,2,...\big\}
\quad \mbox{and} \quad
 \psi_l(y) = \frac{(-1)^l}{\sqrt{l!}}\,D^l_y
F(y) \equiv H_l(y) \, F(y)$}.
\end{equation}
Here $H_l(y)$ denote the standard orthonormal  {\em Hermite
polynomials} (introduced in full generality by Hermite in the
1870s), which are induced by the Gaussian \eqref{heatgauss}. These
Hermite polynomials in 1D were earlier derived by {\em C. Sturm}
in 1836 \cite{St}, where these were used for a classification of
all types of multiple spatial zeros for solutions $u(x,t)$ of the
1D heat equation. This led him to the now famous {\em Sturm
Theorems on zero set}; see a full history in
\cite[Ch.~1]{GalGeom}. The Sturmian results found for the linear
heat equation, have been the basis for other higher and nonlinear
even-order and other related models, where the zero set analysis
of solutions become indeed more involved; see \cite{2mSturm} and
references therein. Whilst solutions to these models retain the
basic symmetric Gaussian structure, they have some oscillatory
tails instead of a pure simple exponential decay, which
essentially change the meaning of Sturmian theorems.

 Thus, the discrete spectrum and Hermite polynomials in \eqref{Herm1S}
are the cornerstone of the  {\em Hermitian spectral theory} of the
rescaled operator $\BB$ in \eqref{Herm1}. Since $\BB$ admits a
symmetric representation:
 $$
 \mbox{$
  \BB \equiv \frac 1 \rho \frac{\mathrm d}{{\mathrm d}y}\big(\rho \frac{\mathrm d}{{\mathrm d}y}
  \big),
  \quad \mbox{where} \quad \rho(y)= {\mathrm e}^{\frac{y^2}4},
  $}
  $$
  this theory is essentially self-adjoint in the weighted space $L^2_\rho(\re)$ including the $N$-dimensional
  case \cite{BS}.
 As an unavoidable step,
we will need to develop a similar theory for the LDE--3 \ef{lde3},
which is not self-adjoint and will go in different lines than the
classic one.

%%%%%%%%%%%%%%%%%%%%%%%%%%%%%%%%%%%%%%%%%%%%%%
\subsection{Very singular similarity solutions}

Around the beginning of the 1980s, study of asymptotics of the {\em
semilinear heat equation with absorption}, given by
\begin{equation}
\label{heat}
   u_t=\Delta u - u^p \quad \text{in} \quad
\mathbb{R}^N\times\mathbb{R}, \quad p>1 \quad (u \ge 0),
\end{equation}
led to a new class of similarity solutions called {\em very
singular solutions} (VSSs). These are self-similar solutions,
which are structurally stable in the evolution and hence, as $t
\to +\infty$, attract wide classes of other more general
solutions. In addition, as $t\to 0$, such a VSS concentrates at
$x=0$ with infinite initial mass and this justifies the term
``very singular".  We refer to the books \cite{StabTechPdeVSS,
quasilin} for extra details and history concerning VSSs. Whilst
VSSs were being studied early on, it wasn't until about 1985 in
papers by Kamin and Peletier \cite{VssHeatAbs, SSHeatAbs} in which
the term VSS was actually used.

The VSSs for the semilinear heat equation \eqref{heat} are given by
\begin{equation*}
u_\ast(x,t) = t^{-\frac{1}{p-1}}f(y), \quad
\mbox{$y=\frac{x}{\sqrt{t}}$},  \quad \mbox{where $f>0$ solves the
elliptic equation}
\end{equation*}
%%%where $f$ solves the elliptic equation
\begin{equation*}
\begin{cases}
\Delta f +\mbox{$\frac{1}{2}\, y\cdot\nabla f +\frac{1}{p-1}\,
f-f^p=0$} \quad \text{in} \quad \mathbb{R}^N,\\ \text{$f(y)$ has
exponential decay as $y\to \infty$}.
\end{cases}
\end{equation*}
Existence of the VSS was established in \cite{AsymEigNonlinPara}
%%by Galaktionov, Kurdyumov,
%%and Samarski\u{\i}
 using a PDE approach,
 and  almost simultaneously, via an ODE
approach,
%% by Brezis, Peletier, and Terman,
in \cite{VssHeatAbs}. Uniqueness of VSS was first proved later
%% by
%%Kamin and V\'{e}ron
 in \cite{ExisUniqVSSPorMed}, by using a PDE
comparison method.  We refer to  \cite{GalANS, VSSSParaPDEs} for
more references and history of VSSs for semilinear parabolic
equation.

%%%%%%%%%%%%%%%%%%%%%%%%%%%%%%%%%%%%%%%%%%%%%%%%%%%
\subsection{Back again to the KdV-type equations}

It is curious that, formally, VSSs can be prescribed for the
classic KdV equation \eqref{KdV}. These have a standard
self-similar form
\begin{equation}
 \label{VSSK}
u_\ast(x,t) = t^{-\frac{2}{3}}f(y), \quad y = x/t^{\frac{1}{3}}.
\end{equation}
Hence \eqref{KdV} reduces to the following ODE:
\begin{equation}
\label{VSSKdV} \mbox{$f^{\prime\prime\prime}+\frac{1}{3}\,
yf^\prime+\frac{2}{3}\, f+ff^\prime$} = 0 \quad \text{in} \quad
\mathbb{R}.
\end{equation}
However, whilst we can seemingly construct a VSS for the KdV
equation, we know that solutions of this type do not exist. The
only natural solution we arrive at is the trivial one, $f\equiv
0$. This nonexistence conclusion directly follows from the global
existence results by Kato \cite{OnKdV, QuasiEvoPDE} and Strauss
\cite{DispEnerWav}. Indeed,  \eqref{KdV} is invariant under
reflection with
\begin{equation}
\label{ReflKdvVSS} t\mapsto T-t \quad \text{and} \quad x\mapsto -x,
\end{equation}
where very singular solutions \eqref{VSSK} are then given by
\begin{equation*}
u_\ast(x,t) = (T-t)^{-\frac{2}{3}}f(y), \quad y =
x/(T-t)^{\frac{1}{3}}.
\end{equation*}
Therefore existence of a nontrivial VSS would mean finite time
blow-up of solutions, since
\begin{equation*}
 \mbox{$
\sup_x|u(x,t)|\sim(T-t)^{-\frac{2}{3}} \to +\infty \quad \text{as}
\quad t\to T^-.
 $}
\end{equation*}
In view of well-known global existence results for the KdV equation
(see references above),  this implies the non-existence of very
singular solutions for \eqref{KdV}.

For the {\em modified KdV equation}
\begin{equation}
\label{GkdV}
 u_t = u_{xxx}+u^pu_x,
\end{equation}
we can find VSS for $p\geq4$, but not for $p<4$. In the case
$p\geq4$, \eqref{GkdV} has self-similar solutions of the standard
form
\begin{equation*}
u_\ast(x,t) = t^{-\frac{2}{3p}}f(y), \quad y = x/t^{\frac{1}{3}},
\end{equation*}
where the rescaled kernel satisfies
\begin{equation}
\label{gkdvODE} \mbox{$f^{\prime\prime\prime}+\frac{1}{3}\,
yf^\prime+\frac{2}{3p}\,f+f^pf^\prime$} = 0.
\end{equation}
The case $p=4$ has been studied more extensively; see papers
\cite{BilinEst,SubcritGKdV} as just two examples.
%%% since here VSSs
%%obeys the mass conservation, so that the ODE \eqref{gkdvODE} can
%%%be integrated once, thus establishing existence and uniqueness of
 %%%%the VSS.
  We thus need to study such similarity solutions for our second
  semilinear model \ef{k1N}.

%%%%%%%%%%%%%%%%%%%%%%%%%%%%%%%%%%%%%%%%%%%%%%%%%%%%%%%%%
\subsection{Four main linear and nonlinear odd-order models to study}

We will treat, in particular, three generalized odd-order models,
which, from our ``spectral-like" and refined  asymptotic points of
view, essentially have never been looked at before. In particular,
we look for similarity solutions and using asymptotic, analytic,
and some numerical methods, we will attempt to find and justify
some local and global properties of the rescaled solutions.  We
use a number of techniques, applied previously for parabolic
even-order problems, for the odd-order ones. However, there are
difficulties that arise, in particular due to the highly
oscillatory nature of fundamental solutions and hence related
linear and nonlinear eigenfunctions.
%% and the anisotrophy of the
%%%operator.
%% as well as non-symmetry about
%%$y=0$, in the rescaled solutions.

Our {\sc first aim} is
%%, using the fundamental solutions,
to develop a  {\em Hermitian spectral theory} of the associated
rescaled operators, for {\em linear dispersion} odd-order
equations such as (cf. \ef{lde3} for $k=1$)
\begin{equation}
\label{airyfunc}
 {\bf (I):}\quad
 u_t = (-1)^{k+1}D_x^{2k+1}u \quad \text{in} \quad
\mathbb{R}\times\mathbb{R}_+ \quad \mbox{for any} \,\,\, k \ge 1.
\end{equation}
In doing this, we gain the {\sc second goal}, an understanding of
its behaviour, which we can then use in the corresponding {\em
semilinear dispersion equation} and its VSSs
\begin{equation}
\label{semilinear pde}
 {\bf (II)}: \quad
 u_t = (-1)^{k+1}D_x^{2k+1}u-|u|^{p-1}u \inB \re \times \re_+,
\quad \mbox{where} \quad p>1.
\end{equation}

In a forthcoming paper \cite{FerGalII},
%% The {\sc second
 %%goal}, to be achieved in \cite{FerGalII},
 these concepts and ideas will be
applied and extended  to PDEs with {\em nonlinear dispersion} (the
NDEs) of the form
\begin{equation}
 \label{NDE1}
  {\bf (III)}: \quad
u_t = (-1)^{k+1}D_x^{2k+1}(|u|^nu) \inB \re \times \re_+ \whereA
n>0.
\end{equation}
 %% is a fixed parameter.
   Concerning self-similar solutions of
\eqref{NDE1}, called ``nonlinear eigenfunctions", asymptotic
behaviour and general properties of solutions such as existence,
uniqueness, shock waves, entropy approaches, etc.,
 turned out to create a number of very difficult questions.
 %%very little is known in mathematical literature.
   First steps of
the study of shock and rarefaction waves for NDEs including the
NDE--3 such as
 $$
 u_t=(u u_x)_{xx} \inB \re \times \re_+,
 $$
 are performed in \cite{GPndeII, GPnde}. As a natural extension of
the NDE \eqref{NDE1}, we also
%%% very briefly
discuss the VSS for the full NDE:
%%%with absorption:
\begin{equation}
 \label{ndekp}
  {\bf (IV)}: \quad u_t =  (-1)^{k+1}D_x^{2k+1}(|u|^n u) -|u|^{p-1}u \inB \re \times
  \re_+
\whereA  p>n+1.
\end{equation}
In both cases of NDEs \ef{NDE1} and \ef{ndekp}, the key idea is to
perform a ``homotopy" limit $n \to 0^+$ to arrive at the rescaled
linear operators and next to use Hermitian spectral theory for
bifurcation-branching analysis of nonlinear eigenfunctions and
VSSs.

%%%%%%%%%%%%%%%%%%%%%%%%%%%%%%%%%%%%%%%%%
%%%%%%%%%%%%%%%%%%%%%%%%%%%%%%%%%%%%%%%%%%%%%%%%%%%%%%%%%%%%%%%
\section{Basic linear dispersion models: LDEs
 and their fundamental solutions}
%% and hermitian spectral theory}
\label{LinChap}

%%%%%%%%%%%%%%%%%%%%%%%%%%%%
 \subsection{Basic LDE: fundamental solutions and rescaled kernels}

Before looking at some more complicated nonlinear PDEs, it is
important to understand in greater details  how the solutions of
linear PDEs behave for large and small times. This spectral theory
of  higher-odd-order linear rescaled operators will be crucial in
the understanding of related nonlinear ones. In particular, as
customary, spectral theory formed in the linear case will play a
large role and will be used in developing understanding of
bifurcations, branching, and asymptotic behaviour for nonlinear
equations.

Thus, we consider the corresponding {\em linear dispersion
equation} \ef{airyfunc}, which we al so call the LDE--$(2k+1)$.
 %%\begin{equation}
%%\label{airyfunc} u_t = (-1)^{k+1}D_x^{2k+1}u \quad \text{in} \quad
%%\mathbb{R}\times\mathbb{R}_+, \quad k \ge 1.
%%\end{equation}
Whilst some lower-order cases for the odd-order linear PDE
\eqref{airyfunc} such as the LDE--3 \ef{lde3} are generally better
understood, the higher-order cases for $k \gg 1$ are not.  Indeed,
it is well-known that the fundamental solution for the lowest
order case $k=1$, i.e., for \ef{lde3}, will lead to the classic
{\em Airy function}. However, we will need  a more general
calculus that is applied for arbitrary $k \ge 1$.

%%%%%%%%%%%%%%%%%%%%%%%%%%%%%%%%%%%%%%%%%%%%%
 %%%\subsection{Fundamental solutions and rescaled kernels}

\ssk

Consider the self-similar fundamental solution of \eqref{airyfunc}
of the standard similarity form
\begin{equation}
\label{kern} b(x,t) = t^{-\frac{1}{2k+1}}F(y), \quad y =
x/t^{\frac{1}{2k+1}}, \quad \mbox{so that} \quad b(x,0)=\d(x),
\end{equation}
in the sense of distributions. Substituting $b(x,t)$ into the PDE
\eqref{airyfunc}, we obtain the ODE for the rescaled kernel $F$,
\begin{equation}
\label{linODE} \mbox{${\bf B}F \equiv (-1)^{k+1}D_y^{2k+1}F
+\frac{1}{2k+1}\,yD_yF + \frac{1}{2k+1}\,F = 0$ \,\,\, in \,\,\,
$\re$, \quad $\int F=1$},
\end{equation}
where  ${\bf B}$ denotes the first key linear rescaled operator
for the LDE. Note that it is possible to integrate the ODE
\eqref{linODE} once, to find that $F(y)$ solves
\begin{equation}
\label{lineqn}
 \mbox{$(-1)^{k+1}F^{(2k)}+ \frac{1}{2k+1}\,y F  =
0$} \quad \text{for} \quad y\in\mathbb{R},
\end{equation}
which is now a linear ODE of order $2k$.
%% From the original linear
%%PDE of order $2k+1$, the problem has now been reduced to a linear
%%ODE of order $2k$.
 Further reductions of this ODE are not possible.
 %%even for $k=1$ .
%The fundamental solution \eqref{kern} has
%initial data $b(x,0) = \delta (x)$, where $\delta$ is the Dirac
%delta function.

For $k=1$, the solution of (\ref{lineqn}) is the classic {\em Airy
function}:
\begin{equation}
 \label{AiNew}
 \mbox{$
 F'' + \frac 13\, y F=0 \LongA
F(y)=\mathrm{Ai}(y).
 $}
\end{equation}
Note that in our case,  due to a difference in sign in the ODE, we
actually have $\mathrm{Ai}(-y)$, but we shall refer to it as the
{Airy} one.

%The Airy equation is given by
%\begin{equation*}
%F^{\prime\prime} - Fy = 0.
%\end{equation*}
%This has linearly independent solutions, $Ai(y)$ and $Bi(y)$ given
%by

%This is the simplest second-order linear differential equation with
%a turning point (a point where the character of the solutions
%changes from oscillatory to exponential).

%The Airy function describes the appearance of a star  a point
%source of light  as it appears in a telescope. The ideal point
%image becomes a series of concentric ripples because of the limited
%aperture and the wave nature of light (Suiter 1994). It is also the
%solution to Schrodinger's equation for a particle confined within a
%triangular potential well.

\subsection{Asymptotic expansions of the rescaled fundamental kernel}
\label{asymptotics}

An important  and standard technique in trying to find the
behaviour of solutions of ODEs is to use asymptotic analysis. This
gives the limiting behaviour of solutions, in particular as
$y\to\pm\infty$. We refer to the book by Bender and Orszag
\cite{AsympPert}  for various asymptotic techniques.  We use here
a method of determining the asymptotic behaviour of the linear
ODE, which corresponds to the classic WKBJ multi-scale analysis of
ODEs, whose basic ideas go back to the 1920s and were reflected in
a number of well-known monographs. Actually, in ODE theory,
asymptotics for ODEs such as \eqref{lineqn} are well-known and
have been classified. However, we will need some more refined
formulae for further applications. Some of them for arbitrary $k
\ge 1$ are not available in standard literature. Looking at the
asymptotic behaviour of the solution of the integrated ODE
\eqref{lineqn}, we now write it as
\begin{equation}
\label{intODE} \mbox{$(-1)^{k+1}F^{(2k)}  = -\frac{1}{2k+1}\,Fy$}.
\end{equation}

%%%%%%%%%%%%%%%%%%%%%%%%%%%%%%%%%%%%%%%%%%%%%%%%%%%
\noi\underline{\em Algebraic decay as $y \to + \iy$}. Let us look
first at the rescaled solution $F(y)$  with {\em oscillatory
algebraic} decay  as $y \rightarrow +\infty$.
 %% It should be noted that all the expansions can be rigorously
 %%%justified, but we do not put special efforts on this here.
 As customary, we set, for  convenience,
\begin{equation}
\label{asympf}
  F(y)=\mathrm{e}^{s(y)} \quad \text{as} \quad
y\to\infty.
\end{equation}
Now assume, as a first approximation, that $s(y)$ is some
polynomial, such that $s(y)\sim a y^{b}$. Then the kernel $F(y)$ and
its derivatives may be given by
\begin{equation*}
 \begin{matrix}
%%\begin{split}
F^{\prime} = s^{\prime}\mathrm{e}^s,\,\,\, F^{\prime\prime}
=\big[s^{\prime\prime}+(s^{\prime})^2\big]\mathrm{e}^s,\,\,\,
F^{\prime\prime\prime}
=\big[s^{\prime\prime\prime}+3s^{\prime}s^{\prime\prime}+(s^{\prime})^3\big]\mathrm{e}^s,\,
...\, , \smallskip \smallskip
\\
%%%&\vdots\\
F^{(2k)}=\big[s^{(2k)}+\hdots+k(2k-1)(s^{\prime})^{2k-2}s^{\prime\prime}+(s^{\prime})^{2k}\big]\mathrm{e}^s,
 \end{matrix}
%%\end{split}
\end{equation*}
where $s^{\prime}\sim ab y^{b -1}, \, s^{\prime\prime}\sim ab(b -1)\
y^{b -2}, \, \hdots, \, s^{(2k)}\sim \frac{ab !}{(b -2k)!}y^{b
-2k}$.
 Substituting this into the ODE \eqref{intODE}, it can easily be seen that all the
$\mathrm{e}^{s(y)}$ terms cancel.  From the resulting equation, we
look to do a dominant balance analysis, in order to determine the
leading order behaviour of $F(y)$. In doing so, we find a first
approximation for the function $s(y)$.  The balance of the equation
depends on the value of the parameter $b$ and we obtain two
different cases.

If $b \leq 0$, then as $y\rightarrow +\infty$, every term on the
left-hand side of the ODE is $o(y)$. Therefore there is no balance
in this case.
 If $b >0$, then $y^{b -\tilde{n}} =o(y^{2k(b -1)})$, for any
$\tilde{n}>1$.  Therefore balancing leading terms, we have that
\begin{equation*}
\mbox{$(-1)^{k+1}a^{2k}b^{2k}y^{2k(b -1)} \sim -\frac{1}{2k+1}\,y$},
\end{equation*}
for all $k\in\mathbb{Z}_+$.  By first equating powers of $y$ and
then coefficients, we find our parameters
\begin{equation*}
\mbox{$b = \frac{2k+1}{2k}$} \quad \mbox{and} \quad \mbox{$a =
-\mathrm{i}\,(2k+1)^{-\frac{1}{2k}} \big(\frac{2k}{2k+1}\big)$}.
\end{equation*}
%%and
%%\begin{equation*}
%%\mbox{$a = -\mathrm{i}\,(2k+1)^{-\frac{1}{2k}}
%%\big(\frac{2k}{2k+1}\big)$}.
%%\end{equation*}
Hence, we now have the first approximation to $s(y)$, with
\begin{equation*}
\begin{matrix}
s(y)=\mbox{$-2k\mathrm{i}\big(\frac{y}{2k+1}\big)^{\frac{2k+1}{2k}}
+ c(y)$},\,\,\,
s^{\prime}(y)=\mbox{$-\mathrm{i}\big(\frac{y}{2k+1}\big)^{\frac{1}{2k}}
+ c^{\prime}(y)$},\ssk\\
s^{\prime\prime}(y)=\mbox{$-\frac{\mathrm{i}}{2k(2k+1)^{\frac{1}{2k}}}\,y^{-\frac{2k-1}{2k}}
+ c^{\prime\prime}(y)$}, \, ...\, .
%%%\\ &\hdots
\end{matrix}
\end{equation*}
Here $c(y)\sim o(y^{\frac{2k+1}{2k}})$ is some function of $y$ and
the next term in the approximation of $s(y)$.  Since $k\geq1$, we
must have that $c^{(m+1)}(y)\sim o(c^{(m)}(y))$ for any $m>1$, which
will be used in determining leading order terms.

We attempt to find this function $c(y)$, in order to improve the
approximation of $s(y)$.  We let
$g(y)=-2k\mathrm{i}\big(\frac{y}{2k+1}\big)^{\frac{2k+1}{2k}}$ for
convenience, to see how the terms are balanced.  Hence balancing
the leading order terms yields
%%\begin{equation*}
%\mbox{$2k(g^{\prime}(y))^{2k-1}c^{\prime}(y) \sim
%-\frac{2k}{2}(2k-1)(g^{\prime}(y))^{2k-2}g^{\prime\prime}(y)$}
%\,\Longrightarrow\,
$c'(y) \sim
\mbox{$-\frac{2k-1}2\,\frac{g^{\prime\prime}(y)}{g^{\prime}(y)}$}=
\mbox{$-\frac{2k-1}{4ky}$}$.
%%%%\end{equation*}
%%so
%%\begin{equation*}
%%c'(y) \sim
%%\mbox{$-\frac{1}{2}\,(2k-1)\,\frac{g^{\prime\prime}(y)}{g^{\prime}(y)}$}=
%%\mbox{$-\frac{2k-1}{4ky}$}.
%%\end{equation*}
Integrating this, the next term in the expansion can be found to
be
 $%\begin{equation*}
 \mbox{$c(y) \sim -\frac{2k-1}{4k}\ln{y}$}.
  $ %%\end{equation*}
The second term in this expansion, $c(y)$, is called the {\em
controlling factor}.  It is an important term, as we will see that
it governs the decay (or any possible growth) of solutions as
$y\to+\infty$.

Whilst with the first two terms, we can have a clear idea of the
leading order behaviour, we look to find a better approximation and
now expand once again with
\begin{equation*}
s(y)=\mbox{$-2k\mathrm{i}\big(\frac{y}{2k+1}\big)^{\frac{2k+1}{2k}}
-\frac{2k-1}{4k}\ln{y} + d(y)$},
\end{equation*}
where $d(y)=o(c(y))$.  Balancing leading order terms once again, we
find that
%%%\begin{equation*}
%\mbox{$2k(g '(y))^{2k-1}d '(y) \sim -\frac{2k}{2}(2k-1)2(g '
%(y))^{2k-3}\big(g''(y)\big)^2$} \, \Longrightarrow\,
%% \mbox{
$d'(y) \sim -\frac{2k-1}{4k^2y^2}$.
%%%}.
%%%%\end{equation*}
%%so
%%\begin{equation*}
%%\mbox{$d^{\prime}(y) \sim -\frac{2k-1}{4k^2y^2}$}.
%%%\end{equation*}
Then by integrating, we find that the third term in the expansion
is given by
   $ %%\begin{equation*}
\mbox{$d(y) \sim \frac{2k-1}{4k^2y}$}.
 $ %%\end{equation*}
Note that this is the third term in the expansion and hence all
following lower-order terms are $o(1)$, and so, by \eqref{asympf},
this does not affect the exponential.

From this, we find that the required asymptotic behaviour of
$F(y)$
 %% as $y\rightarrow +\infty$ can be
 is given by
\begin{equation*}
\mbox{$F(y) \sim
y^{-\frac{2k-1}{4k}}\exp\big\{-2k\mathrm{i}\big(\frac{y}{2k+1}\big)^{\frac{2k+1}{2k}}\big\}$}
\quad \text{as} \quad y \to +\infty.
\end{equation*}
Therefore, for real solutions, we have that, for some constant
$\hat c \in \re$,
\begin{equation}
\label{posasympt}
 \mbox{$F(y) \sim y^{-\frac{2k-1}{4k}}\cos
{\big(d_ky^{\frac{2k+1}{2k}}+\hat{c}\big)}$} \quad \text{as} \quad
y \to +\infty, \quad \mbox{where} \quad
\mbox{$d_k=2k\,\big(\frac{1}{2k+1}\big)^{\frac{2k+1}{2k}}$}.
\end{equation}
One can see that \ef{posasympt} guarantees the convergence of the
integral of $F$ in \ef{linODE} (but not in the absolute sense), so
that the fundamental solution $b(x,t)$ at $t=0$ takes Dirac's
delta in the sense of bounded measures. Note that, {\em a priori},
 this is not guaranteed in general, since  delicate ODE
 asymptotics are necessary to charge such a property.

%%where
%%\begin{equation}
%%\label{ExpConstd}
%%\mbox{$d_k=2k\,\big(\frac{1}{2k+1}\big)^{\frac{2k+1}{2k}}$}
%%%\end{equation}
%%and $\hat{c}$ is some constant.

\ssk

%%%%%%%%%%%%%%%%%%%%%%%%%%%%%%%%%%%%%%%%%%%%%%%%%%%%%%%
\noi\underline{\em Exponential decay as $y \to - \iy$}. The same
analysis can be applied for $y\to-\infty$, by letting $y\mapsto-y$
and performing the same calculations.  Hence we find that the
rescaled kernel decays exponentially fast in the opposite
direction,
\begin{equation}
 \label{ymin}
%%\begin{split}
F(y) \sim
%%&
|y|^{-\frac{2k-1}{4k}} \cos\big(d_k|y|^{\frac{2k+1}{2k}}\sin
{b_k}+\hat{c}\big)
%%% \\[1mm]
%%& \times
\, \exp\big\{-\hat d_k |y|^{\frac{2k+1}{2k}}\big\} \quad \text{as}
\quad y \to -\infty,
%%%\end{split}
\end{equation}
where $d_k$ is as in \ef{posasympt} and
\begin{equation}
 \label{yminNN}
\hat d_k= d_k |\cos b_k|>0, \quad b_k=
\begin{cases}
\mbox{$\frac{\pi}{k}\lfloor\frac{k+1}{2}\rfloor$} \quad \text{for
even $k$},\\
\mbox{$\frac{\pi}{k}\frac{k+1}{2} + \frac{\pi}{2k}$} \quad \text{for
odd $k$}.
\end{cases}
\end{equation}
We note that the coefficient $a$ has many roots,
 as follows from the algebraic equation:
 $$
  \mbox{$
  F(y) \sim {\mathrm e}^{a|y|^{(2k+1)/2k}}
   \LongA a^{2k}= \frac{(-1)^{k+1} (2k)^{2k}}{(2k+1)^{2k+1}},
    $}
    $$
 which represent
different solutions of the ODE.  However, we only want roots such
that there is exponential decay, rather than growth.  Therefore we
exclude the roots where there is growth, which corresponds to
$\mathrm{Re}\, a > 0$.  The asymptotics here show the behaviour of
the first roots, such that $\mathrm{Re}\, a \leq 0$. The method of
finding these roots and all other roots may be seen in
 Appendix A,
 %%Section
%%\ref{SectRadcond},
 which also explains more carefully which roots we
need to look at.
%It is noted that the principal roots have been
%taken, when calculating the behaviour.  The method of finding these
%roots and all other roots may be seen in Section \ref{SectRadcond}.

\ssk

Thus, according to \ef{posasympt}, as $y \to +\infty$, we have a
purely imaginary root  for all values of $k$ and this gives slow
decaying oscillatory behaviour. The decay of the oscillations also
increases for larger $k$, from $y^{-\frac 14}$ for $k=1$  towards
$\sim y^{-\frac 12}$ as $k \to \iy$. On the contrary, as $y \to
-\infty$, \ef{ymin} shows exponential decay with  oscillations for
any $k>1$ and these oscillations are always exponentially fast.
For the only case $k=1$, we have $b_k=\pi$, hence $\sin b_k = 0$.
This gives the pure exponential non-oscillatory behaviour of the
Airy function:
\begin{equation}
 \label{Airymin}
F(y)\sim |y|^{-\frac{1}{4}}\mathrm{e}^{-2 \cdot
3^{-\frac{3}{2}}|y|^{\frac{3}{2}}} \quad \text{as} \quad
y\to-\infty.
\end{equation}
This is the reason why the Airy function is the only odd-order
linear case  where there is no exponential {\em infinitely
oscillatory} behaviour as $y\to-\infty$. In a natural sense, this
mimics the positive behaviour of the Gaussian \ef{heatgauss}, but
only partially, in a one-sided limit $y \to - \iy$. Anyway, the
exponential decay in \ef{ymin} simplifies a part of spectral
analysis of the operator $\BB$, for which  we will need to pay the
main attention to the opposite limit $y \to + \iy$, where most of
technical difficulties appeared from.

\ssk

We now summarize the results of our asymptotic analysis to get a
sharp uniform bound of the rescaled kernel $F$, which will
eventually define the weighted functional space for $\BB$:

\begin{proposition}
The rescaled kernel $F(y)$ of the fundamental solution
\eqref{kern} of the LDE \eqref{airyfunc} satisfies
\begin{equation}
\label{initFest} |F(y)| \leq
\begin{cases}D_0(1+y^2)^{-\frac{2k-1}{8k}}\mathrm{e}^{- \hat d_k|y|^\alpha}
\quad \text{for} \quad y\leq 0,\\ D_0(1+y^2)^{-\frac{2k-1}{8k}}
\quad \text{for}\quad y\geq 0,
\end{cases}
\end{equation}
where $\hat d_k$ is as in $(\ref{yminNN})$, $D_0$ is a positive
constant dependent on $k$, and
\begin{equation*}
\alpha = \mbox{$\frac{2k+1}{2k}$} \in (1,2) \quad \text{for all}
\quad k\geq1.
\end{equation*}
\end{proposition}

%%This follows from the above asymptotic analysis.

%%%%%%%%%%%%%%%%%%%%%%%%%%%%%%%%%%%%%%%%%%%%%%%%%%%%%%%%%%%%%%%%%%%
\subsection{Numerical construction of fundamental kernels}
\label{NumConstr}

%%Numerics play an important part in the theory of differential
%%equations.  They are an important way of looking at the behaviour
%%of the solutions and to check any results found.
The  results here
%%in this section
 were obtained by using the {\tt
MatLab bvp4c} solver, to look at the singular solutions of the
linear ODE \eqref{linODE}.

\begin{figure}[htbp]
\centering \subfigure[Oscillations of $F(y)$ for $y \in (0,20)$.]{
\label{k3a}
\includegraphics[scale=0.45]{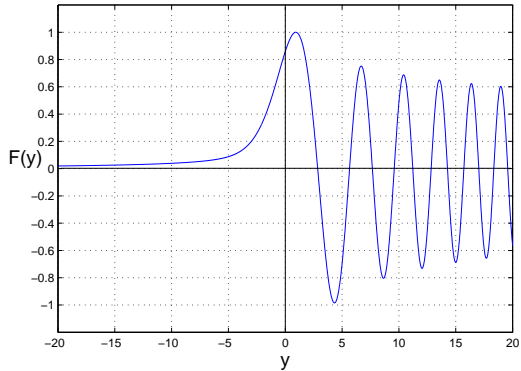}}
\subfigure[Oscillatory ``tail" of $F(y)$ for $y\gg 1$.]{
\label{k3b}
\includegraphics[scale=0.45]{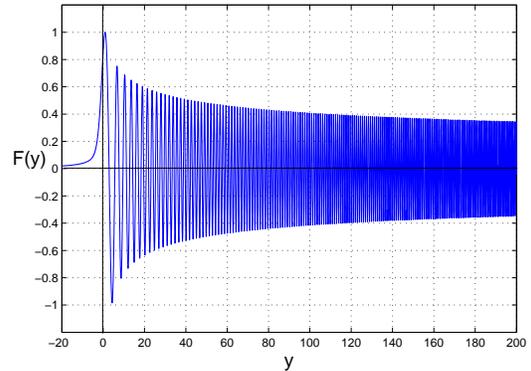}}
\caption{\small The rescaled kernel $F(y)=\mathrm{Ai}(y)$ of the
fundamental solution \eqref{kern} of (\ref{lde3})  ($k=1$)
obtained by two-sided shooting.} \label{k3plot}
\end{figure}

Obviously, if $F$ is a solution of \eqref{lineqn}, then $cF$ is
also a solution, for all $c \in \mathbb{R}$.  Due to such
``non-uniqueness" (and, to some extent, an ``instability") of the
solutions, the zero solution is likely to be found using numerical
methods.  Therefore in order to ensure that a non-zero solution is
obtained, we set the (normalisation) constraint
\begin{equation}
\label{constraint} \max |F| = 1,
\end{equation}
which is attained at some point $y=\hat{a}$. We then solved the
ODE for different right and left solutions at this maximum point
$\hat{a}$, using the {\tt bvp4c} solver.  Since we have fixed
$\hat{a}$ as a maximum, the first derivative is also zero at this
point. The other boundary condition placed was to ensure that
$F(y)=0$, at an end point that is sufficiently removed from
$\hat{a}$. In order for a match of right and left solutions at the
point $\hat{a}$, we needed to have the second derivative of the
solution to be the same here, in order for $F(y)$ to be continuous
at $y=\hat{a}$. So the value of $\hat{a}$ was moved in order to
match the second derivative for the right and left solutions. See
Figure \ref{k3plot}, which was obtained by this shooting method.

\begin{figure}[htbp]
\begin{center}
\includegraphics[scale=0.5]{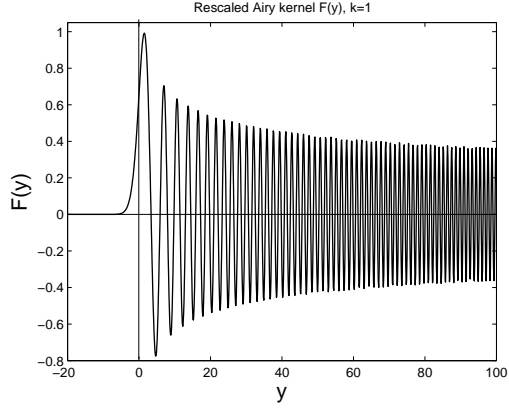}
\caption{\small The rescaled kernel $F(y) =\mathrm{Ai}(y)$ of the
ODE \eqref{AiNew}  ($k=1$)
%%%. This is denoted as $\mathrm{Ai}_5(y)$
obtained by one-sided shooting via the bundle (\ref{Airymin}).}
 \label{F1side}
\end{center}
\end{figure}

Indeed, for $k=1$, as a simple alternative, one can use ``shooting
from the left" in the second-order ODE \ef{AiNew} by using the
exponential decay asymptotics \ef{Airymin} for $y \ll -1$ to get a
proper behaviour \ef{posasympt} for $y \gg 1$, since no other are
available in the version \ef{AiNew}.
 Such a simpler one-sided shooting is
 shown in Figure \ref{F1side}, which is rather similar to Figure
 \ef{k3plot}(b). What is most important is that the sizes of the
 tails in both Figures coincide. However, such an approach is
 hardly applied in higher-order cases $k \ge 2$, where extra
 solutions selecting procedures such as multi-dimensional shooting
 are unavoidable.

\begin{figure}[htbp]
\begin{center}
\includegraphics[scale=0.5]{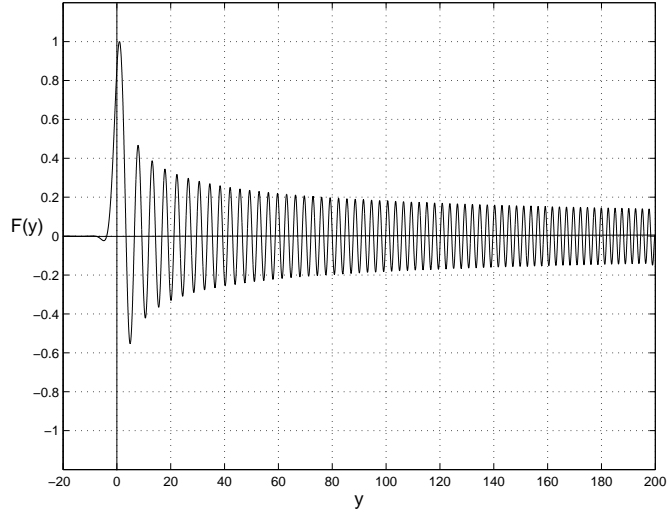}
\caption{\small The rescaled kernel $F(y) =\mathrm{Ai}_5(y)$ of
the fundamental solution \eqref{kern} of the LDE--5 $u_t =
-u_{xxxxx}$ ($k=2$)
%%%. This is denoted as $\mathrm{Ai}_5(y)$
obtained by two-sided shooting.}
 \label{k5a}
\end{center}
\end{figure}

A similar method of two-sided shooting was applied to the fifth
order ($k=2$) equation, where values of the left-hand solution
were used as boundary conditions for the right solution, and the
correct value for $\hat{a}$ was found by matching the fourth
derivative. See Figure \ref{k5a}.
%%    Obviously due to the method used, we cannot guarantee that the exact
%%solutions were found, but these show the behaviour of the solutions.
Note that, in view of \ef{constraint}, the plots do not show the
fundamental solutions of the ODE, such that $\int F =1$, but
rescaled profiles, for which the normalization  \eqref{constraint}
is satisfied.

As one can see from comparing Figures \ref{k3plot} and \ref{k5a},
the oscillations are smaller as $k$ increases.  This follows from
the asymptotic analysis done in Section \ref{asymptotics}; see
\ef{posasympt}. Decay of the algebraic envelope is very slow as $y
\to +\infty$ also, especially for the {\em Airy} function when
$k=1$, and is faster for $k=2$ in Figure \ref{k5a}.

For convenience we denote these rescaled kernels as higher-order
{\em Airy} functions by
\begin{equation}
\label{Ai}
   F(y) = {\rm Ai}_{2k+1}(y) \quad \text{for} \quad
k=1,2,\hdots, \quad \mbox{so that ${\rm Ai}={\rm Ai}_3$.}
\end{equation}

%\newpage

%%%%%%%%%%%%%%%%%%%%%%%%%%%%%%%%%%%%%%%%%%%%%%%%%%%%%%%%5
%%%%%%%%%%%%%%%%%%%%%%%%%%%%%%%%%%%%%%%%%%%%%%%%%%%%%

%%%%%%%%%%%%%%%%%%%%%%%%%%%%%%%%%%%%%%%%%%%%%%%%%%%%%%%%%%%%%%%
 %%%%%%%%%%%%%%%%%%%%%%%%%%%%%%%%%%%%%%%%%%%%%%%%%%%%%%%%%%%%%%5

%\newpage

%%%%%%%%%%%%%%%%%%%%%%%%%%%%%%%%%%%%%%%%%%%%%%%%%%%%%%%%%%%%%%%
 %%%%%%%%%%%%%%%%%%%%%%%%%%%%%%%%%%%%%%%%%%%%%%%%%%%%%%%%%%%%%%5

%%%%%%%%%%%%%%%%%%%%%%%%%%%%%%%%%%%%%%%%%
%%%%%%%%%%%%%%%%%%%%%%%%%%%%%%%%%%%%%%%%%
%%%%%%%%%%%%%%%%%%%%%%%%%%%%%%%%%%%%%%%%%
%%%%%%%%%%%%%%%%%%%%%%%%%%%%%%%%%%%%%%%%%

%%%%%%%%%%%%%%%%%%%%%%%%%%%%%%%%%%%%%%%%%%%%%%%%%%%%%%%%%
\section{Explicit semigroup representation: first spectral properties}

Here, we begin to develop first aspects of a Hermitian spectral
theory, which our further analysis will depend crucially on.

%%%%%%%%%%%%%%%%%%%%%%%%%%%%%%%%%%%%%%%%%%%%%%%%%%%%%%
\subsection{Operator {\bf B}: formal series}
\label{SemiOpB}

Let $u(x,t)$ be the solution of the Cauchy problem for the
LDE--$(2k+1)$ \eqref{airyfunc}, with bounded measurable initial
data $u(x,0)= u_0(x)$ from a weighted $L^2$ space to be specified.
The solution is then represented by the convolution of initial
data with the fundamental solution given by a Poisson-type
integral:
\begin{equation}
\label{convolution}
  \mbox{$
u(x,t) = b(t)\ast u_0 \equiv
t^{-\frac{1}{2k+1}}\int\limits_{\mathbb{R}}
\mbox{$F\big((x-z)t^{-\frac{1}{2k+1}}\big)u_0(z) \, \mathrm{d}z$}.
 $}
\end{equation}

Let us next introduce further rescaling  corresponding to the
variables of the  fundamental solution \eqref{kern},
\begin{equation}
 \label{wwu1}
u(x,t) = t^{-\frac{1}{2k+1}}w(y,\tau) \whereA
y=x/t^{\frac{1}{2k+1}}, \quad \tau = \ln t : \mathbb{R}_+
\rightarrow\mathbb{R},
\end{equation}
where we now have scaling with respect to time as well.  The
rescaled solution $w(y,\tau)$ then satisfies the evolution equation
\begin{equation}
\label{rescale1}
   w_{\tau} = {\bf B}w,
\end{equation}
where ${\bf B}$ is the linear operator for the rescaled kernel
described in \eqref{linODE}. More precisely, {\bf B} is still a
formal linear differential expression to be equipped with proper
``boundary conditions" to get its actual domain.

 Here, $w(y,\tau )$ satisfies the Cauchy
problem for \eqref{rescale1} in $\mathbb{R}\times\mathbb{R}_+$,
with initial data at $\tau = 0$, i.e., at $t=1$ and not $t=0$. So
now the initial data are given by
\begin{equation}
\label{init1}
   w_0(y) = u(y,1) \equiv b(1)\ast u_0 = F\ast u_0.
\end{equation}
Hence, the linear operator $\frac{\partial}{\partial\tau}  - {\bf
B}$ is the rescaled version of the original linear dispersion
operator
 $%%\begin{equation*}
\mbox{$\frac{\partial}{\partial t} +(-1)^{k}D^{2k+1}_x$}.
  $ %%%\end{equation*}
      Therefore, the corresponding semigroup $\mathrm{e}
^{\bf{B} \tau}$ admits an explicit integral representation. This
helps to establish some properties of ${\bf B}$, including key
spectral ones, and describe other evolution features of the linear
flow.

 Thus,
rescaling convolution \eqref{convolution} gives the explicit
representation of the semigroup
\begin{equation}
 \label{ww21}
 \mbox{$
w(y,\tau ) = \int\limits_{\mathbb{R}} \mbox{$F\big(y - z\mathrm{e}
^{-\frac{\tau }{2k+1}}\big)u_0(z) \, \mathrm{d}z$} \equiv
\mathrm{e} ^{{\bf B}\tau}u(y,1) \quad \text{for} \quad \tau \geq
0. $}
\end{equation}
For any $y \in \mathbb{R}$, Taylor's power series   for the
analytic kernel $F$ can be used to expand the convolution and to
obtain
\begin{equation}
 \label{Fexp1}
%%\begin{split}
 \mbox{$
\mbox{$F\big(y-z\mathrm{e} ^{-\frac{\tau }{2k+1}}\big)$}
%% &
= \sum\limits_{(\beta )} \mbox{$\mathrm{e} ^{-\frac{|\beta |\tau
}{2k+1}}\frac{(-1)^{|\beta |}}{\beta !}D_y ^{\beta }F(y)z^{\beta
}$}
%%\\ &
\equiv \sum\limits_{(\beta )} \mbox{$\mathrm{e} ^{-\frac{|\beta
|\tau}{2k+1}}\frac{1}{\sqrt{\beta !}}\psi _{\beta }(y)z^{\beta
}$}.
%%\end{split}
$}
\end{equation}
%% where $|\beta|$ stands for $\beta \ge 0$
Here, for the first time, we reveal the {\em eigenvalues} and {\em
eigenfunctions} of the operator ${\bf B}$ as
\begin{equation}
\label{Eigfuncs}
 \mbox{$
  \l_\b= - \frac {|\b|}{2k+1} \andA
 $}
 \mbox{$\psi _{\beta }(y) =
\frac{(-1)}{\sqrt{\beta !}}^{|\beta |}D_y ^{\beta }F(y)$} \quad
\mbox{for any (multiindex)} \,\,\,\beta, \,\,\,
%% \mbox{with}
 %%\,\,
 |\b| \ge 0.
\end{equation}

\ssk

\noi{\bf Remark: on extensions to $\ren$.} We note that whilst
$|\beta|\equiv\beta$ here, this is not the case in the
multi-dimensional $\mathbb{R}^N$, where $\beta$ stands for a
multiindex. However, we often still keep the notation of
$|\beta|$, to show that the theory here can be extended to
multi-dimensional spaces. We bear in mind that Hermitian spectral
theory can be properly developed for the typical linear dispersion
operator in $\mathbb{R}^N$,
\begin{equation}
\label{multidimop}
   \mbox{${\bf B} =
(-1)^{k+1}\frac{\partial}{\partial y_1}\Delta^k  +
\frac{1}{2k+1}\,y\cdot\nabla_y  + \frac{N}{2k+1}\,I$},
\end{equation}
which appears after the  scaling (similar to \ef{wwu1}) of the
following LDE:
\begin{equation}
 \label{LDEN}
\mbox{$ u_t=(-1)^k \frac{\partial}{\partial x_1} \Delta^k u, $}
\end{equation}
where analogies of the operator on the right-hand side occur in
completely integrable equation theory (cf. the
Kadomtsev--Petviashvili equation for $N=2$ and others from
higher-order hierarchies). The multi-index is then $\beta =
(\beta_1,\hdots,\beta_N)$, with the length $|\beta| =
\beta_1+\hdots+\beta_N$.  With this notation, some of our basic
results can be directly translated to operators such as
\eqref{multidimop} in $\mathbb{R}^N$; see \cite{SemiPDE, 2mSturm}
for necessary technical details.
%% From now onwards, this notation
%%is not used.

\ssk

Thus, without loss of generality, we continue to develop a 1D
theory. Looking back at the expansion of the convolution \ef{ww21}
via \ef{Fexp1}, the solution of \eqref{rescale1} is given by the
series (understood formally still):
\begin{equation}
 \mbox{$
\label{wsol}
 w(y,\tau ) = \sum\limits _{(\beta )} \mathrm{e}
^{-\frac{\beta \tau}{2k+1}}M_{\beta }(u_0)\psi _{\beta }(y). $}
\end{equation}
It then follows that  $\lambda _{\beta } = -\frac{\beta }{2k+1}$
and $\psi _{\beta }(y)$ are the eigenvalues and eigenfunctions of
${\bf B}$ as in \ef{Eigfuncs}
%%(with an obscure functional
%%%setting),
 and
\begin{equation}
\label{EigMoments}
 \mbox{$
M_{\beta }(u_0) = \mbox{$\frac{1}{\sqrt{\beta !}}$}
\int\limits_{\mathbb{R}}z^{\beta }u_0(z) \, \mathrm{d}z,
 $}
\end{equation}
then turns out to be
%% Here $M_{\beta }(u_0)$ are the
 corresponding
moments of the initial data $w_0$, i.e., ``scalar products" of
$w_0$ with some ``adjoint polynomials" (eigenfunctions of the
``adjoint" operator $\BB^*$) to be detected shortly together with
a proper metric involved.

%%%%%%%%%%%%%%%%%%%%%%%%%%%%%%%%%%%
\subsection{$\BB$: convergence of series}

We next need to check suitable metrics of convergence of formal
series \ef{wsol}. Those questions were addressed in \cite{SemiPDE}
and in the most general case in \cite{2mSturm} in the purely
``elliptic" (``parabolic") case, when the rescaled kernel $F(y)$
always have exponential decay as $y \to \iy$. In view of
\ef{posasympt}, this is not currently the case as $y \to + \iy$.
We concentrate on this case and refer to the analysis in
\cite[\S~5]{2mSturm} of convergence of the  integrals as $y \to
-\iy$.

First of all, as in \cite{SemiPDE, 2mSturm}, we have that the ODE
 \ef{linODE} for $F$ implies the following estimate of the
eigenfunctions:
 \be
  \label{psi1}
   \mbox{$
    |\psi_\b(y)| \le \frac c{\sqrt{\b !}}\, (1+y)^{\frac {\b}{2k}-
    \frac{2k-1}{4k}} \a_k^{\b}
    \quad \mbox{for} \quad y \ge 0 \whereA \a_k=\,\frac{2k+1}{2k}\, d_k \sin
    b_k.
    $}
    \ee
Here, the main growing factor $\sim y^{\frac{\b}{2k}}$ as $y \to +
\iy$, as well as the multiplier $\a_k^{\b}$, come from
differentiating $\b$ times inside the $\cos$ function in
\ef{posasympt}.

\ssk

We next claim that the series in \ef{wsol} converges {\em
uniformly on
 compact subsets} for data from the  weighted $L^2$ space (e.g.,
 with  fast exponential decay at infinity):
  \be
  \label{u01}
   \mbox{$
  u_0 \in L^2_{\hat\rho}(\re), \quad \mbox{where} \quad \hat \rho(y)={\mathrm
  e}^{a|y|^\a}, \quad \a= \frac{2k+1}{2k},
  %% \quad (y>0),
   $}
    \ee
 and $a>0$ is a sufficiently small constant. As in the parabolic case
 \cite{SemiPDE},
  in view of the exponential decay \ef{ymin}, for $y<0$ the
  condition  $a \in (0,2 \hat d_k)$ is necessary, but this case  will not be treated here as being
 standard.

 %%Again, \ef{u01} assumes $y>0$, while
 %%for $y<0$, it suffices to fix a usual weight from
 %%\cite{SemiPDE}:
 %%\be
 %% \label{u01+}
 %%  \mbox{$
%%  u_0 \in L^2_{\hat \rho}, \quad \mbox{where} \quad \hat \rho(y)={\mathrm
%%  e}^{-a|y|^\a}, \quad \a= \frac{2k+1}{2k} \quad (y<0),
%%   $}
%%    \ee
%% and $a \in (0,2d_k)$, which we will not treat as being
%% standard.

The proof of uniform convergence on compact subsets uses typical
H\"older estimates such as (for $y>0$ again)
 \be
 \label{H01}
  \mbox{$
   \int |z|^\b |u_0| \equiv \int \frac 1{\sqrt {\hat \rho}}\, |z|^\b
   \sqrt {\hat \rho} |u_0| \le \sqrt{\int \hat \rho |u_0|^2} \,\, \sqrt{\int \frac 1
   {\hat \rho} \, |z|^{2\b}}
 $}
  \ee
  and, using standard properties of the $\Gamma$-function and Stirling's formula,  we estimate
  the last integral as follows:
  \be
  \label{H01N}
   \tex{
  %%%   \,\,\mbox{and}\,\,
\int \frac 1
   {\hat \rho} \, |z|^{2\b} \sim \Gamma \big(\frac{2k\b}{2k+1}\big)
   \sim \big(\frac{2k\b}{2k+1}\big)! \sim  \big(\frac{2k\b}{(2k+1){\mathrm e}}\big)^{\frac{2k\b}{2k+1}}.
    }
 \ee
  Overall, this gives the following majorizing series for
\ef{wsol} for $y \in (0,L)$:
 \be
 \label{H02}
  \mbox{$
 \sum_{(\b)} \frac 1 {\sqrt{\b !}}\, \sqrt{\big(\frac{2k
 \b}{2k+1}\big)!} \,\, L^{\frac \b{2k}},
  $}
  \ee
 which obviously converges for arbitrary fixed $L >0$.

 \ssk

We next show that the series \ef{wsol} converges  {\em in the
mean} in the metric of the weighted space
 \be
 \label{H03}
  \mbox{$
  L^2_{\rho}(\re)                                                                                                                    , \quad \mbox{where}
 \quad
 \rho(y) = \begin{cases} \mathrm{e}^{a|y|^{\alpha}} &\text{for $y\leq
-1$},
 \\ \mathrm{e}^{-ay^{\alpha}} &\text{for $y\geq 1$},
\end{cases}
 $}
 \ee
 %% Note that the standard ``adjoint relation" of the  weights
 %% $\rho(y)= \frac 1{\rho^*(y)}$
 %%\quad \mbox{for both $y>0$ and $y<0$}
 %% $}
 %%  \ee
   cf. \ef{u01}.
   %% and \ef{u01+}.
   Using the same
   estimates \ef{psi1} and \ef{H01}, we then derive the following
   more sensitive geometric majorizing series:
    \be
    \label{H04}
     \tex{
     \sum_{(\b)} \g_k^{\frac{\b}{2k+1}}, \quad \mbox{where}
      \quad \g_k=  \frac
      1{2k+1}\big(\frac{2k}{2k+1}\big)^{2k-1}<1,
       }
       \ee
 so the series converges. We do not include in $\g_k$
 %%%%into account
 the extra multipliers $\sim |d_k \sin b_k|<1$ shown in
 \ef{psi1}, since the first one $\frac{2k}{2k+1}<1$ already suffices
 for the
 convergence of \ef{H04}.

%%%%%%%%%%%%%%%%%%%%%%%%%%%%%%%%%%%
\subsection{$\BB$: alternative representation of semigroup and series, generalized
Hermite polynomials}

Let us discuss an equivalent explicit representation of the
semigroup for ${\bf B}$, which, though being more difficult, will
more clearly determine the eigenfunctions of the adjoint operator
${\bf B}^{\ast}$ to be introduced and studied  next  by a simpler
direct approach. We perform another rescaling to exclude the
relation \eqref{init1}, in order to find the correct semigroup,
corresponding to the initial data at $t=0$, i.e., $w_0=u_0$:
\begin{equation*}
u = (1+t)^{-\frac{1}{2k+1}}w, \quad y = x(1+t)^{-\frac{1}{2k+1}},
\quad \tau = \ln{(1+t)}:\mathbb{R}_+\rightarrow\mathbb{R}_+.
\end{equation*}
Then rescaling the convolution gives
\begin{equation}
%%\begin{split}
\label{rescconv}
 \mbox{$
w(y,\tau )
%%&
= \mathrm{e} ^{{\bf B}\tau}u_0
%%\\
%%&
 \equiv(1-\mathrm{e} ^{-\tau })^{-\frac{1}{2k+1}}\int\limits_{\mathbb{R}}
\mbox{$F\big((y-z\mathrm{e} ^{-\frac{\tau }{2k+1}})(1-\mathrm{e}
^{-\tau})^{-\frac{1}{2k+1}}\big)w_0(z) \, \mathrm{d} z$}.
%%\end{split}
$}
\end{equation}
Once again we can look to find explicit representations for the
eigenfunctions, eigenvalues, and adjoint eigenfunctions, given in
the dual products $\langle u_0, \psi_\beta ^\ast\rangle$, in the
standard metric of $L^2$ (actually, a more delicate indefinite
metric will be needed; see below). It will be shown later that we
can actually determine the adjoint eigenfunctions,
$\{\psi_\beta^\ast\}$, using a much easier method.

Meantime, let us perform some convenient manipulations. Looking at
our rescaled equation \eqref{rescconv}, by Taylor's expansion we
have
\begin{equation*}
%%\begin{split}
 \mbox{$
F\big((y-z\mathrm{e} ^{-\frac{\tau}{2k+1}})(1-\mathrm{e}
^{-\tau})^{-\frac{1}{2k+1}}\big) = \sum _{(\mu)
}\mbox{$\frac{(-1)}{\mu !}^{\mu }D_y ^{\mu }F\big(y(1-\mathrm{e}
^{-\tau })^{-\frac{1}{2k+1}}\big)(\mathrm{e} ^{\tau
}-1)^{-\frac{\mu }{2k+1}}$} z^\mu
 %%\end{split}
  $}
\end{equation*}
%%and
\begin{equation*}
 \mbox{$
\mbox{and} \quad \mbox{$F\big(y(1-\mathrm{e}
^{-\tau})^{-\frac{1}{2k+1}}\big)$} = \sum _{(\nu)
}\mbox{$\frac{1}{\nu !}\, (D_y ^{\nu }F)(0)y^{\nu}(1-\mathrm{e}
^{-\tau })^{-\frac{\nu }{2k+1}}$}.
 $}
\end{equation*}
Then, using these expansions, our solution is given by
\begin{equation}
\begin{split}
\label{expwsol}
 \mbox{$
w(y,\tau ) = (1-\mathrm{e} ^{-\tau})^{-\frac{1}{2k+1}}\sum _{(\mu
,\, \nu) }
 $}
&\mbox{$\frac{(-1)}{\mu !\nu !}^{\mu }\, (D_y^{\nu
}F)(0)D_y^{\mu}(y^{\nu})(\mathrm{e} ^{\tau }-1)^{-\frac{\mu
}{2k+1}}$}
\\&
\mbox{$\times (1-\mathrm{e} ^{-\tau })^{-\frac{\nu }{2k+1}}$}
 \mbox{$
\int_{\mathbb{R}}
 $}
z^{\mu }w_0(z) \, \mathrm{d} z.
\end{split}
\end{equation}
Rearranging this, we have
\begin{equation}
 \label{BwFull}
 %%\begin{split}
  \mbox{$
w(y,\tau ) = \sum _{(\mu ,\, \nu) }\mbox{$ \mathrm{e} ^{
-\frac{\mu \tau}{2k+1}}(1-\mathrm{e}
^{-\tau})^{-\frac{\mu+\nu+1}{2k+1}}\,\frac{(-1)}{\mu !\nu !}^{\mu
} \,(D_y^{\nu}F)(0)D_y^{\mu}(y^{\nu})$}\,
 \mbox{$
\int_{\mathbb{R}}
 $} z^{\mu }w_0(z) \, \mathrm{d} z.
 %%%\end{split}
  $}
\end{equation}
 Further expansion  in the exponential term containing $(1-\mathrm{e}
^{-\tau})^{-\frac{\mu+\nu+1}{2k+1}}$ in terms of ${\mathrm
e}^{-\tau}$ yields an alternative more
%% more complicated and
complete eigenfunction expansion of the semigroup with a full
representation of all the eigenfunctions and the generalized
Hermite polynomials $\{\psi_\b^*\}$ in the moments $M_\b$.
 Moreover, being re-written in the form
  \be
  \label{wwNew}
   \tex{
  w(y,\t)= \sum_{(\b)} {\mathrm e} ^{- \frac {\b \t}{2k+1}}
  \psi_\b(y) \, \langle w_0, \psi_\b^* \rangle_*,
   }
   \ee
 this expansion, besides the adjoint polynomials
 $\Phi^*=\{\psi_\b^*\}$, also defines the corresponding
 {\em indefinite} metric $\langle \cdot, \cdot \rangle_*$ to be
 carefully introduced and explained in Section \ref{SectB*}, where
  those polynomials will be obtained  via a
 simpler direct approach.

%%% Finally, we note that it is not difficult

%We see also that
%\begin{equation*}
%\begin{split}
%\mbox{$\big(-\frac{(\mu+\nu+1)}{2k+1}\big)$}&\mbox{$\hdots\big(-\frac{(\mu+\nu+1)}{2k+1}-\phi+1\big)$}\\
%&=\mbox{$(-1)^{\phi}\big(\frac{(\mu+\nu+1)+(2k+1)\phi-(2k+1)}{2k+1}\big)\hdots\big(\frac{\mu+\nu+1}{2k+1}\big)$}\\
%&=\mbox{$\frac{(-1)^{\phi}}{(2k+1)^{\phi}}\big[(\mu+\nu+1)+(2k+1)\phi-(2k+1)\big]\hdots\big(\mu+\nu+1\big)$},
%\end{split}
%\end{equation*}
%hence we can simplify our solution further to

%\begin{equation*}
%\begin{split}
%w(y,\tau ) = \sum _{(\mu ,\, \nu ,\, \phi) }&\mbox{$
%\mathrm{e}^{\lambda_l\tau}\frac{(-1)^{2\phi}}{(2k+1)^{\phi}}\big[(\mu+\nu+1)+(2k+1)\phi-(2k+1)\big]\hdots
%(\mu+\nu+1)\,\frac{(-1)^\mu}{\mu! \phi !}$}\\
%& \times \mbox{$D_y^{\mu}(\sum y^{\nu}D_y^{\nu}
%F(0))$}\int_{\mathbb{R}} z^{\mu }w_0(z) \, \mathrm{d} z.
%\end{split}
%\end{equation*}

%%The moments $\int_{\mathbb{R}} z^{\mu}w_0 \, \mathrm{d}z$ are finite
%%for all continuous data $w_0$ with sufficient decay at infinity
%%(say, with compact support).  Whilst due to the expansion of terms
%%using Taylor's Series, it suggests an infinite sum, we expect that
%%the interaction between the summations will give us a finite sum.

%The precise meaning of convergence
%(uniformly on compact subsets) is under review.\\[5mm]

%%%%%%%%%%%%%%%%%%%%%%%%%%%%%%%%%%%%%%%%%
 \subsection{Classification of large-time asymptotics}
  \label{S3.4}

 As a by-product, for our LDEs, we have described the large-time
asymptotic behaviour for the problem
\begin{equation}
%%\begin{cases}
\label{fulldeflinpde} u_t = (-1)^{k+1}D^{2k+1}_x u \inB \re \times
\re_+ \whereA u_0\in L^2_{\hat \rho},
%%% \quad \rho(y)=\mathrm{e}^{a|y|^{\frac{2k+1}{2k}}},
%%%\end{cases}
\end{equation}
%%with initial data $u_0$ and
and the weight  $\hat \rho$ is as in \ef{u01}.
%%% and \ef{u01+}.

%%%and the constant $a>0$ are to be properly defined later on.

\begin{theorem}
 \label{Th.B}
For the problem \eqref{fulldeflinpde},
 the eigenfunction expansion
 \eqref{wsol} implies that, for any $u_0\in L^2_{\hat \rho}(\re)$, $u_0 \not
=0$, there exists a finite $l$ such that, as $t \to +\infty$,
\begin{equation*}
 \tex{
u(x,t)= t^{-\frac{1}{2k+1} +
\lambda_\beta}\big[c_\beta\psi_\beta(x/t^{\frac{1}{2k+1}})+o(1)\big]
\whereA  \l_\b= -\frac l{2k+1}
 }
\end{equation*}
 and $l=|\beta|$ is the first eigenvalue index, with  the
corresponding  moment $c_\b=M_\beta(u_0)\neq 0$.
\end{theorem}

The following uniqueness conclusion is straightforward and we keep
this as a simple illustration for further results:

\begin{corollary}
Assume that, for any constant $K >0$, the solution of
$\eqref{fulldeflinpde}$ satisfies
\begin{equation*}
 \mbox{$
\sup_x |u(x,t)| = o(t^{-K}) \quad \text{as} \quad t\to\infty.
%%%\forall k.
$}
\end{equation*}
%%Hence as $t\to\infty$ we have that
Then $u(x,t) \equiv 0$.
%%%\\[5mm]
\end{corollary}

Such results belong to classic  Carleman--Agmon-type estimates in
operator theory: if a solution of a linear equation, under proper
conditions on operators involved, decays super-exponentially fast
(in terms of $\tau= \ln t$) as $t \to +\iy$, then it is trivial.
For elliptic equations $P(x,D)u=0$, this has the natural
counterpart on {\em strong unique continuation property} saying
that nontrivial solutions cannot have zeros of infinite order; a
result first proved by Carleman in 1939 for $P=-\D + V$, $V \in
L^\iy_{\rm loc}$, in $\re^2$ \cite{Carl39}; see \cite{Ion04,
Dos05, Tao08} for further references and modern extensions.

\ssk

 Thus, for convenience,  we fix again some conclusions achieved
 above:
 %%%%conclude with the following suggestions:
\begin{itemize}
\item
There exists point spectrum \{$\lambda_\beta = -\frac{\beta}{2k+1},
\, \beta \ge 0$\} of non-symmetric operator ${\bf B}$.
\item
We have no integral terms in the expansion, hence the spectrum is
expected to be discrete in a proper functional setting.
\item
The set of eigenfunctions $\{\psi_\beta\}$ seems to be complete
and closed in a certain weighted $L^2$-space.
\item
Traces of the polynomials, which give rise to the adjoint
eigenfunctions, $\{\psi_\beta^\ast\}$. %- to be studied and eventually
%proved.
\end{itemize}

%%%%%%%%%%%%%%%%%%%%%%%%%%%%%%%%%%%%%%%%%%%%%%
%%%%%%%%%%%%%%%%%%%%%%%%%%%%%%%%%%%%%%%%%%%%%%
%%%%%%%%%%%%%%%%%%%%%%%%%%%%%%%%%%%%%%%%%%%%%%
%%%%%%%%%%%%%%%%%%%%%%%%%%%%%%%%%%%%%%%%%%%%%%
%%%%%%%%%%%%%%%%%%%%%%%%%%%%%%%%%%%%%%%%%%%%%%
%%%%%%%%%%%%%%%%%%%%%%%%%%%%%%%%%%%%%%%%%%%%%%

%%%%%%%%%%%%%%%%%%%%%%%%%%%%%%%%%%%%%%%%%%%%%%%%%%%%%%%%%%%%%%%%%%%%%%%%
\subsection{Semigroup of the adjoint operator ${\bf B^\ast}$ and its eigenfunction expansion}
 \label{S3.2}

We now find the explicit representation of the semigroup $\mathrm{e}
^{{\bf B^{\ast} }\tau}$, where ${\bf B}^*$ is obtained from the LDE
in \eqref{fulldeflinpde} by using other {\em blow-up} rescaling. Let
us introduce the rescaled variables
\begin{equation}
 \label{ScBl1}
u(x,t) = w(y,\tau ), \quad y=x/(1-t)^{\frac{1}{2k+1}}, \quad \tau
=-\ln{(1-t)}:(0,1)\rightarrow \mathbb{R} _+,
\end{equation}
where, for convenience, the blow-up time is $T=1$. Then $w$ now
solves the problem
\begin{equation*}
w_{\tau } = {\bf B^{\ast}}w \quad \text{for} \quad \tau > 0,
\end{equation*}
with data $w(y,0)=w_0(y) \equiv u_0(y)$.  Here the ``adjoint"
operator ${\bf B^\ast}$ is given by
\begin{equation}
 \label{adjoint}
{\bf B^\ast}=\mbox{$(-1)^{k+1}D_y^{2k+1}-\frac{1}{2k+1}\,yD_y$}.
\end{equation}
By rescaling the convolution \eqref{convolution}, we have
\begin{equation*}
 \mbox{$
w(y,\tau )=(1-\mathrm{e} ^{-\tau
})^{-\frac{1}{2k+1}}\int_{\mathbb{R}} \mbox{$F\big((y\mathrm{e}
^{-\frac{\tau}{2k+1}}-z)(1-\mathrm{e}
^{-\tau})^{-\frac{1}{2k+1}}\big)w_0(z) \, \mathrm{d} z$}.
 $}
\end{equation*}
Using Taylor's expansion yields
\begin{equation*}
 \mbox{$
\begin{split}
&F \mbox{$\big((y\mathrm{e} ^{-\frac{\tau }{2k+1}}-z)(1-\mathrm{e}
^{-\tau })^{-\frac{1}{2k+1}}\big)$}
\\ = &
  \mbox{$
   \sum_{(\beta)}
    $}
    \mbox{$\frac{(-1)}{\beta !}^{\beta
}(1-\mathrm{e} ^{-\tau })^{-\frac{\beta }{2k+1}}D^{\beta
}F\big(y\mathrm{e} ^{-\frac{\tau }{2k+1}}(1-\mathrm{e}
^{-\tau})^{-\frac{1}{2k+1}}\big)z^{\beta }$}
\end{split}
$}
\end{equation*}
and expanding in $y$ leads to
\begin{equation*}
 \mbox{$
\mbox{$F\big(y\mathrm{e} ^{-\frac{\tau }{2k+1}}(1-\mathrm{e}
^{-\tau})^{-\frac{1}{2k+1}}\big)$} = \sum _{(\nu)}
\mbox{$\frac{1}{\nu !}D^{\nu }F(0)y^{\nu}(\mathrm{e} ^{\tau
}-1)^{-\frac{\nu }{2k+1}}$}.
 $}
\end{equation*}
So the solution is represented as
\begin{equation}
 %% \mbox{$
\begin{split}
\label{expusol}
%%% \mbox{$
 \mbox{$w(y,\tau)=(1-\mathrm{e} ^{-\tau
})^{-\frac{1}{2k+1}}$}&
  \mbox{$
\sum_{(\beta ,\, \nu)}
 $}
\mbox{$\frac{(-1)^{\beta }}{\beta !\nu !}(1-\mathrm{e}
^{-\tau})^{-\frac{\beta }{2k+1}}(\mathrm{e} ^{\tau
}-1)^{-\frac{\nu }{2k+1}}$}\\ \times&\mbox{$ D^{\nu
}F(0)\frac{1}{(\nu - \beta )!}y^{\nu - \mu }$} \mbox{$
\int_{\mathbb{R}} $} \mbox{$z^{\beta }w_0(z) \, \mathrm{d}z$}.
 %%$}
\end{split}
%% $}
\end{equation}
%As with the case of $w(y,\tau)$ in \eqref{expwsol}, the convergence
%of $u(x,t)$ in \eqref{expusol} is still under review, with the exact
%definition of adjoint eigenfunctions $\{\psi_{\beta}^{\ast}\}$ here
%to be determined.

%\begin{center}
%\begin{tabular}{|c|c|c|}
%\hline $l$ & $\lambda $ & $\psi ^*$\\ \hline
%0 & 0 & 1\\ \hline
%1 & $-\frac{1}{2k+1}$ & $y$\\ \hline
%2 & $-\frac{-}{2k+1}$ & $y^2$\\ \hline
%$\vdots$ & $\vdots$ & $\vdots$ \\ \hline
%$m < 2k+1$ & $-\frac{m}{2k+1}$ & $y^m$\\ \hline
%$2k+1$ & -1  & $y^{2k+2}+c_1$ \\ \hline
%$2k+2$ & $-\frac{(2k+2)}{2k+1}$ & $y^{2k+1}+c_2y$\\ \hline
%$2k+n$ & $-\frac{(2k+n)}{2k+1}$ & $y^{2k+n}+c_ny$\\ \hline
%\end{tabular}
%\end{center}

Again, similar to the expansion \ef{wwNew} for the operator ${\bf
B}$,
%% in the previous subsection,
  (\ref{expusol}) can be viewed as
an eigenfunction expansion of the solution that can reveal many
key spectral properties of ${\bf B}^*$.
 In particular, we obtain the same real spectrum and eigenfunctions:
  \be
  \label{SpB*}
  \tex{
 \BB^*: \quad\l_\b= - \frac \b{2k+1} \andA \mbox{Hermite-type polynomial  eigenfunctions
 $\psi_\b^*(y)$}.
 }
 \ee
 However, further refining of this expansion will lead to more
complicated formulae, which are not that effective and useful.
Therefore, we return later on to polynomial eigenfunctions of ${\bf
B}^*$ on the basis of a simpler direct approach.

%%%%%%%%%%%%%%%%%%%%%%%%%%%%%%%%%%%%%%%%%%%%%%%%%%%%
\subsection{Classification of multiple zeros: micro-scale blow-up
asymptotics}

Due to the blow-up character of the scaling \ef{ScBl1}, the
polynomial eigenfunction set $\Phi^*=\{\psi_\b^*\}$ of the adjoint
operator $\BB^*$ is able to give insight into the structure of
multiple spatial zeros of solutions $u(x,t)$ of the corresponding
LDEs. This blow-up asymptotic theory is a natural counterpart of
the large-time one in Section \ref{S3.4} induced by eigenfunctions
$\Phi$ of $\BB$. The key principles of such a classification for
LDEs can be found in \cite[\S~9]{2mSturm}, where other blow-up
application of such a spectral analysis can be found.

%%%%%%%%%%%%%%%%%%%%%%%%%%%%%%%%%%%%%%%%%%%%%%%%%
%%%%%%%%%%%%%%%%%%%%%%%%%%%%%%%%%%%%%%%%%%%%%%%%%
%%%%%%%%%%%%%%%%%%%%%%%%%%%%%%%%%%%%%%%%%%%%%%%%%
%%%%%%%%%%%%%%%%%%%%%%%%%%%%%%%%%%%%%%%%%%%%%%%%%
%%%%%%%%%%%%%%%%%%%%%%%%%%%%%%%%%%%%%%%%%%%%%%%%%
%%%%%%%%%%%%%%%%%%%%%%%%%%%%%%%%%%%%%%%%%%%%%%%%%
%%%%%%%%%%%%%%%%%%%%%%%%%%%%%%%%%%%%%%%%%%%%%%%%%

%%%%%%%%%%%%%%%%%%%%%%%%%%%%%%%%%%%%%%%%%%%%%%%%%%%%%%%%%%%%%%
\section{Hermitian spectral theory:
%%discrete point spectrum of
  operator {\bf B}}
 \label{S2.3}

 We now start more systematically to develop necessary spectral
theory for the operator pair $\{{\bf B},\,{\bf B}^*\}$ introduced
above. In some aspects, this theory repeats standard steps of
self-adjoint theory for the  Hermite classic operator
(\ref{Herm1}), which since the nineteenth century, is associated
with the names of Sturm, Hermite, and other famous mathematicians.

%%%%%%%%%%%%%%%%%%%%%%%%%%%%%%%%%%%%%%%%%%%%%%%%%%%%%%%%%%%
\subsection{Bounded operator $\BB$}
 According to \ef{H03} and \ef{u01},
 %%% \ef{u01+},
we calculate the spectrum of the linear operator ${\bf B}$  in the
weighted space $L^2_{\rho}(\mathbb{R})$, with an exponential
weight given in \ef{H03},
%%\begin{equation}
%% \label{rr*1}
%%\rho(y) = \begin{cases} \mathrm{e}^{a|y|^{\alpha}} &\text{for $y\leq
%%-1$},
%% \\ \mathrm{e}^{-ay^{\alpha}} &\text{for $y\geq 1$},
%%\end{cases}
%%\end{equation}
where, as usual, we assume $\rho (y)$ to be sufficiently smooth in
the complement interval $[1,-1]$. Here we have that $\rho (y)>0$
and $a\in (0,2 \hat d_k)$ (this is necessary for $y \ll -1$ only)
is a sufficiently small positive constant, where $\hat d_k$ is
defined as before, in \eqref{posasympt}, \ef{yminNN}.

We introduce a Hilbert space of functions
$H^{2k+1}_{\rho}(\mathbb{R})$ with the inner product and the
induced norms:
\begin{equation*}
  \mbox{$
\langle v, w \rangle_{2k+1, \rho} =
\int\limits_{\mathbb{R}}\rho(y)\sum\limits _{r=0}
^{2k+1}D_y^rv(y)D_y ^r w(y) \, \mathrm{d} y, \quad
 $}
%%\end{equation*}
 %%%%and  the induced norm is
%%\begin{equation*}
  \mbox{$
\|v\|^2_{2k+1, \rho} = \int\limits_{\mathbb{R}}\rho(y)\sum\limits
_{r=0} ^{2k+1}|D_y^rv(y)|^2 \, \mathrm{d} y.
 $}
\end{equation*}
 It turns out that it suffices to restrict to the real case.
 Then $H^{2k+1}_{\rho}\subset L^2_{\rho}(\mathbb{R})\subset
L^2(\mathbb{R})$. The first conclusion is standard as in
\cite{SemiPDE}:

\begin{lemma}
\label{bddop} {\bf B} is a bounded linear operator from
$H^{2k+1}_{\rho}(\mathbb{R})$ to $L^2_{\rho}(\mathbb{R})$.
\end{lemma}

\noindent {\em Proof.} Consider
 the differential expression
 %%%%the linear operator given by
\eqref{linODE},
\begin{equation*}
\mbox{${\bf B}v = (-1)^{k+1}v^{(2k+1)}
+\frac{1}{2k+1}\,yv^{\prime}+\frac{1}{2k+1}\,v$}.
\end{equation*}
For ${\bf B}$ to be bounded, it is necessary to look at the second
term with the unbounded coefficient $y$.  In order to do this, we
want to show that
\begin{equation*}
  \mbox{$
\int\rho (yv^{\prime})^2 \leq C\int \rho (v^{(2k+1)})^2 \,
\mathrm{d} y,
 $}
\end{equation*}
for some constant $C>0$. To show this, we look at the non-negative
integral
\begin{equation*}
  \mbox{$
0 \leq \int\rho(v^{\prime} + y^{\gamma}v)^2 \, \mathrm{d} y
=\int\rho\big((v^{\prime})^2 + y^{2\gamma}v^2 +
2y^{\gamma}v^{\prime}v\big)\, \mathrm{d}y,
 $}
\end{equation*}
where $\gamma>0$ is some unknown exponent. Integrating by parts
 (allowed for functions in $H^{2k+1}_{\rho}(\mathbb{R})$) in the last
term yields
\begin{equation*}
  \mbox{$
2\int (\rho y^{\gamma})v^{\prime}v\, \mathrm{d}y = \int (\rho
y^{\gamma})(v^2)^{\prime}\, \mathrm{d}y
 = -\int v^2 (\rho y^{\gamma})^{\prime}\, \mathrm{d}y.
  $}
\end{equation*}
%%We can integrate by parts for functions in the spaces
%%%$H^{2k+1}_{\rho}(\mathbb{R})$ and $L^2_{\rho}(\mathbb{R})$.
Then for the exponential weight,
%%% we have
%%\begin{equation*}
%%\rho (y) = \mathrm{e}^{-ay^{\alpha}}
%%\end{equation*}
we have to have  that
\begin{equation*}
\rho (y) = \mathrm{e}^{-ay^{\alpha}}
 \quad \Longrightarrow \quad
 (\rho y^{\gamma})^{\prime}
\sim \mathrm{e}^{-ay^{\alpha}}y^{\gamma + \alpha -1},
\end{equation*}
for $y\gg 1$ (and also for  $y\ll -1$ that is similar, where we
replace $y\mapsto|y|$). Hence
\begin{equation*}
 \mbox{$
\int \rho (v^{\prime})^2\, \mathrm{d}y + \int \rho y^{2\gamma}v^2\,
\mathrm{d}y +C_1\int \rho y^{\gamma +\alpha - 1}v^2\, \mathrm{d}y
\geq 0.
 $}
\end{equation*}
By equating powers of $y$, this yields
%%\begin{equation*}
$\gamma = \alpha -1$.
%%%\end{equation*}
Substituting this $\gamma$, we have that
\begin{equation*}
 \mbox{$
\int \rho y^{2(\alpha -1)}v^2 \, \mathrm{d}y\leq C_2\int \rho
(v^{\prime})^2\, \mathrm{d}y.
 $}
\end{equation*}

In particular, for $k=1$, by the {\em Hardy-type inequality}, the
following holds:
\begin{equation*}
 \mbox{$
\int \rho y^{4(\alpha -1)}(v^{\prime})^2\, \mathrm{d}y \leq
C_3\int \rho y^{2(\alpha -1)}(v^{\prime\prime})^2 \, \mathrm{d}y
\leq C_3^2\int\rho (v^{\prime\prime\prime})^2\, \mathrm{d}y.
 $}
\end{equation*}
However, we want that
\begin{equation*}
 \mbox{$
\int \rho y^2(F^{\prime})^2 \, \mathrm{d}y\leq C_4\int\rho
(F^{\prime\prime\prime})^2\, \mathrm{d}y.
 $}
\end{equation*}
Hence we have
%\begin{equation*}
$4(\alpha -1) = 2$,
%%\end{equation*}
which gives
%%\begin{equation*}
$\alpha = \mbox{$\frac{3}{2}$}$.
%%\end{equation*}

For the general case, we have $2k$ iterations and so
\begin{equation*}
4k(\alpha -1) =2 \quad \Longrightarrow \quad \alpha \, =
\mbox{$\frac{2k+1}{2k}$}. \qed
\end{equation*}
%Now we have that the space has exponential weight
%\begin{equation}
%\rho(y) = \begin{cases} \mathrm{e}^{a|y|^{\frac{2k}{(2k+1)}}} &\text{for
%$y\leq -1$}
%\\ \mathrm{e}^{-ay^{\frac{2k}{(2k+1)}}} &\text{for $y\geq 1$}.
%\end{cases}
%\end{equation}

%%\hfill$\Box$
%%\end{proof}

%%%\smallskip

%%%%%%%%%%%%%%%%%%%%%%%%%%%%%%%%%%%%%%%%%%%%%%%%%%
 \subsection{Completeness and compact resolvent}

%%For the spectral results below, we always mean that the
%%differential form for ${\bf B}$ is equipped with the proper
%%``radiation conditions" at infinity, which shall be explained in
%%Section \ref{SectRadcond}. In other words, this precisely defines
%%the necessary domain $\tilde{\mathcal D}({\bf B})$ of the operator
%% ${\bf B}$.

 The following conclusions are also pretty general \cite{SemiPDE,
 2mSturm}:

\begin{lemma}
 \label{L.Compl}
%%%\\
\noindent {\rm (i)} The set of eigenfunctions
$\Phi=\{\psi_{\beta}\}$ given in $(\ref{Eigfuncs})$ is complete in
$L^2(\re)$ and in $L^2_{\hat \rho}(\mathbb{R})$.

\ssk

 \noindent {\rm (ii)}
 %% For any $\lambda
%%\not\in \sigma({\bf B})$,
 The resolvent $({\bf B}-\lambda I)^{-1}$
is a compact operator in $L^2_{\hat \rho}(\mathbb{R})$.
\end{lemma}

\noindent {\em Proof.} (i) In fact, due to the differential
structure of the eigenfunctions in \ef{Eigfuncs}, this is a
standard result in functional analysis; see analytic function
theory approach in Kolmogorov--Fomin \cite[p.~431]{KolmF}. Let us
show that the system of eigenfunctions $\{D^{\beta}F\}$ is
complete in $L^2_{\hat \rho}(\mathbb{R})$; the case of $L^2(\re)$
is simpler. By the Riesz--Fischer theorem, we have to show that,
given a function $G \in L^2_{\hat\rho}(\mathbb{R})$, then
\begin{equation*}
  \mbox{$
\int D^{\beta}F(x)G(x)\, \mathrm{d} x = 0 \quad \text{for any} \quad
\beta,
 $}
\end{equation*}
implies that $G=0$.  Let $\hat F=\mathcal{F}(F)$ and $\hat G
\mathcal{F}(G)$ be the Fourier transforms of $F$ and $G$.
 Then substituting $F={\mathcal F}^{-1}(\hat F)$, $G={\mathcal F}^{-1}(\hat
 G)$,
 and using that $D^\b {\mathcal F}^{-1}(\hat F) \sim \xi^\b {\mathcal F}^{-1}(\hat
 F)$, after integration by parts via ${\mathcal F}^{-1}(\hat
 G)(\xi) \sim {\mathcal F}(\hat
 G)(-\xi)$, we obtain that
\begin{equation}
 \label{FG1}
  \mbox{$
\int \xi^{\beta}\hat{F}(\xi)\hat{G}(-\xi) \, \mathrm{d} \xi = 0
\quad \text{for any} \quad \beta.
 $}
\end{equation}

Let us calculate $\hat F$. Applying the Fourier transform to the
equation \eqref{linODE} yields
\begin{equation}
 \label{Fhat}
  \mbox{$
\mbox{$\mathrm{i}\,|\xi|^{2k}\xi \hat{F} + \frac{1}{2k+1}\,\xi
D\hat{F} = 0$} \LongA \hat{F}(\xi) =
\mathrm{e}^{-\mathrm{i}|\xi|^{2k}\xi}.
%%% \quad \mbox{so that}
 $}
\end{equation}
%%Hence $\hat{F}(\xi) = \mathrm{e}^{-\mathrm{i}|\xi|^{2k}\xi} $. So,
Substituting \ef{Fhat} into \ef{FG1} yields
\begin{equation}
\label{compl}
 \mbox{$
\int \xi^{\beta}\mathrm{e}^{-\mathrm{i}|\xi|^{2k}\xi}
\hat{G}(-\xi) \, \mathrm{d} \xi = 0 \quad \text{for any} \quad
\beta.
 $}
\end{equation}

 Finally, consider the function
\begin{equation*}
  \mbox{$
M(z) = \int \mathrm{e}^{-\mathrm{i}|\xi|^{2k}\xi}\hat
G(-\xi)\mathrm{e}^{\mathrm{i}z\xi}\,\mathrm{d}\xi,
 $}
\end{equation*}
which is entirely analytic in $\mathbb{C}$, \cite{FuncAnaly}. So
\eqref{compl} means that $D^{\beta}M(0)=0$ for any $\beta$.
Therefore, $M(z) \equiv 0$.  Hence $\mathcal{G}(\xi)=0$ almost
everywhere and $G=0$.

\smallskip

\noindent (ii) The proof follows that of the $2m$th-order case,
set out in  \cite[\S~2.3]{SemiPDE}. In fact,  replacing
 $ %%\begin{equation*}
m \mapsto \mbox{$k+\frac{1}{2}$},
 $ %%%\end{equation*}
will yield the same result.
 This also directly follows from the compact
 embedding of the corresponding spaces $H_\rho^{2k+1} \subset
 L^2_\rho$; see Maz'ya's classic monograph on Sobolev spaces
 \cite[p.~40]{Maz}. $\qed$
%%\hfill$\Box$
%%\end{proof}

\smallskip

%%%%%%%%%%%%%%%%%%%%%%%%%%%%%%%%%%%%%%%%%%%%%%%%%
%%%%%%%%%%%%%%%%%%%%%%%%%%%%%%%%%%%%%%%%%%%%%%%%%
%%%%%%%%%%%%%%%%%%%%%%%%%%%%%%%%%%%%%%%%%%%%%%%%%
%%%%%%%%%%%%%%%%%%%%%%%%%%%%%%%%%%%%%%%%%%%%%%%%%
%%%%%%%%%%%%%%%%%%%%%%%%%%%%%%%%%%%%%%%%%%%%%%%%%

%%%%%%%%%%%%%%%%%%%%%%%%%%%%%%%%%%%%%%%%%%%%%%%%%%%%%%%%%%%%%%%%%%%%%
\subsection{Space of closure $\tilde L^2_\rho(\re)$}
 \label{S4.3}

We next apply another direct approach to constructing the
necessary domain of $\BB$; cf. that in \cite{SemiPDE}. Such a
rather ``artificial" construction, though not being completely
unavoidable for such odd-order operator $\BB$ (a more traditional
way of doing this along the lines in \cite{SemiPDE, 2mSturm} is
still available at least partially), will essentially simplify in
Section \ref{SectB*} the way of introducing the set of linear
functionals $\Phi^*=\{\psi_\b^*\}$ from the adjoint space being
generalized Hermite polynomials as eigenfunctions of the adjoint
operator $\BB^*$.

Namely, as in \cite[\S~5]{2mSturm}, given the complete set of
eigenfunctions
 $\Phi=\{\psi_\b\}$
 we
introduce the {\em space of $\Phi$-closure} denoted by $\tilde
L_\rho^2(\re) \subset L^2_\rho(\re)$ as follows:
 \be
 \label{cl1}
\tilde L_\rho^2(\re): \quad \mbox{closure of finite sums
$\sum^M_{\b=0} c_\b \psi_\b$ in the metric of $L^2_\rho(\re)$}.
 \ee
The notation the ``space of closure" is then justified naturally:
 \be
 \label{cl2}
 \mbox{any $v \in \tilde L_\rho^2(\re)$ has the unique
 eigenfunction expansion via $\Phi$: \, $v= \sum_{(\b \ge 0)} c_\b
 \psi_\b$},
  \ee
i.e., the set $\Phi$ is closed therein.  Note that, for  the whole
space $ L_\rho^2(\re)$, the property \ef{cl2} is not easy to check
and, most plausibly, it fails. A functional characterization of
such a space is not that straightforward, and some estimates of
the sequences $\{c_\b\}$, for which the series in \ef{cl2}
converges, can be derived as in \cite[\S~5]{2mSturm}, so we will
not treat those questions here.
 In Appendix A, we present an alternative explanation of some of
 key ``radiation condition" properties for functions from the
 space of $\Phi$-closure.
  Appendix B shows spectral theory for a majorizing integral operator
 for \ef{convolution}, which is order-preserving so some
 comparison features are available.

 \ssk

  We now complete our spectral theory of $\BB$ as follows:

%%%%%%%%%%%%%%%%%%%%%%%%%
 \begin{lemma}
  \label{L.spB}
 %%\noindent{\rm (i)}
  The operator ${\bf B}$  in $(\ref{linODE})$ with the domain $\tilde {\mathcal
  D}(\BB)= \tilde L_\rho^2(\re)\cap H^{2k+1}_\rho(\re)$
  has the point spectrum
%%defined in $(\ref{linODE})$ in $\tilde{\mathcal D}({\bf B})$
comprising  real eigenvalues only:
\begin{equation}
\label{Eigs}
 \sigma({\bf B}) = \big\{\lambda_{\beta} =
-\mbox{$\frac{\beta}{2k+1}$},\, \beta = 0,1,2,...\big\}.
\end{equation}
Eigenvalues are simple with eigenfunctions $(\ref{Eigfuncs})$.
%%\begin{equation*}
%%\mbox{$\psi _{\beta }(y) = \frac{(-1)}{\sqrt{\beta !}}^{\beta }D_y
%%^{\beta }F(y)$}.
%%\end{equation*}
 \end{lemma}

%%%%%%%%%%%%%%%

\noi{\em Proof.} Eigenvalues and eigenfunctions are found by
applying  $D_y ^{\beta}$ to \eqref{linODE}:
\begin{equation*}
\mbox{$D_y ^{\beta}{\bf B}F \equiv {\bf B}D_y ^{\beta}F +
\frac{\beta}{2k+1}\,D_y ^{\beta}F=0$}.
\end{equation*}
By the construction of the space of closure \ef{cl1}, this
completes the proof. $\qed$

\ssk

 Note also that as
 follows from the asymptotic expansion of \eqref{wsol} and
\eqref{expwsol} as $\tau \rightarrow \infty$,  no other
eigenfunctions exist (at least for such data). This again confirms
that all eigenvalues are real and are given in \eqref{Eigs}.

%%%\smallskip

A sharper using Stirling's series  shows that (see estimates in
\cite[\S~5.1]{2mSturm} that are applied to the present case)
 %%It follows that
% there holds
 \be
\label{cc11N}
 \mbox{$
  v= \sum c_\b \psi_\b \in \tilde L^2_\rho \,\,
\Longrightarrow \,\, c_\b = o(|\b|^{|\b|(\nu + \e)}) \quad
\mbox{with an $\e > 0$},
 $}
 \ee
for ``almost all" $|\b| \gg 1$.
 A sharper and almost optimal estimate  reads:
 \be
 \label{opt1}
  \mbox{$
   \sum c_\b
\psi_\b \in
 \tilde L^2_\rho \quad \mbox{if} \quad
 c_\b = O(\d^{|\b|} |\b|^{|\b| \nu} ) \whereA \nu= \frac{2-\a}{2\a}=
 \frac {2k-1}{2(2k+1)},
  $}
  \ee
   with a sufficiently small
 constant
 $\d>0$.
%%Since estimates of the leading terms in (\ref{rho11N}) are sharp,

By $\tilde H^{2k+1}_\rho \subset \tilde L^2_\rho$, we denote the
dense linear subspace obtained as the closure in the norm of
$H^{2k+1}_\rho$ of the subset of eigenfunction expansions with
coefficients satisfying (\ref{cc11N}) or \ef{opt1}. $\tilde
H^{2k+1}_\rho$ with the scalar product of $H^{2k+1}_\rho$ becomes
a Hilbert space and can be considered as the domain of $\BB$ in
$H^{2k+1}_\rho$. There holds
 \be
  \label{11.20}
   \tilde H^{2k+1}_\rho \subseteq
H^{2k+1}_\rho \cap \tilde
 L^2_\rho.
 \ee
%%Note that (\ref{cc11}) does not apply for $m=1$ since then $\a=2$
%%and hence $\nu = 0$. Actually, a natural optimal analogy of
%%$\tilde H^{2m}_\rho$ for $m=1$ is $H^{2}_\rho$, the domain of
%%$\BB$ in $L^2_\rho$.

%%%%%%%%%%%%%%%%%%%%%%%%%%%%%%%%%%%%%%%%%%%%%%%%
\subsection{Little Hilbert spaces}

 We will need a subspace of $\tilde L^2_\rho$ introduced
as a {\em little} Hilbert space
 $l^2_{\rho}$  of functions $v = \sum c_\b
\psi_\b \in \tilde L^2_\rho$ with coefficients satisfying
 \begin{equation}
 \label{abeta1}
  \mbox{$
 \sum |c_\b|^2 < \infty,
  $}
 \end{equation}
where the scalar product  and the induced norm are given by
 \begin{equation}
 \label{vwin}
  \mbox{$
  (v,w)_{0} = \sum c_\b a_\b \,\,\, \mbox{for} \,\, w = \sum a_\b \psi_\b \in l^2_\rho,
  \quad \mbox{and} \,\,\,\,
\|v\|^2_0 = (v,v)_0.
 $}
 \end{equation}
Obviously, $l^2_\rho$ is isomorphic to the Hilbert space $l^2$ of
sequences
 $\{c_\b\}$ with the same inner product, and hence
 \be
 \label{Phiort}
\Phi \quad \mbox{is orthonormal in} \,\,\, l^2_\rho.
 \ee
 It follows from \ef{opt1} that
  \be
  \label{opt2}
   l^2_\rho \subset \tilde L^2_\rho \LongA l^2_\rho \cap \tilde
   L^2_\rho= l^2_\rho.
    \ee

We next define a little Sobolev space $h^{2k+1}_\rho$ of functions
$v \in l^2_\rho$ such that $\BB v \in l^2_\rho$, i.e.,
 $$
  \mbox{$
  \sum |\l_\b c_\b|^2 < \infty.
   $}
  $$
  The scalar product and the induced norm in $h^{2k+1}_\rho$ are
   \be
   \label{5656}
    \mbox{$
  (v,w)_{1} = (v,w)_0 + (\BB v,\BB w)_0 \quad \mbox{and}
  \quad  \|v\|_{1}^2 =
  (v,v)_{1} \equiv \sum (1+|\l_\b|^{2})|c_\b|^2.
   $}
  \ee
This norm is equivalent to the graph norm induced by the positive
operator $(-\BB + a I)$ with  $a>0$. Then $h^{2k+1}_\rho$ is the
domain of $\BB$ in $l^2_\rho$.
% (see also Proposition 2.1 in
% \cite{Eg4}).
 We also have a Sobolev embedding theorem,
 \begin{equation}
  \label{SobE}
 h^{2k+1}_\rho \subset l^2_\rho \quad \mbox{compactly},
 \end{equation}
 which follows from the well-known criterion of compactness in $l^p$:
a  $T \subset l^p$ is  compact  iff
  \be
  \label{lp11}
  \forall \,\,\, \e> 0 \,\,\, \exists \,\,\mbox{integer} \,\, K=K(\e)>0 \,\, \mbox{such that}
   \,\,\forall \,\, \{c_\b\} \in T \,\Longrightarrow \,
    \mbox{$\sum_{|\b| \ge K} |c_\b|^p<\e;
    $}
     \ee
  see,
 e.g., \cite{LustSob}.
 %%%the generalized  Ascoli-Arzel\'a theorem in
 %%\cite[Ch.~2]{KolF}.
 %%In  the self-adjoint case $m=1$, the little space $l^2_\rho$
%%coincides with the big one,
%% \begin{equation}
 %% \label{lLL}
 %%l^2_\rho = L^2_\rho \quad \mbox{for} \,\,\, m=1 \,\,\, \mbox{if\,
%% $a=\frac 14$\, in\, $(\ref{rho11})$}.
%% \end{equation}
%%Then $h^2_\rho$ is the domain $H^2_\rho$ of ${\bf B}$ (recall that
%%if $a \not = \frac 14$, then ${\bf B}$ is not self-adjoint in
%%%$L^2_\rho$).
% and, in general, (\ref{lLL}) is not true, $\tilde
% L^2_\rho \not = L^2_\rho$, even for $m=1$.

 Since the
orthonormality of $\Phi$ is known to be  of importance in operator
theory and applications, in some linear and nonlinear problems
dealing with operators like ${\bf B}$, the little space $l^2_\rho$
can play a special role in comparison with the big one $L^2_\rho$.

It follows from (\ref{vwin}) that  $\BB$ is {\em self-adjoint} in
$l^2_\rho$ with the domain $h^{2k+1}_\rho$,
 \be
 \label{Bsad}
 ( \BB v, w )_0 =  ( v, \BB w )_0 \quad
 \mbox{for all} \,\,\, v,w \in h^{2k+1}_\rho.
 \ee
 Notice that this  {\em a posteriori}  conclusion in a special functional setting  is
 obtained
after establishing all the necessary spectral properties of the
operator.

%% Let us state other
%% straightforward
%% consequences (this list can be easily extended).

%%%%%%%%%%%%%%%%%%%%%%%%%%%%%%%%%%%%%%%%%%%%%%%%%%%%%%%%%%%%%%%%%%%%%%%%%%%%
\section{Spectrum and polynomial eigenfunctions  of the adjoint operator ${\bf B^{\ast}}$}
 \label{SectB*}

We now look to explicitly describe the eigenfunctions of the
``adjoint" operator \ef{adjoint}.
 %%\begin{equation}
 %%\label{adjoint} \mbox{${\bf B^\ast } = (-1)^{k+1} D_y ^{2k+1} -
 %%\frac{1}{2k+1}\,yD_y$}.
 %%\end{equation}
 We recall that $\BB^*$ was derived in Section \ref{S3.2} via the
 blow-up scaling \ef{ScBl1}, and yet has nothing to do with the
 standard adjoint differential form $(\BB)^*$ in the metric of
 $L^2$, which, though  is easily  derived, has no further applications.

%%%%%%%%%%%%%%%%%%%%%%%%%%%%%%%%%%%%%%%%%%%%%%%%%%%%%%%%%%%%%
\subsection{Indefinite metric}

 Thus, before we look at the operator ${\bf B^\ast}$ however, we first
make the following easy observation:

\begin{proposition}
${\bf B^\ast}$ is not adjoint to {\bf B} in the standard metric of
$L^2(\mathbb{R})$.
\end{proposition}

\noindent {\em Proof.} Let $v, w \in C_0^{\infty}(\mathbb{R})$,
then integration by parts yields
\begin{equation*}
 \mbox{$
\begin{split}
\langle {\bf B}v, \, w\rangle &\equiv
 \mbox{$\int\limits_{\mathbb{R}}
 $}
 \mbox{$\big((-1)^{k+1}v^{(2k+1)}+\frac{1}{2k+1}\,(yv)^{\prime}\big)w$}\,\mathrm{d} y\\
&=
  \mbox{$
\int\limits_{\mathbb{R}}
 $}
\mbox{$\big((-1)^{k}v^{(2k)}-\frac{1}{2k+1}\,yv\big)w^{\prime}$}\,\mathrm{d}
y\\
 %%%%&=\quad \vdots\\
&=...=
  \mbox{$
\int\limits_{\mathbb{R}}
 $}
\mbox{$v\big((-1)^kw^{(2k+1)}-\frac{1}{2k+1}
\,yw'\big)$}\,\mathrm{d} y\\ &= \langle v, \, ({\bf B})^\ast
w\rangle, \quad \mbox{where} \quad ({\bf B})^\ast =
\mbox{$(-1)^kD_y^{2k+1}-\frac{1}{2k+1}\,yD_y\neq{\bf B^\ast}$}.
\qed
\end{split}
 $}
\end{equation*}
%%where
%%\begin{equation*}
%%{\bf \tilde{B}^\ast} =
%%\mbox{$(-1)^kD_y^{2k+1}-\frac{1}{2k+1}\,yD_y\neq{\bf B^\ast}$}.
%%\qed
%%\end{equation*}
%%\hfill$\Box$
%%\end{proof}

\smallskip

Thus, in order to get the correct adjoint operator ${\bf B^\ast}$,
it is necessary to use another metric. The scalar product of this
{\em indefinite metric} of the space $\bar L^2$ is given by
\begin{equation}
\label{indefprod}
 \mbox{$
 \langle v, \, w\rangle _\ast = \int\limits
_{\mathbb{R}}v(y)\overline{w(-y)}\, \mathrm{d} y.
 $}
\end{equation}
Since ${\bf B}$ and ${\bf B^\ast}$ have real point spectrum (see
\ef{SpB*}, a full theory to be developed), we may omit the complex
conjugate.

\begin{proposition}
${\bf B^\ast}$ is adjoint to ${\bf B}$ in the indefinite metric of
$\bar{L}^2(\mathbb{R})$, with the indefinite scalar product \eqref{indefprod}.
\end{proposition}

\noindent {\em Proof.} For our operator, taking $v, w \in
C_0^{\infty}(\mathbb{R})$
\begin{equation*}
\begin{split}
\langle {\bf B}v, \, w\rangle _\ast &=
 \mbox{$
\int\limits _{\mathbb{R}}
 $}
\big((-1)^{k+1}v^{(2k+1)}(y)+\mbox{$\frac{1}{2k+1}$}(yv(y))^{\prime}\big)w(-y)
 \,\mathrm{d} y\\ &=
 \mbox{$
\int\limits_{\mathbb{R}}
 $}
 \mbox{$\big((-1)^{k+1}v^{(2k)}(y)+\frac{1}{2k+1}\,yv\big)w^{\prime}(-y)$}\,\mathrm{d} y\\
%%%%&=\quad \vdots\\
 &=...=
 \mbox{$
\int\limits_{\mathbb{R}}
 $}
 \mbox{$v(y)\big((-1)^{k+1}w^{(2k+1)}(-y)-\frac{1}{2k+1}\,(-y)w'(-y)\big)$}\,\mathrm{d} y
 %\\
 %&
 = \langle v, \, {\bf B^\ast} w\rangle_\ast. \qed
\end{split}
\end{equation*}
%%Hence ${\bf B^\ast}$ is adjoint to the operator ${\bf B}$ in the
%%metric of $\bar{L}^2(\mathbb{R})$.

%%\hfill$\Box$
%%\end{proof}

\smallskip

%%We let our indefinite inner product space be equipped with a
%%decomposition into subspaces
%%\begin{equation}
%%\label{indef} E = E_+\oplus E_-
%%\end{equation}
%%and $E_+\perp E_-$. Here $E_+$ represents the space of odd functions
%%and $E_-$ represents the space of even functions. Then
%%\begin{equation*}
%%v(y)=\mbox{$\frac{v(y)-v(-y)}{2}+\frac{v(y)+v(-y)}{2}$}.
%%%\end{equation*}

Thus, ${\bf B}^*$ is adjoint to ${\bf B}$
%%%in the sense of  (\ref{Badj1})  %%%in the same space $L^2_\rho$
in the given
%% via the skew-symmetric ``dual" product
{indefinite metric}, which we write down again as
\begin{equation}
\label{nn1}
  \mbox{$\langle v,w \rangle_*= \int v(y) \overline{w(-y)}
\, {\mathrm d}y \equiv \langle v,\overline{ J w}\rangle, \quad
 \mbox{for any} \,\,\, v
\in L^2_{\rho}, \,\, Jw \in L^2_{1/\rho}.
    $}
\end{equation}
%%The spectrum of ${\bf B}^*$ is real, so we  omit the complex conjugation and use the space over the
%%field of real numbers.
Here, the {\em canonical symmetry operator} $J w(y)= w(-y)$ is
bounded, self-adjoint, and unitary (it is the {\em Gramm operator}
of this metric).
 Moreover, one can see that, for the given anisotropic weight
 $\rho(y)$ in \ef{H03},
  $$
Jw \in L^2_{1/\rho}, \quad \mbox{if} \quad w \in L^2_\rho,
 $$
 so that the adjoint space in the indefinite metric \ef{nn1} is
  \be
   \label{KK*}
  (L^2_\rho)^*_*= L^2_\rho.
   \ee
%% Note that the simultaneous conditions $w \in L^2_\rho$ and
%%% $J w \in L^2_\rho$ %%%$v_\pm \in L^2_\rho$
%%determine the corresponding  space $L^2_{ \rho^\ast}$

%%Later on, for  some extended linear functionals (polynomials),
 Let us next introduce another natural weighted space,
 with a symmetric exponentially decaying
weight,
\begin{equation}
 \label{rho*}
  \tex{
  L^2_{\rho^*}(\re) \whereA
\rho^*(y) = \frac 1{\hat \rho(y)}= {\mathrm e}^{-a|y|^{\alpha}}
\quad \text{in} \quad \re, \,\,\, \mbox{so} \,\,\, (L^2_{\hat
\rho})^*_{L^2}= L^2_{\rho^*}.
 }
%% |y|
%%\ge 1.
\end{equation}
Note that the standard ``adjoint relation" of the weights
 $
 \tex{
  \rho^*(y)= \frac 1{\rho(y)}
  }
  $
  exists for $y \ll -1$ only and is wrong for $y \gg 1$. Indeed,
  this reflects a strong anisotropic  behaviour of the
  eigenfunctions $\psi_\b(y)$ of $\BB$ as $ y \to \pm \iy$.

The set of even functions $E_+=\{v(-y) \equiv v(y)\}$ is a {\em
positive  lineal} (a linear manifold) of the metric \eqref{nn1},
\begin{equation*}
\langle v,v \rangle_* > 0 \quad \mbox{for} \quad v \in E_+ \subset
L^2_{\hat \rho}, \quad v \not = 0,
\end{equation*}
and odd functions $E_-=\{v(-y) \equiv -v(y)\} \subset L^2_{\hat
\rho}$ give the corresponding {\em negative lineal}.
%%%Since, obviously, each function $w(y)$ admits  the expansion
Therefore, $L^2_{\hat \rho}$ with this metric is {\em
decomposable}:
\begin{equation*}
\mbox{$ v=v_+ + v_- \equiv \frac{v(y)+v(-y)}{2} +
\frac{v(y)-v(-y)}{2}, \quad \text{where} \quad v_\pm \in E_\pm \,
\Longrightarrow \, L^2_{\hat \rho} = E_+ \oplus E_-, $}
\end{equation*}
where, in addition, $ E_+ \,\bot \, E_-$ in the metric \eqref{nn1}.
The corresponding positive {\em majorizing} metric is given by
\begin{equation*}
|\langle v,v\rangle_*| \le [v,v]_* \equiv  \langle
v_+,v_+\rangle_\ast -  \langle v_-,v_-\rangle_*,
\end{equation*}
etc. This case of the decomposable space with an indefinite metric
having a straightforward majorizing one is treated as rather
trivial; see Azizov--Iokhvidov \cite{linopindefmet} for  linear
operators theory in spaces with indefinite metric.  Metric
\eqref{nn1} is widely used therein; see \cite[p.~13, 17, 23,
114]{linopindefmet}. Then the domain of $ \BB^\ast$ is defined as
$H^3_{ \rho^\ast}$, etc.

\smallskip

\noindent{{\bf Historical Remark:}} As we mentioned, basic results
of linear operator theory in spaces with indefinite metrics can be
found in Azizov and Iokhvidov \cite{linopindefmet}. It was not
until about 1944 that L.S. Pontryagin published the article on
``Hermitian operators in spaces with indefinite metric"
\cite{HermOpIndefMet}. A new area of operator theory had been
formed from Pontryagin's studies, which, during the time of the
WWII,
 were originated and associated  with some missile-type military research, \cite{Led08}.
 This work set by
Pontryagin was continued from 1949 and in the 1950s by M.G.~Krein
\cite{HelCurv} and I.S.~Iokhvidov \cite{UnitOpIndef}.
%%%\\[5mm]

%%%%%%%%%%%%%%%%%%%%%%%%%%%%%%%%%%%%%%%%%%%%%%%%%
%%%%%%%%%%%%%%%%%%%%%%%%%%%%%%%%%%%%%%%%%%%%%%%%%
%%%%%%%%%%%%%%%%%%%%%%%%%%%%%%%%%%%%%%%%%%%%%%%%%
%%%%%%%%%%%%%%%%%%%%%%%%%%%%%%%%%%%%%%%%%%%%%%%%%
%%%%%%%%%%%%%%%%%%%%%%%%%%%%%%%%%%%%%%%%%%%%%%%%%

%%%%%%%%%%%%%%%%%%%%%%%%%%%%%%%%%%%%%%%%%%%%%%%%%%%%%%%%%%
\subsection{Discrete spectrum and polynomial eigenfunctions of ${\bf B}^\ast$}

\noindent We consider the spectrum of the linear adjoint operator
${\bf B^\ast}$ in the weighted space $L^2_{\rho^\ast}(\mathbb{R})$,
with exponentially decaying weight
 \eqref{rho*}.
%%, which has been already introduced above,
%%\begin{equation*} \rho^\ast(y) = \mathrm{e}^{-a|y|^{\alpha}}
%%\quad \text{for $|y|\geq 1$}.
%%\end{equation*}
Here we have that $\rho^\ast (y)>0$ and $a\in (0, 2 \hat d_k)$
(see \ef{yminNN})   is a sufficiently small constant. The proof of
the following results does not differ from that of Lemma
\ref{bddop}.

\begin{lemma}
 \label{L.B*1}
${\bf B^\ast}$ is a bounded linear operator from
$H^{2k+1}_{\rho^\ast}(\mathbb{R})$ to $L^2_{\rho^\ast}(\mathbb{R})$.
\end{lemma}

%%%%%%%%%%%%%%%%%%%%%%%%%%%%%%%%%%%%%%%%%%%%%%
%%%%%%%%%%%%%%%%%%%%%%%%%%%%%%%%%%%%%%%%%%%%%%
%%%%%%%%%%%%%%%%%%%%%%%%%%%%%%%%%%%%%%%%%%%%%%
%%%%%%%%%%%%%%%%%%%%%%%%%%%%%%%%%%%%%%%%%%%%%%
%%%%%%%%%%%%%%%%%%%%%%%%%%%%%%%%%%%%%%%%%%%%%%
%%%%%%%%%%%%%%%%%%%%%%%%%%%%%%%%%%%%%%%%%%%%%%

As we have mentioned before, we now, in an easier and direct
manner, construct polynomial eigenfunctions of $\BB^*$.

%%%%%%%%%%%%%%%%%%%%%%%%%%%%
 \begin{lemma}
  \label{L.B*2}
 \noi {\rm (i)} The operator $(\ref{adjoint})$ in
$L^2_{\rho^*}$ has the same point spectrum $\s(\BB^*)=\s(\BB)$ as
in $(\ref{Eigs})$ and the eigenfunctions $\psi^*(y)$ are
generalized Hermite polynomials given by
\begin{equation}
\label{poleigfunc}
 \mbox{$
  \mbox{$\psi ^{\ast}_{\beta}(y) =
\frac{1}{\sqrt{\beta !}}$}\Big[\mbox{$y^{\beta } +$}
(-1)^{k+1}\sum\limits _{j=1}^{\lfloor \frac{\beta
}{2k+1}\rfloor}\mbox{$\frac{1}{j!}D^{(2k+1)j}y^{\beta }$}\Big].
 $}
\end{equation}

 \noi{\rm (ii)} The eigenfunction subset $\Phi^*=\{\psi_\b^*\}$ is
 complete in $L^2$ and in $L^2_{\rho^*}$.

 \end{lemma}

\noi{\em Proof.} (i) The construction is straightforward using the
fact that the operator $y D_y$ is homogeneous. Consider the
eigenfunction equation
\begin{equation}
\label{eigvalprob} {\bf B^\ast } \psi_\b^* = \lambda \psi_\b^*.
\end{equation}
%% Therefore, by
%%distribution theory (see e.g.  Vladimirov \cite[Sect.~8]{Vlad72}),
%%any solution $u(y)$ must be a polynomial. If its degree is $l$,
 Assume that $\psi_\b^*(y)$ is a polynomial of degree $l \ge 0$.
 Then we look for a solution of \ef{eigvalprob} in the form
\begin{equation*}
 \mbox{$
\psi^*_\b(y) = \sum ^s_{j=0}P_j(y),
 $}
\end{equation*}
where each $P_j(y)$ is a homogeneous polynomial and
$s=\big\lfloor\frac{l}{2k+1}\big\rfloor$. From the eigenvalue
problem \eqref{eigvalprob}, we can work out all terms of the
polynomial, for a given degree $l$.  Then
\begin{equation}
\label{ptspec}
 \mbox{$(-1)^{k+1}D_y^{2k+1} \psi_\b^* - \frac{1}{2k+1}\,y D_y
\psi_\b^* = \lambda _\b \psi_\b^*$}
 \quad \Longrightarrow \quad
 \mbox{$\lambda _\b = -\frac{l}{2k+1}$}, \quad
l=0,1,2...\, ,
\end{equation}
%%Then
%%\begin{equation}
%%\label{ptspec} \mbox{$\lambda _l = -\frac{l}{2k+1}$}, \quad
%%k=0,1,2...
%%\end{equation}
and we can define all other polynomials $P_j(y)$ by
\begin{equation}
 \label{PP1}
\mbox{$P_j(y) = \frac{(-1)^{(k+1)}}{j!}\, D^{(2k+1)j}P_0(y), \quad
j=1,...,\,s$}.
\end{equation}
Fixing the leading term $P_0(y)=y^{\beta }$ calculating all the
terms via \ef{PP1} yields \ef{poleigfunc}.

\noi (ii) The proof of completeness is much simpler than that in
Lemma \ref{L.Compl}(ii). In fact, it is a common fact in
functional analysis that polynomials are complete in any
reasonable weighted spaces; see the famous text-book
\cite[p.~431]{KolmF}. $\qed$

\ssk

As in \ef{cl1}, we next define the space of $\Phi^*$-closure
$\tilde L^2_{\rho^*}(\re)$ with the corresponding property
\ef{cl2}. Since, obviously, the resolvent now acts as follows:
 \be
 \label{Res44}
  \tex{
   \mbox{for $v= \sum c_\b \psi_\b^*$, \,\,\, $(\BB^*-\l I)^{-1}v=
    \sum \frac {c_\b}{\l_\b-\l} \, \psi_\b^*$},
     }
 \ee
 i.e., the convergence of the series is improved, it is not
 difficult to conclude that the resolvent is a compact operator in
 $\tilde L^2_{\rho^*}$, though we are not going to use such facts.
 By construction, \ef{Eigs} is the full spectrum of
 $\BB^*$ in $\tilde L^2_{\rho^*}$, with the domain $\tilde {\mathcal D}
  (\BB^*)=\tilde L^2_{\rho^*}(\re) \cap H^{2k+1}_{\rho^*}(\re)$.

We need a few more properties. Firstly, we note that, as in
\cite[\S~5.2]{2mSturm},
 %%\begin{proposition}
%%\label{Pr.L3}
 given an arbitrarily small $\e>0$, there holds
 \be
\label{c152}
 \mbox{$
c_\b = o(|\b|^{-|\b|(\nu +\e)}) \,\, \Longrightarrow \,\, \sum
c_\b \psi_\b^* \in \tilde L^2_{\rho^*}.
 $}
 \ee
%%\end{proposition}
%%%One can see from such estimates that
 %%It also follows that
In addition, if $v= \sum c_\b \psi_\b^* \in \tilde L^2_{\rho^*}$,
then for arbitrarily small fixed $\e>0$,
 \be
\label{c142} \mbox{$ c_\b = o(|\b|^{-|\b|(\nu -\e)}) \quad
\mbox{for ``almost all"} \,\,\, |\b| \gg 1 \quad \bigl(\nu =
\frac{2-\a}{2 \a}= \frac{2k-1}{2(2k+1)}\bigr).
 $}
 \ee

%% \smallskip

%%By $\tilde H^{2k+1}_{\rho^*}$, we then denote the closure in the
%%norm of $H^{2m}_{\rho^*}$ of the linear subspace of eigenfunction
%%expansions with coefficients satisfying (\ref{c142}) for some $\e
%%>0$.  Being equipped with the scalar product (\ref{vnorm}) with $\rho \mapsto \rho^*$,  $\tilde
%%H^{2m}_{\rho^*}$ is a Hilbert space becoming the domain of $\BB^*$
%%in $H^{2m}_{\rho^*}$. We have
%% \be
%%  \label{11.20N}
%%\tilde H^{2m}_{\rho^*} \subseteq H^{2m}_{\rho^*} \cap \tilde
%% L^2_{\rho^*}.
%% \ee

In view of the fast decay (\ref{c152}) of the expansion
coefficients, similar to $l^2_\rho$, we introduce the adjoint
little Hilbert space $l^2_{\rho^*}$ of eigenfunction expansions
$v= \sum c_\b \psi_\b^* \in \tilde L^2_{\rho^*}$ with the scalar
product $(\cdot,\cdot)_{0*}$ and the norm $\|\cdot\|_{0*}$ defined
as in (\ref{vwin}).
 Then the domain of $\BB^*$ in $l^2_{\rho^*}$ is
  the corresponding little Sobolev space $h^{2k+1}_{\rho^*}$
compactly embedded  into $l^2_{\rho^*}$. Here,
$(\cdot,\cdot)_{1*}$ and $\|\cdot\|_{1*}$  denote the scalar
product and the induced norm.
 Then $\BB^*$ is self-adjoint in $l^2_{\rho^*}$, and $\tilde L^2_{\rho^*}$, $\tilde
H^{2k+1}_{\rho^*}$ are dense subspaces of $l^2_{\rho^*}$.

%%%%%%%%%%%%%%%%%%%%%%%%%%%%%%%%%%%%%%%%%%%%%%%%%%%%
\subsection{Bi-orthonormality of the bases $\Phi$ and $\Phi^*$}
 \label{S5.3}

This is a principal issue for several applications of
eigenfunction expansions via eigenfunction sets $\Phi$ and the
adjoint one $\Phi^*$ to be developed later on. First of all, we
note that the generalized Hermite polynomials $\psi_\b^*(y) \in
\Phi^*$ initially  appeared in the semigroup expansions \ef{wsol},
\ef{EigMoments} and \ef{wwNew}, as more or less standard
%%(up to an
 %%indefinite metric to be get rid of below)
 {\em linear functionals}
defined in the weighted space $L^2_{\hat \rho}(\re)$ as given in
\ef{u01}. Therefore, by Hilbert space theory, these polynomials
belong to the true adjoint space: for any $\b$,
 \be
  \label{adj01}
   \tex{
  \psi_\b^* \in L^{2*}_{\hat \rho}(\re)=L^2_{\rho^*}(\re) \whereA
  \rho^*= \frac 1{\hat \rho} \quad (\mbox{cf. (\ref{rho*})).}
 }
   \ee
 Moreover, the whole set $\Phi^*$ of such polynomials is complete
 in $L^2_{\rho^*}$,  and is closed in the corresponding space of
 $\Phi^*$-closure $\tilde L^2_{\rho^*}$.

  However, this is not enough  for further applications.
 Namely, let us  return to the space of $\Phi$-closure $\tilde
L^2_\rho(\re)$ defined for $\BB$ in Section \ref{S4.3}.
 According to the construction therein,
we have defined linear functionals from the adjoint space $\tilde
L^{2*}_\rho$.
  Thus, as
we have seen earlier from the eigenfunction expansion of the
semigroup \ef{wsol} and \ef{BwFull} (for $\t=0$), those linear
functionals are well defined  for data $u_0 \in L^2_{\hat\rho}$
with the standard (not in the {\em v.p.} sense) definition of the
integral; cf. \ef{EigMoments}. Therefore, we then need to extend
those functionals  to the whole space of $\Phi$-closure.

Before doing that, let us get rid of the indefinite metric
$\langle \cdot,\cdot \rangle_\ast$, applying of which could rise
some natural questions, and return to the standard (dual)
$L^2$-one.
 %% However, before we look at this
%%bi-orthonormality condition of the dual space, we first note the
%%following remark: Due to the indefinite metric $\langle
%%\cdot,\cdot \rangle_\ast$,
 Actually, we then need slightly revise our definition
of eigenfunctions $\{\psi_\b^*(y)\}$ of $\BB^*$. Looking back  at
the expansion of the convolution, where the rescaled solution is
given by \eqref{wsol}, we note the definition of the moments of
the initial data, which give rise to the adjoint eigenfunctions,
$\{\psi^\ast_\beta(y)\}$. These moments \eqref{EigMoments}, in the
$L^2$-metric, are assumed to be given by
\begin{equation*}
\mbox{$M_\beta(u_0) \equiv \langle w_0(z), \psi^\ast_\beta (z)
\rangle \equiv \frac{1}{\sqrt{\beta !}}$}\,
 \mbox{$
\int\limits_{\mathbb{R}}
 $}
 \mbox{$z^\beta u_0(z) \, \mathrm{d}z$} \quad \mbox{for any $\b$}.
\end{equation*}
However, we must have the representation in the indefinite metric
that actually yields
\begin{equation*}
%%\begin{split}
\mbox{$M_\beta(u_0) \equiv \langle w_0(z), \psi^\ast_\beta (z)
\rangle_\ast$}  \equiv \mbox{$ \frac{1}{\sqrt{\beta !}}$}\,
 \mbox{$
\int\limits_{\mathbb{R}}
 $}
 \mbox{$(-z)^\beta u_0(z) \, \mathrm{d}z$} =
\mbox{$\frac{(-1)^\beta}{\sqrt{\beta !}}$}\,
 \mbox{$
\int\limits_{\mathbb{R}}
 $}
 \mbox{$z^\beta u_0(z) \, \mathrm{d}z$}.
%%\end{split}
\end{equation*}
Therefore, in the expansion there must be an extra multiplier
$(-1)^\beta$ for the dual products. So, for purposes of
convenience, instead of \ef{poleigfunc},  we now have to take  the
adjoint eigenfunctions $\psi_\b^*$ in the following form (i.e.,
the extra multiplier $(-1)^{|\beta|}$ is included):
\begin{equation}
\label{RevPsi}
%%\begin{equation}
%%\label{poleigfunc}
 \mbox{$
  \mbox{$\hat\psi^{\ast}_{\beta}(y) =
\frac{(-1)^\b}{\sqrt{\beta !}}$}\Big[\mbox{$y^{\beta } +$}
(-1)^{(k+1)}\sum\limits _{j=1}^{\lfloor \frac{\beta
}{2k+1}\rfloor}\mbox{$\frac{1}{j!}D^{(2k+1)j}y^{\beta }$}\Big].
 $}
\end{equation}
 %% \mbox{$\psi_\beta(y) = \frac{1}{\sqrt{\beta
%%!}}\,D^\beta_y F(y)$}, \quad |\beta| \ge 0.
%%\end{equation}

 Thus,
as the next step, according to our construction,
 we first note that, by standard properties of Hilbert spaces, the space
 adjoint to $L^2_\rho(\re)$ is
\be
 \label{adj1}
 \tex{
 %%\tilde
L^{2*}_\rho(\re) = L^2_{\bar \rho}(\re) \whereA \bar \rho = \frac
1 \rho \not = \rho^*  \LongA \psi_\b^* \not \in L^{2*}_\rho(\re).
 }
  \ee
   Similarly, in the metric \ef{nn1}, for which \ef{KK*} holds,
   there is no inclusion of polynomials either.
 In other words, unlike the more standard ``parabolic" case
 \cite{SemiPDE, 2mSturm} (where always $\rho^*= \frac 1 \rho$), the
 polynomials $\psi_\b^*$ do not belong to the adjoint space.
 However, we deal with the space of closure $\tilde L^2_\rho
 \subset L^2_\rho$, with the adjoint one $\tilde L^{2*}_\rho$ being wider
 than the pure and standard  $L^{2*}_\rho$.
 %%We stop at this moment discussing such intriguing aspects of
 %%functional spaces involved that do not play a role later on, and
 Finally, completing this discussion of such intriguing aspects of
 functional spaces involved, which will not play a role later on, we
  mention that
one needs to extend and define such ``extended" linear functionals
 $\{\psi_\b^*\}$ from the adjoint space $L^{2*}_\rho$ to the space
 $L^2_{\rho^*}$, which are defined in
 %%%%%with values in
 the whole $\tilde L^2_\rho$.
 On one hand, this looks related
 %% well corresponds
 %%% In order to deal later on with eigenfunction expansions for both
 %$%%$\BB$ and $\BB^*$, one needs to introduce a formalism of a weak
%%topology in the adjoint space.  In functional analysis, the
 to a standard
procedure of extension of  uniformly convex functional in linear
normed spaces
%%onto $\LL$
%%% are expressed
 by the  classic
Hahn--Banach theorem\footnote{``If $X$ is a linear normed space,
$L$ is a linear manifold, and $f$ is a linear continuous
functional defined on $L$, then $f$ can be extended to $F$ on $X$
and $\|F\|_X=\|f\|_L$".}, \cite{KolmF}.
%% thought there is a
 %%%difference indeed.

%% (we then mean the
%%%maximal extension). ??? ???? ????
 %%We perform this procedure as follows.
%% Namely, as usual, we consider
%%polynomials $\psi_\b^*$ in \ef{psi**1} as linear functionals
%%naturally defined in the adjoint space $\LLL$. Next, the
%%eigenfunctions $\psi_\b$ in \ef{eigen} will be then treated as
%%linear functionals in the second adjoint space denoted by
%%$(\LL)^{**}$, so that the basic space $\LL$ is not assumed to be
%%reflexive.
%% We will not carefully specify the topology of such spaces, and
%% just mention that related questions of weak and weak-* topologies are
%% classic in text books on functional analysis; see e.g.,
%%   Kolmogorov--Fomin
%% \cite[p.~199]{KolF}.
%% In particular,

On the other hand,
 naturally following our direct definition of the space of closure \ef{cl1},
 which assumes the {\em a priori} knowledge of the expansion
 coefficients $\{c_\b\}$
 %%%of any its element,
   for any $v \in \tilde L_\rho^2$, such
extended linear functionals for every $\b$ act as follows:
 \be
 \label{MEAN1}
 \mbox{
 $\langle v, \psi_\b^* \rangle_* \equiv \langle
v, \hat \psi_\b^* \rangle$\, denote  the expansion coefficient
$c_\b$ of $v$ in (\ref{cl2}).
 }
 \ee
 In view of the performed construction of
$\tilde L_\rho^2$ via closure of finite sums, it is not difficult
to see that such continuous extended linear functionals are
defined uniquely (by the density of finite sums in (\ref{cl1})),
and in $L^2_{\hat \rho}$ are given by the standard integrals as in
\ef{EigMoments} and \ef{wwNew}. Defining a countable and complete
set $\Phi^*$ of such extended linear functionals
 suffices for our applications.
 %% so we stop at this moment studying
 %%%the adjoint space.

Thus,  overall, in the sense of \ef{MEAN1},
%% this makes obvious
 the standard bi-orthonormality of the bases $\{\psi_\b\}$ and $\{\psi^*_\b\}$ becomes trivial:
%%%% of these bases:
 \be
 \label{bi1}
%% Then the orthonormality condition holds
 %%$$
 \langle \psi_\b,  \psi_\g^* \rangle_* \equiv
\langle \psi_\b, \hat \psi_\g^* \rangle = \d_{\b\g} \quad
\mbox{for any (multi-indices)} \,\,\, \b \,\, \mbox{and} \,\, \g,
 \ee
 where $\langle \cdot,\cdot \rangle$  is the usual duality
product in $L^2(\re)$ and $\d_{\b\g}$ is the Kronecker delta.
%%%This defines
%%% such functionals $\langle v, \psi_\b^* \rangle$
%%%for any $v \in \tilde L^2_\rho$. Note that by the density of
 %%$L^2_{\rho^*}$, the functional extension used in \ef{bi1} is
 %%uniquely defined.
%% Similarly, using the subset $\Phi^*=\{\psi_\b^*\}$
%%of the generalized Hermite polynomials \ef{poleigfunc} (or
%%\ef{RevPsi} to drop the indefinite metric),
%% from
%%%Proposition \ref
%%  we are obliged to define
%% the corresponding space of closure $\tilde L^{2}_{\rho^*}$
%%  of eigenfunction expansions. Then, eventually, we  treat
Then, we treat similarly
   the adjoint linear functionals
$\langle w, \psi_\b \rangle_*$ for any $w \in \tilde
L_{\rho^*}^{2}$, where $\Phi=\{\psi_\b\}$ is now a countable
(complete) set of linear functionals, which have been properly
extended from the adjoint space $L^{2*}_{\rho^*}(\re)= L^2_{\bar
\rho^*}(\re)$, with $\bar \rho^*= \frac 1{\rho^*}= \hat \rho$.
 %%\equiv \langle v, \bar\psi_\b^* \rangle$
  %% In addition, note that, as usually can
 %%happen,
%%an alternative, but
%%seems equivalent, way of representing such results,
 %%if $\tilde L^2_\rho$ is non-reflexive,
 %%i.e.,
%% \be
%% \label{bi22}
%% (\tilde L^2_\rho)^{**} \not = \tilde L^2_\rho,
%%  \ee
%%  the bi-orthonormality property \ef{bi1} can be interpreted as values of the
%%  linear functionals $\psi_\b \in (\tilde L^2_\rho)^{**}$ on the elements $\psi_\g^* \in
%%  \tilde
%%   L^{2}_{\rho,*}$. This requires extra study of those spaces,
%% which is too detailed for the present need of applications, so
%%    we will not use such an issue in what follows.

%% (we then mean the
%%%maximal extension). ??? ???? ????
 %%We perform this procedure as follows.
%% Namely, as usual, we consider
%%polynomials $\psi_\b^*$ in \ef{psi**1} as linear functionals
%%naturally defined in the adjoint space $\LLL$. Next, the
%%eigenfunctions $\psi_\b$ in \ef{eigen} will be then treated as
%%linear functionals in the second adjoint space denoted by
%%$(\LL)^{**}$, so that the basic space $\LL$ is not assumed to be
%%reflexive.
%% We will not carefully specify the topology of such spaces, and
%% just mention that related questions of weak and weak-* topologies are
%% classic in text books on functional analysis; see e.g.,
%%   Kolmogorov--Fomin
%% \cite[p.~199]{KolF}.
%% In particular,

\ssk

According to the above analysis, we now show some convenient
formal calculus with non-convergent integrals, which illustrate
standard bi-orthonormality and other properties. For instance,
according to these formal integration techniques, we have:

\begin{proposition}
\label{Ortho1}
\begin{equation*}
\langle \psi_\beta, \psi_\beta^\ast \rangle_\ast = 1 \quad
\text{for all} \quad \beta\geq 0.
\end{equation*}
\end{proposition}

\noindent {\em Proof.} We start by looking at the scalar product,
in the indefinite metric, defined by
\begin{equation*}
\mbox{$\langle \psi_\beta, \psi^\ast_\beta \rangle_\ast = -\int
\psi_\beta(y) \psi^\ast_\beta(-y)\, \mathrm{d}y$}.
\end{equation*}
By our definitions of $\psi_\beta(y)$ \eqref{RevPsi} and
$\psi^\ast_\beta(y)$ \eqref{poleigfunc}, we substitute to find that
\begin{equation*}
 \mbox{$
\langle \psi_\beta, \psi^\ast_\beta \rangle_\ast = -\int
\mbox{$\frac{1}{\sqrt{\beta!}}\,D^\beta_y F(y)$}\,
\mbox{$\frac{1}{\sqrt{\beta !}}$}\Big[\mbox{$(-y)^{\beta } +$}
(-1)^{(k+1)}\sum\limits _{j=1}^{\lfloor \frac{\beta
}{2k+1}\rfloor}\mbox{$\frac{1}{j!}D^{(2k+1)j}(-y)^{\beta }$}\Big]
\,\mathrm{d}y.
 $}
\end{equation*}
We now apply to $\psi^\ast_\beta$ the identity operator
 %%\begin{equation*}
$(D_y^{\beta})^{-1}D_y^{\beta} = I$,
 %%%\end{equation*}
with a standard definition and construction of the inverse
integral operator $(D_y^{\beta})^{-1}=D_y^{-\beta}$,  such that
\begin{equation*}
\mbox{$\langle \psi_\beta , \psi^\ast_\beta
\rangle_\ast$}=\mbox{$\langle \psi_\beta , (D_y^{-\beta}\,
D_y^{\beta}) \psi^\ast_\beta \rangle_\ast$}.
\end{equation*}
Integrating by parts, $\beta$ times, we find that
\begin{equation*}
\begin{split}
\mbox{$\langle \psi_\beta , \psi^\ast_\beta \rangle_\ast$} &=
 \mbox{$
\int
 $}
\mbox{$\frac{1}{\sqrt{\beta!}}\,D^\beta_y F(y)(D^{-\beta}_y\,
D_y^\beta)\, \frac{1}{\sqrt{\beta !}}$}\, \mbox{$\big((-y)^{\beta
} + o(y^\beta)\big)\,\mathrm{d}y$}\\ &=
%%\mbox{$\frac{1}{\beta!}$}\,
%% \mbox{$
%%\int
%% $}
%%  -\big[ \int\mbox{$ D^\beta_y F(y)\,$}\mathrm{d}y\,
%%D_y(D^{\beta}_y)^{-1}D_y^\beta\, \mbox{$\big((-y)^{\beta } +
%%o(y^\beta)\big)$}\big]\,\mathrm{d}y\\ &=
\mbox{$\frac{1}{\beta!}$}\,
 \mbox{$
\int
 $}
 (-1)\big[ \mbox{$ D^{\beta-1}_y F(y)$}\,
D^{-\beta+1}_y \,D_y^\beta\, \mbox{$\big((-y)^{\beta } +
o(y^\beta)\big)$}\big]\,\mathrm{d}y
%%%%\\ &
%%\vdots
%%%\cdots
\\ &=...= \mbox{$\frac{1}{\beta!}$}\,
 \mbox{$
\int
 $}
 (-1)^\beta\,\big[ \mbox{$ F(y)$}\,D_y^\beta\,
\mbox{$\big((-y)^{\beta } + o(y^\beta)\big)$}\big]\,\mathrm{d}y\\
&= \mbox{$\frac{1}{\beta!}$}
 \mbox{$
\int
 $}
  (-1)^\beta\,\mbox{$ F(y)\,(-1)^\beta\beta!$}\,\mathrm{d}y.
\end{split}
\end{equation*}
One can see that, according to these formal calculus of
integration by parts, we each time improve the convergence
properties of the integrals involved, meaning using a
distributional treatment of those integrals as values of certain
linear functionals as generalized functions (distributions). As
customary, this corresponds to a regularization of divergent
integrals.

 Hence, eventually,  it follows that
\begin{equation*}
 \mbox{$
\mbox{$\langle \psi_\beta , \psi^\ast_\beta \rangle_\ast$} = \int
F(y)\, \mathrm{d}y = 1,
 $}
\end{equation*}
 so finally we arrive at a convergent integral, but not absolutely, since
 $F$ is not Lebesgue measurable in $\re$, as we have seen.
 $\qed$

%%\hfill$\Box$
%%\end{proof}

\begin{proposition}
\label{Ortho2}
   In terms of the above formal calculus,
\begin{equation}
%\label{betagamma0}
 \label{ort77}
\langle \psi_\beta, \psi_\gamma^\ast \rangle_\ast = 0 \quad
\text{for all} \quad \beta\neq \gamma.
\end{equation}
\end{proposition}

\noindent {\em Proof.} The first part of the proof follows that of
Proposition \ref{Ortho1}. First, consider the case where
$\beta>\gamma$.  After integration by parts, it can be seen that
$\langle \psi_\beta, \psi_\gamma^\ast \rangle_\ast$ may be written
as
\begin{equation*}
  \mbox{$
\mbox{$\langle \psi_\beta , \psi^\ast_\gamma \rangle_\ast$}
=(-1)^\beta \int \mbox{$\frac{(-1)^\beta}{\sqrt{\beta!}}
(D^{\beta}_y)^{-1} D^\beta_y F(y)\frac{1}{\sqrt{\gamma
!}}$}\,D_y^\beta \mbox{$\big((-y)^{\gamma } +...
\big)\,\mathrm{d}y$}.
 $}
\end{equation*}
However, since $\beta>\gamma$, then it is known that
%%\begin{equation*}
$D_y^\beta \mbox{$(-y)^{\gamma }$} = 0$.
%%\end{equation*}
So it follows that \eqref{ort77} holds.

Now consider the case when $\beta < \gamma$.  In this case it can
easily be seen that the above argument will not work.  Rather than
attempting to use a similar argument, we instead use another proof
which encompasses both cases of $\beta>\gamma$ and $\beta<\gamma$.

By the definitions of the linear operator ${\bf B}$ and the adjoint
operator ${\bf B}^\ast$, we know that
\begin{equation}
\label{eigvalsdualprob}
\begin{cases}
{\bf B}\psi_\beta = \lambda_\beta\psi_\beta,\\
{\bf B^\ast}\psi_\gamma^\ast = \lambda_\gamma\psi_\gamma^\ast,
\end{cases}
\end{equation}
which defines the eigenvalue problems for these two operators.
Taking the inner product of \eqref{eigvalsdualprob}, in the
indefinite metric, with $\psi_\gamma^\ast$ and $\psi_\beta$,
respectively, yields
\begin{equation*}
\begin{cases}
\langle {\bf B}\psi_\beta, \psi_\gamma^\ast \rangle_\ast = \lambda_\beta\langle \psi_\beta, \psi_\gamma^\ast \rangle_\ast, \\
\langle \psi_\beta, {\bf B^\ast}\psi_\gamma^\ast \rangle_\ast =
\lambda_\gamma\langle \psi_\beta, \psi_\gamma^\ast \rangle_\ast.
\end{cases}
\end{equation*}
However from the definition of the adjoint operator ${\bf B}^\ast$,
we know that
\begin{equation*}
\langle {\bf B}\psi_\beta, \psi_\gamma^\ast \rangle_\ast = \langle
\psi_\beta, {\bf B^\ast}\psi_\gamma^\ast \rangle_\ast.
\end{equation*}
Hence, if $\beta\neq \gamma$ (i.e., $\lambda_\beta\neq
\lambda_\gamma$), (\ref{ort77}) follows.
%% it must follow that
%%\begin{equation*}
%%\langle \psi_\beta, \psi_\gamma^\ast \rangle_\ast = 0 \quad
%%\text{for all} \quad \beta\neq \gamma.
%%\end{equation*}
 %%Again, when necessary, we assume a proper regularisation of the
 %%%integrals, which are treated as values of some linear functionals.
 $\qed$

%%\hfill$\Box$
%%\end{proof}

%%\begin{corollary}
%%The orthonormality condition
%%\begin{equation*}
%%\langle \psi_\beta, \psi_\gamma^\ast \rangle_\ast =
%%\delta_{\beta\gamma}
%%\end{equation*}
%%holds true, in the indefinite metric, for all $\beta,\gamma \ge
%%0$.
%%\end{corollary}
 %%This follows directly from Propositions \ref{Ortho1} and
%%%\ref{Ortho1}.

%%%%%%%%%%%%%%%%%%%%%%%%%%%%%%%%%%%%%%%%%%%%%%%%%%%%%%%%%%%%%%%%
%%%%%%%%%%%%%%%%%%%%%%%%%%%%%%%%%%%%%%%%%%%%%%%%%%%%%%%%%%%%%%%%
%%%%%%%%%%%%%%%%%%%%%%%%%%%%%%%%%%%%%%%%%%%%%%%%%%%%%%%%%%%%%%%%
%%%%%%%%%%%%%%%%%%%%%%%%%%%%%%%%%%%%%%%%%%%%%%%%%%%%%%%%%%%%%%%%
%%%%%%%%%%%%%%%%%%%%%%%%%%%%%%%%%%%%%%%%%%%%%%%%%%%%%%%%%%%%%%%%
%%%%%%%%%%%%%%%%%%%%%%%%%%%%%%%%%%%%%%%%%%%%%%%%%%%%%%%%%%%%%%%%
%%%%%%%%%%%%%%%%%%%%%%%%%%%%%%%%%%%%%%%%%%%%%%%%%%%%%%%%%%%%%%%%
%%%%%%%%%%%%%%%%%%%%%%%%%%%%%%%%%%%%%%%%%%%%%%%%%%%%%%%%%%%%%%%%

%%%%%%%%%%%%%%%%%%%%%%%%%%%%%%%%%%%%%%%%%%
%%%%%%%%%%%%%%%%%%%%%%%%%%%%%%%%%%%%%%%%%%
%%%%%%%%%%%%%%%%%%%%%%%%%%%%%%%%%%%%%%%%%%
%%%%%%%%%%%%%%%%%%%%%%%%%%%%%%%%%%%%%%%%%%
%%%%%%%%%%%%%%%%%%%%%%%%%%%%%%%%%%%%%%%%%%
%%%%%%%%%%%%%%%%%%%%%%%%%%%%%%%%%%%%%%%%%%%%%%%%%%%%%%%%%%%%%%%%%%%%%
\section{Semilinear dispersion PDEs: VSSs, numerics,  and preliminaries}
 \label{SSemD}

\subsection{Semilinear model}

We now consider the odd-order problem, but now with a nonlinear
absorption. We look at the Cauchy problem for the semilinear
odd-order equation \ef{semilinear pde}, for $k = 1, 2, 3,\hdots$,
%%\begin{equation}
%%\label{semilinear pde} u_t = (-1)^{k+1}D _x^{2k+1}u - u^p \quad
%%\text{in} \quad \mathbb{R}\times\mathbb{R}_+,
%%\end{equation}
with sufficiently good initial data $u(x,0)=u_0(x)$.
%% Here $p>1$ is a fixed absorption exponent.
%% For
%%convenience, we write
%%\begin{equation*}
%%u^p := |u|^{p-1}u
%%%%,\quad \text{with} \quad p>1,
%%\end{equation*}
%%in order to avoid any singularities due to the power $p$, due to the
%%changing sign of the solution $u(x,t)$.
 As we have already commented on, this non-integrable model is somehow connected to the KdV equation
\eqref{KdV}, with the difference being that the extra operator
$-|u|^{p-1}u$, which corresponds to absorption, is simpler and of
zero differential order.
 This kind of nonlinearity also allows us to avoid
 entering remarkable classes of integrable PDEs, which obey various
extremely strong and quite specific properties that are illusive
for more general nonlinear odd-order PDEs.

Let us mention again, as in the linear case, even-order semilinear
problems have been studied over recent years and are fairly well
understood, unlike similar odd-order ones.  The related
generalized semilinear even-order model is given by
\begin{equation}
 \label{Par1}
u_t = -(-\Delta)^m u \pm |u|^{p-1}u,
\end{equation}
which includes both the cases of absorption and reaction.  The
lower order absorption case, with $m=1$, corresponds to the heat
equation with absorption given in \eqref{heat}.

%This case where we have
%****** has blow-up solutions. The lower order case where $m=1$
%corresponds to the heat equation with absorption *********.

%%%%%%%%%%%%%%%%%%%%%%%%%%%%%%%%%%%%%%%%%%%%%%%%%
%%%%%%%%%%%%%%%%%%%%%%%%%%%%%%%%%%%%%%%%%%%%%%%%%
%%%%%%%%%%%%%%%%%%%%%%%%%%%%%%%%%%%%%%%%%%%%%%%%%
%%%%%%%%%%%%%%%%%%%%%%%%%%%%%%%%%%%%%%%%%%%%%%%%%
%%%%%%%%%%%%%%%%%%%%%%%%%%%%%%%%%%%%%%%%%%%%%%%%%
%%%%%%%%%%%%%%%%%%%%%%%%%%%%%%%%%%%%%%%%%%%%%%%%%

%%%%%%%%%%%%%%%%%%%%%%%%%%%%%%%%%%%%%%%%%%%%%%%%%%%%%%%%%%%%%%%%
\subsection{Similarity solutions of semilinear equations}

%We look at papers \cite{SemiPDE}, \cite{VSSSParaPDEs}

As is customary in PDE theory, there exists a critical exponent for
\eqref{semilinear pde} given by
\begin{equation}
\label{fujita}
\mbox{$p = p_0 = 1+ \frac {2k+1}{N}\big|_{N=1} = 2k+2$}. %\big_{N=1}= 2k+2.
\end{equation}
It can be called the critical {\em Fujita} exponent; see further
comments below. As for the parabolic equation \eqref{Par1}, where
\begin{equation*}
\mbox{$ p_0=1+ \frac {2m}{N} $}
\end{equation*}
(as in \eqref{fujita}, $2m$ stands for the order of the differential
operator involved), the critical Fujita exponent characterizes
parameter ranges of blow-up and non-blow-up solutions and changing
of the stability of the trivial zero solutions for the PDEs under
consideration.
%%this means that

As usual, we consider self-similar solutions of the very singular
type of \eqref{semilinear pde}
\begin{equation*}
u_{\ast}(x,t) = t^{-\frac{1}{p-1}}f(y), \quad y =
x/t^{\frac{1}{2k+1}},
\end{equation*}
where $f$ solves the ODE
\begin{equation}
\label{conv}
\mbox{$(-1)^{k+1}f^{(2k+1)}+\frac{1}{2k+1}\,f'y+\frac{1}{p-1}\,f-|f|^{p-1}f=0$}
\quad \text{in} \quad \mathbb{R}.
\end{equation}
 Unlike the linear case, we cannot integrate  to reduce the
order of the ODE and hence we remain with an equation of order
$2k+1$. We see that \eqref{conv} is a difficult higher-order
equation and so begin with numerical results.

%%%%%%%%%%%%%%%%%%%%%%%%%%%%%%%%%%%%%%%%%%
%%%%%%%%%%%%%%%%%%%%%%%%%%%%%%%%%%%%%%%%%%
%%%%%%%%%%%%%%%%%%%%%%%%%%%%%%%%%%%%%%%%%%
%%%%%%%%%%%%%%%%%%%%%%%%%%%%%%%%%%%%%%%%%%
%%%%%%%%%%%%%%%%%%%%%%%%%%%%%%%%%%%%%%%%%%

%%%%%%%%%%%%%%%%%%%%%%%%%%%%%%%%%%%%%%%%
\subsection{Numerical results for the semilinear equation: $k=1$}
\label{NumerSem}

 %%\subsection

%% \noi{\sc First basic numerics.}

We look at similarity profiles of the semilinear equation
\eqref{k1N}, for $k=1$,
\begin{equation}
 \label{fk1}
\mbox{$f^{'''}+\frac{1}{3}\, f'y+\frac{1}{p-1}\,f-|f|^{p-1}f=0$}
\quad \text{in} \quad \mathbb{R}.
\end{equation}
 %%where, as we have said, $f^p=|f|^{p-1}f$.
Once again we use the {\tt MatLab bvp4c} solver to plot the
profiles. This however gives us the same problem, as in Section
\ref{NumConstr}, of solving an initial value problem using the BVP
solver. We take $f$ and $f'$ to be zero, as boundary conditions at
the left-hand boundary point of a fixed interval and $f$ to be
zero as the condition at the right-hand point.
 This corresponds well to a proper BVP setting for odd-order
 operators; see \cite{Fam02} as a guide.
 We then try to shift the point at the right boundary in
order to coincide with a point where the oscillations go through
zero and hence find the best profile.  We look in particular for
convergence and reflectional (about the axis $f=0$) symmetry of
the tail.

\begin{figure}[htbp]
\centering \subfigure[General view of $f(y)$.]{
\includegraphics[scale=0.45]{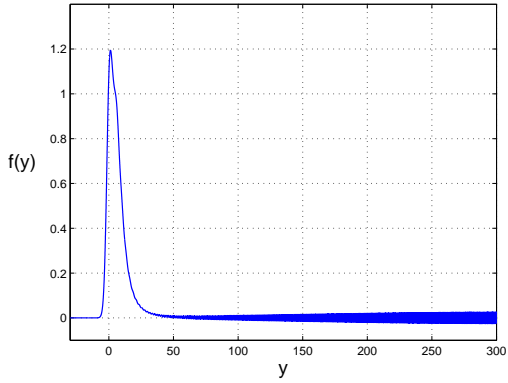}}
\subfigure[``Tail" of the solution, for $y\gg 1$]{
\includegraphics[scale=0.45]{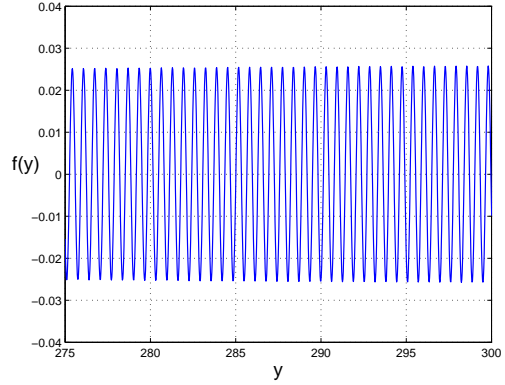}}
\caption{\small A VSS profile of (\ref{fk1})  for $p=1.9$.}
 \label{p19}
\end{figure}

\begin{figure}[htbp]
\centering \subfigure[General view of $f(y)$]{
\includegraphics[scale=0.45]{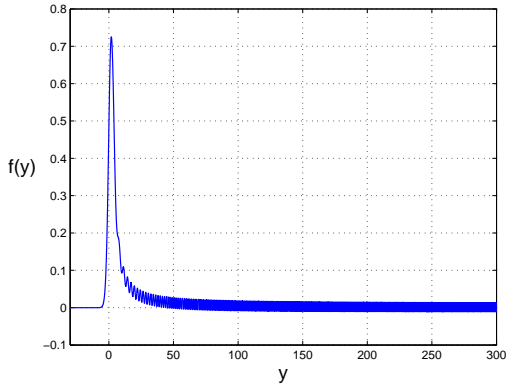}}
\subfigure[``Tail" of the solution, for $y\gg 1$]{
\includegraphics[scale=0.45]{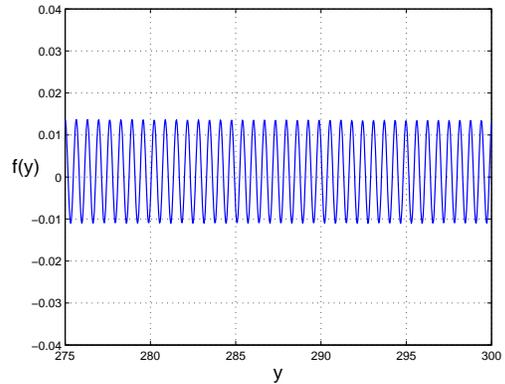}}
\caption{\small A VSS profile of (\ref{fk1}) for $p=3.3$.}
 \label{p33}
\end{figure}

Figures \ref{p19} and \ref{p33} show two sets of profiles for
$p=1.9$ and $p=3.3$.
%% We also refer to Chapter \ref{AppenNumer} for
%%more profiles for values
%%\begin{equation*}
%%1.8 \leq p \leq 3.2.
%%\end{equation*}
Whilst we cannot guarantee the accuracy of the size of the tail, the
profile close to the origin, where $\max|f|$ is obtained, is very
stable. We see that as $p$ decreases, $\max|f|$ increases and this
seems to justify the term ``Very Singular", with the mass
concentrated close to $y=0$ .

In the linear case set out in Section \ref{NumConstr}, due to the
instability of the solutions, we had to use a ``matching
technique'' in order to find reliable profiles.  However, in the
semilinear case, this is not necessary, as the profiles are
stable, since we do not have a scaling group of solutions
$\{cf\}$, due to the nonlinear term $|f|^{p-1}f$. In fact, we find
that whilst the tail of the profile may differ, given different
boundary points, the rest of the structure is extremely stable and
rarely changes. For, approximately,
\begin{equation*}
2\le p \le 3,
\end{equation*}
we can find profiles easily, for almost every boundary point value
we use and they only really differ from the tail.  But for other
values it is more
difficult to find profiles, especially reliable ones.
%%%%\\[3mm]

In particular, we look to see the behaviour of the solutions as $p
\to 1^+$ and as $p \to 4^-$, which is the critical {\em Fujita}
exponent \eqref{fujita} for $k=1$.
%This work is ongoing and we hope to find reliable profiles within the region $1<p<4$.\\[5mm]

\ssk

 %%\subsection

 %% \noi{\sc Extra aspects of finding reliable profiles.}
Since we are using the boundary value solver of {\tt MatLab}, we
face numerous difficulties in obtaining the correct numerical
results. As we are solving the initial value problem of the
semilinear equation as a boundary value problem, difficulties
arise in finding the correct right-hand end-point of the interval.

As we do not know exactly at which points $f$ and $f^\prime$ are
both zero, we attempt to approximate the point by looking at
various profiles and finding the best.  If the far right boundary
point is incorrect, then we produce artificial oscillations.  We
know which profiles are most likely to be false from the analysis
we have done on odd-order linear PDEs. Since the semilinear
equation can be naturally treated as the linear one with a
perturbation  $-|f|^{p-1}f$, we expect similar behaviour.  With a
wrong boundary point, we can end up with non-symmetric tails,
since we are forcing the oscillations through a specific point.
Hence, the most important condition we look for in reliable
profiles are that the oscillations are symmetric.

We also look to ensure that the oscillations become symmetric as
quickly as possible.   This also follows from having forced
oscillations due to wrong boundary points.  Finally, we look for
oscillations that decay as $y$ increases.

The solver that we use, {\tt bvpinit}, requires as the name
suggests an initial guess for the function.  Whilst this guess
does not necessarily have to be very accurate, occasionally a
wrong guess can lead to wrong profiles being found.  We have yet
to find any evidence so far, to suggest that the initial guess is
bad enough to affect the profiles though. However for less stable
profiles where $p\sim2$, the solutions blow-up and the initial
guess for $\max|f|$ may not be accurate. The values needed for the
initial guess for the maximum value are much larger than the
maximum value found in the profile.  Whilst the profiles look
reliable, it is unclear as to whether the inaccurate initial guess
affects the output.
%%\\[5mm]

\ssk

  Thus, the conditions we look at to ensure the best
possible profiles can be briefly summarised as follows:
\begin{itemize}
\item
Symmetry (reflectional, $f\mapsto -f$) of tail for $y\gg1$.
\item
Symmetry of tail occurs as close to $0$ as possible.
\item
Minimisation of symmetric tail.
\end{itemize}

%%%%%%%%%%%%%%%%%%%%%%%%%%%%%%%%%%%%

%%%%%%%%%%%%%%%%%%%%%%%%%%%%%%%%%%%%%%%%%%%%%%
%%%%%%%%%%%%%%%%%%%%%%%%%%%%%%%%%%%%%%%%%%%%%%
%%%%%%%%%%%%%%%%%%%%%%%%%%%%%%%%%%%%%%%%%%%%%%
%%%%%%%%%%%%%%%%%%%%%%%%%%%%%%%%%%%%%%%%%%%%%%
%%%%%%%%%%%%%%%%%%%%%%%%%%%%%%%%%%%%%%%%%%%%%%

\begin{figure}
%  \vskip -.3cm
\centering
\includegraphics[scale=0.6]{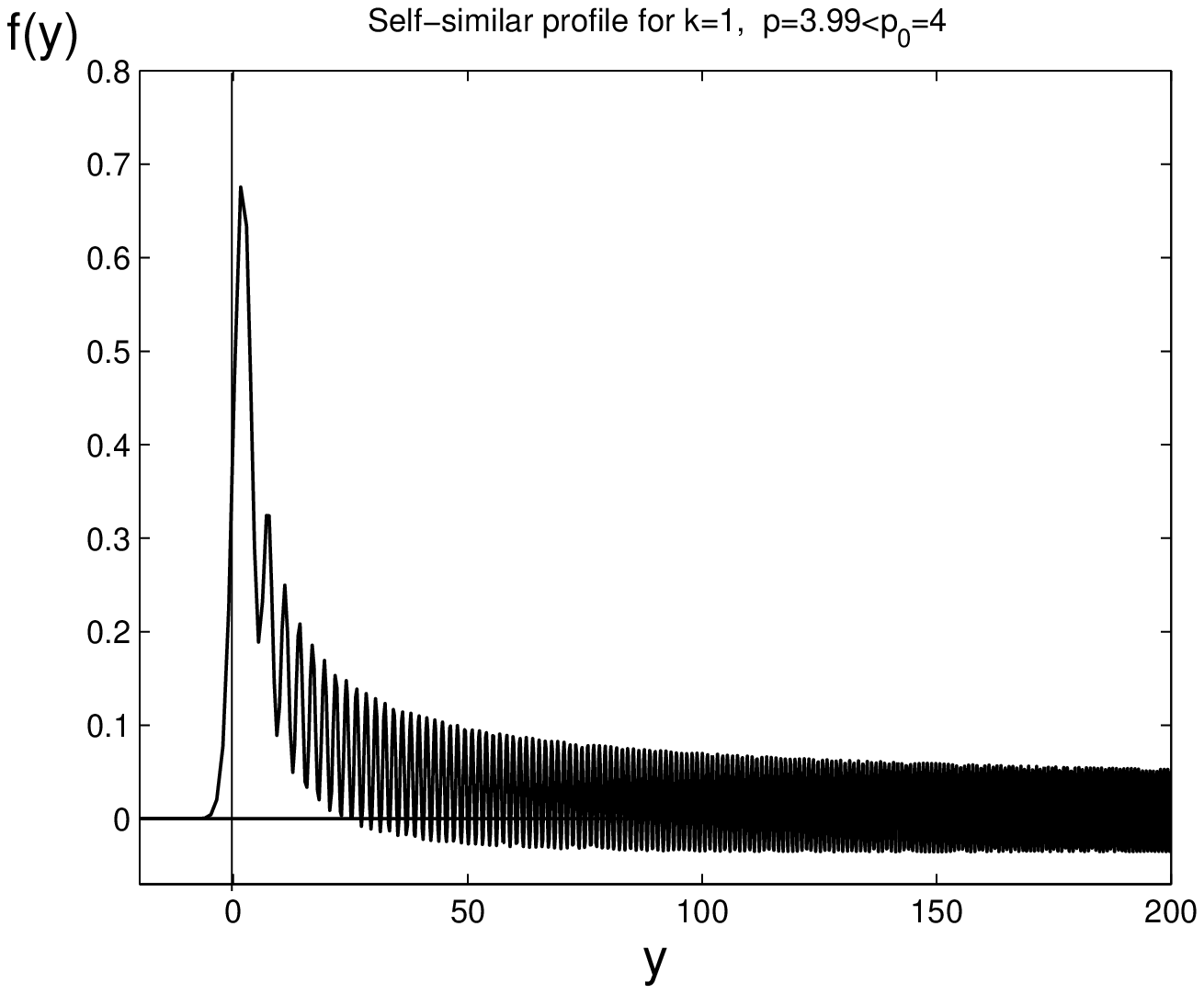}
\vskip -.5cm \caption{\small A VSS profile $f(y)$ of (\ref{fk1})
for $p =3.99<p_0=4$.}
 %%%%%%%%%%%%%%%%%%% Here $z \mapsto -z$.}
   \vskip -.3cm
 \label{Fp3.99}
\end{figure}

Finally, in Figure \ref{Fp3.99} we present the VSS profile for
$p=3.99$, which is sufficiently close to the bifurcation critical
Fujita exponent $p_0=4$; see Section \ref{S7.4} for derivation.

%%FIG%%%%%%%%%%%%%%%%%%%%%%%

%%%%%%%%%%%%%%%%%%%%%%%%%%%%%%%%%%%%%%%%%%%%%%%%%%%%
\subsection{Numerical results: $k \ge 2$}

\begin{figure}
%  \vskip -.3cm
\centering
\includegraphics[scale=0.6]{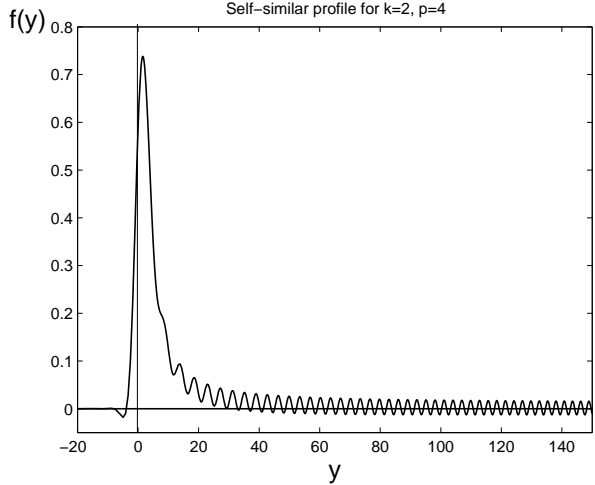}
\vskip -.5cm \caption{\small A VSS profile $f(y)$ of (\ref{fk2})
for $p =4<p_0=6$.}
 %%%%%%%%%%%%%%%%%%% Here $z \mapsto -z$.}
   \vskip -.3cm
 \label{Fpk2}
\end{figure}

\begin{figure}
%\vskip -.3cm
\centering
 % \vskip -.4cm
\subfigure[$p=3$]{
\includegraphics[scale=0.52]{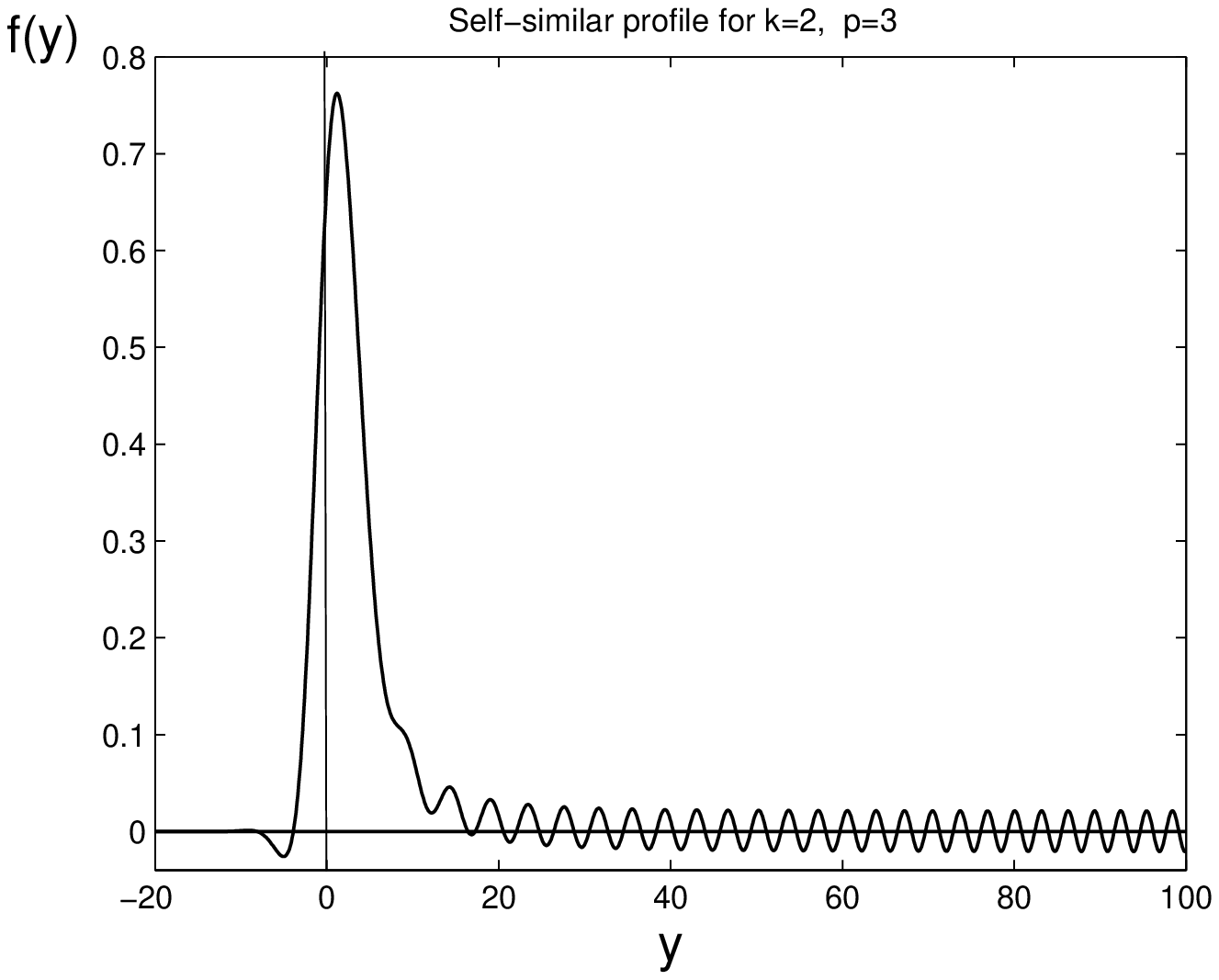}
}
 %\vskip -.4cm
\subfigure[$p=5.9$]{
\includegraphics[scale=0.52]{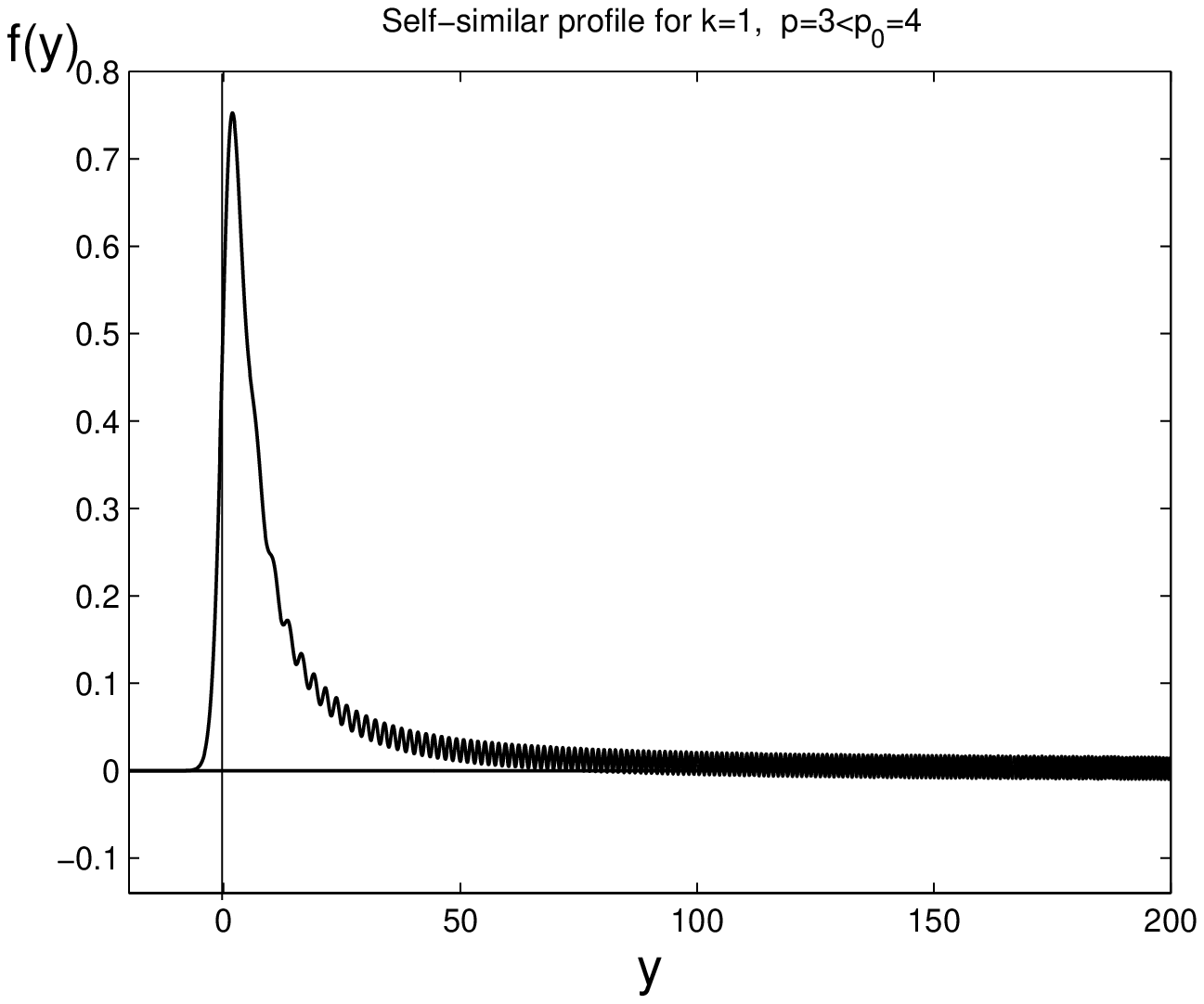}
}
 \vskip -.2cm
\caption{\rm\small VSS profiles $f(y)$ of (\ref{fk2}) for $p=3$
(a) and $p=5.9<p_0=6$ (b).}
 %%\vskip -.2cm
 \label{Fpk21}
\end{figure}

For $k=2$,  scaling (\ref{kern}) in the PDE yields the following
equation for $f$:
 \be
 \label{fk2}
  \mbox{$
 u_t=-u_{xxxxx}-|u|^{p-1}u \LongA -f^{(5)}+ \frac 15\, f'y + \frac
 1{p-1}\,f-|f|^{p-1}f=0 \inB \re.
  $}
  \ee
  In Figure \ref{Fpk2}, we show the similarity solution of
 the ODE \ef{fk2} for $p=4$, where the tail is not that
 oscillatory as used to be for $k=1$.
 Figure \ref{Fpk21} shows VSS profiles for $p=3$ and $p=5.9$,
 which is slightly lower than the critical exponent $p_0=6$ for
 $k=2$.

%%%%%%%%%%%%%%%%%%%%%%%%%%%%%%%%%%%%%%%%%%%%%%%%%%%%%%%

%%%%%%%%%%%%%%%%%%%%%%%%%%%%%%%%%%%%%%%%%%%%%%%%%%%
%%%%%%%%%%%%%%%%%%%%%%%%%%%%%%%%%%%%%%%

%%%%%%%%%%%%%%%%%%%%%%%%%%%%%%%%%%%%%%%%%%%%%%%%%%%%%%%%%

\ssk

\begin{figure}
%  \vskip -.3cm
\centering
\includegraphics[scale=0.65]{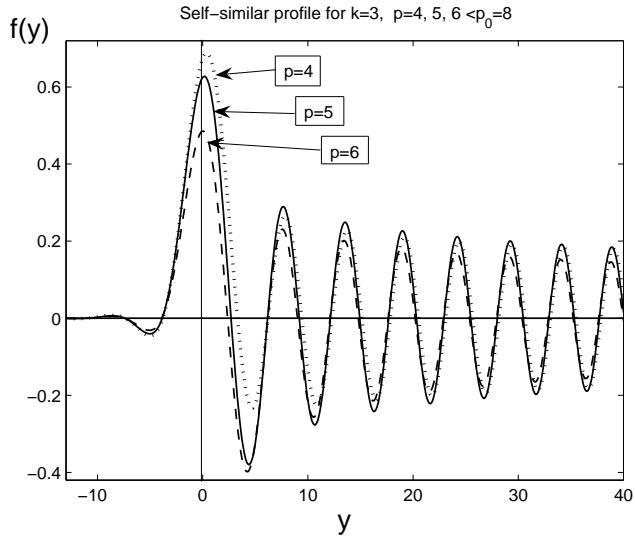}
\vskip -.5cm \caption{\small VSS profiles of the ODE (\ref{fk3})
for $p =4,\, 5,\, 6$.}
 %%%%%%%%%%%%%%%%%%% Here $z \mapsto -z$.}
   \vskip -.3cm
 \label{Fpk3}
\end{figure}

Finally, we present a few VSS profiles for the seventh-order
equation for $k=3$:
\be
 \label{fk3}
  \mbox{$
 u_t=u_{xxxxxxx}-|u|^{p-1}u \LongA f^{(7)}+ \frac 17\, f'y + \frac
 1{p-1}\,f-|f|^{p-1}f=0 \inB \re.
  $}
  \ee
 These are given in Figure \ref{Fpk3} for $p=4, 5,$ and 6. Note
 that the first critical exponent is $p_0=8$ for $k=3$, so these
 values are well below the $p_0$-bifurcation point; see Section
 \ref{S7.4}.

 %%%%%%%%%%%%%%%%%%%%%%%%%%%%%%%%%%%%%%%%%%%%%%%%%%%%%%%

%%%%%%%%%%%%%%%%%%%%%%%%%%%%%%%%%%%%%%%%%%%%%%
 %%%%%\newpage
\section{Linearized stability analysis and centre subspace behaviour:
applications of Hermitian spectral theory}
 \label{SectLinStab}

\subsection{Linearized stability and critical (Fujita) exponent}

We begin with some applications of the spectral analysis obtained
before, in the theory for LDEs.

%%In view of the highly oscillatory character of the eigenfunctions
%%$\{\psi_\beta\}$ of ${\bf B}$, with no exponential decay, some
%%functional formalities of ``dual" space $ \tilde L^2(\mathbb{R})$
%%with the indefinite metric \eqref{indefprod} are not easy at all
%%and still remain obscure. So some of our future conclusions will
%%be formal and we will need clearly indicate which ones are.

As in the linear case, we introduce the following similarity
scaling in the semilinear equation \eqref{semilinear pde}, with
\begin{equation*}
u(x,t) = (1+t)^{-\frac{1}{2k+1}}v(y,\tau), \quad
y=x/(1+t)^{\frac{1}{2k+1}}, \quad \tau = \ln(1+t).
\end{equation*}
Our rescaled equation is then given by
\begin{equation}
%%\begin{split}
\label{rescsem} v_\tau = (-1)^{k+1}D_y^{2k+1}v
+\mbox{$\frac{1}{2k+1}\,yD_yv+\frac{1}{p-1}\,v-|v|^{p-1}v$}
\equiv{\bf B}_1v -|v|^{p-1}v.
%%\end{split}
\end{equation}
Here, the linear operator ${\bf B}_1$ is defined by
\begin{equation*}
{\bf B}_1= {\bf B}+ d_1I, \quad \text{where} \quad
d_1=\mbox{$\frac{1}{p-1}-\frac{1}{2k+1}=\frac{p_0-p}{(2k+1)(p-1)}$},
\end{equation*}
and ${\bf B}$ is our canonical rescaled  one  \ef{linODE}.
%% is from the linear
%%theory in Section \ref{S2.3}.

\begin{lemma}
For $p>p_0=2k+2$, zero is exponentially linearly stable for
\eqref{rescsem} in $H^{2k+1}_\rho(\mathbb{R})\cap \tilde
L^2_\rho(\re)$.
\end{lemma}

\noindent {\em Proof.}  We look at the linearized problem of
\eqref{rescsem} about zero,
\begin{equation*}
v_\tau={\bf B}_1v.
\end{equation*}
Then the spectrum is given by
%%\begin{equation*}
$\mbox{$\sigma({\bf B}_1) = \big\{d_1-\frac{l}{2k+1} \quad
\text{for} \quad l\geq 0 \big\}$}$.
%%\end{equation*}
Therefore, for $l=0$,
 %%%\begin{equation*}
$\lambda_0 = d_1\equiv\mbox{$\frac{1}{p-1}-\frac{1}{2k+1}$} <0$,
for all $ p>2k+2$.
%%\end{equation*}
So, in our space,
\begin{equation*}
\|v(\tau)\|_{2k+1,\rho}\sim \mathrm{e}^{\lambda_0\tau} \to 0 \quad
\text{as} \quad \tau \to \infty.
\end{equation*}
Hence zero is linearly stable for $p>p_0$. $\qed$

%%\hfill$\Box$

%%\end{proof}

\begin{lemma}
For $1<p<p_0$, zero is exponentially linearly unstable for
\eqref{rescsem}, for small data $v_0\in H^{2k+1}_\rho\cap \tilde
L^2_\rho(\re)$.
\end{lemma}

\noindent {\em Proof.} For $l=0$,
 %%\begin{equation*}
$\lambda_0 = d_1\equiv\mbox{$\frac{1}{p-1}-\frac{1}{2k+1}$} >0$,
 %%\end{equation*}
so that
\begin{equation*}
\|v(\tau)\|_{2k+1,\rho}\sim \mathrm{e}^{\lambda_0\tau} \to +\infty.
\end{equation*}
Hence zero is linearly unstable for $p<p_0$. $\qed$

\smallskip

Thus, $p_0=2k+2$ plays almost a full role of the {\em critical
Fujita exponent} for \eqref{semilinear pde}, that mimics the
standard one $p_0=1+ \frac {2m} N$ ($=2m+1$ for $N=1$) for the
semilinear parabolic equations \eqref{Par1}; see \cite{SemiPDE}
for extra details.

%%\hfill$\Box$
%%\end{proof}

%%%%%%%%%%%%%%%%%%%%%%%%%%%%%%%%%%%%%%%%%%%%%%
%%%%%%%%%%%%%%%%%%%%%%%%%%%%%%%%%%%%%%%%%%%%%%
%%%%%%%%%%%%%%%%%%%%%%%%%%%%%%%%%%%%%%%%%%%%%%
%%%%%%%%%%%%%%%%%%%%%%%%%%%%%%%%%%%%%%%%%%%%%%
%%%%%%%%%%%%%%%%%%%%%%%%%%%%%%%%%%%%%%%%%%%%%%
\subsection{Centre subspace behaviour}

 Let us show that, at some
critical values of $p$, there exists some quite special
asymptotics of solutions.

\ssk

\noindent\underline{\sc First critical exponent $p_0$}. We look
again at the rescaled equation given by \eqref{rescsem}.
%\begin{align*}
%w_\tau &= (-1)^{k+1}D_y^{2k+1}w +\mbox{$\frac{1}{2k+1}yD_yw+\frac{1}{p-1}w-w^p$}\\
%&={\bf B}w -\mbox{$\big(\frac{1}{2k+1}-\frac{1}{p-1}\big)w-w^p$}.
%\end{align*}
%Where as before
%\begin{equation*}
%u(x,t) = t^{-\frac{1}{(2k+1)}}w(y,\tau), \quad y=xt^{-\frac{1}{(2k+1)}}, \quad \tau
%= \ln(1+t).
%\end{equation*}
Let us look at the case where $l=0$, with critical exponent
$p=p_0=2k+2$. We check the behaviour close to the centre subspace
of ${\bf B}$, i.e.,
%%% we set
\begin{equation*}
v(\tau)=c_0(\tau)\psi_0 + v_0^{\perp}(\tau),
\end{equation*}
where $v_0^{\perp}$ is asymptotically small in comparison with the
first term and is orthogonal to $\psi_0$, i.e., $\langle
v_0^{\perp}, \psi_k^{\ast}\rangle_\ast = 0$, for any $ k\geq 1$.
Then, since ${\bf B}\psi_0=0$, we multiply the equation  by
$\psi_0^\ast\equiv 1$ to get the following leading term:
\begin{equation*}
c_0^{\prime}=-|c_0|^{2k+1} c_0\langle |\psi_0|^{2k+1}\psi_0+...\,,
\psi_0^\ast\rangle_\ast.
\end{equation*}
We need to assume that (see a justification below)
\begin{equation*}
\gamma_0 = \langle |\psi_0|^{2k+1}\psi_0,\psi_0^\ast\rangle_\ast
>0.
 %%%\not = 0.
\end{equation*}
 Note that analytically, proving that $\gamma_0 \not =0$ is
difficult and even checking this numerically is also questionable.
So, assuming that $\g_0>0$,
\begin{equation*}
 \mbox{$
|c_0|^{-(2k+1)} \frac{c_0^\prime}{c_0} = -\gamma_0+...  \quad
\Longrightarrow \quad
 -\mbox{$\frac{1}{2k+1}$}\, |c_0|^{-(2k+1)}=-\gamma_0\tau+\hdots
 \quad \text{as} \quad \tau
\to \infty.
 $}
\end{equation*}
%%and therefore
%%\begin{equation*}
%%-\mbox{$\frac{1}{2k+1}$}c_0^{-(2k+1)}=-\gamma_0\tau+\hdots
%%\end{equation*}
Finally, this yields the following rate of decay:
\begin{equation*}
c_0(\tau)\approx\big[(2k+1)\gamma_0\tau\big]^{-\frac{1}{2k+1}}
\quad \text{for} \quad \tau\gg 1.
\end{equation*}

Thus,  the centre subspace behaviour of $u(x,t)$ is given, as
$t\to\infty$, by
\begin{equation*}
u(x,t)\approx
(1+t)^{-\frac{1}{2k+1}}\big[(2k+1)\gamma_0\ln(1+t)\big]^{-\frac{1}{2k+1}}
\psi_0\big(x/(1+t)^{\frac{1}{2k+1}}\big),
%%%\\[5mm]
\end{equation*}
i.e., contains a typical extra logarithmic factor. A full
justification of such a behaviour remains open (even existence of
an invariant centre manifold is obscure).
%%% existence

\ssk

%%%%%%%%%%%%%%%%%%%%%%%%%%%%%%%%%%%%%%%%
\noindent\underline{\sc Other critical exponents $p_l$}.
%% stable
%%subspace behaviour.} \noindent
 Let us now look at the general case
for $l$, with critical point $p=p_l=1+\mbox{$\frac{2k+1}{l+1}$}$,
where we check the behaviour close to the one-dimensional kernel
of ${\bf B}-\lambda_lI$, by setting
\begin{equation*}
v=c_l(\tau)\psi_l + v_l^{\perp}.
\end{equation*}
Here, $v_l^{\perp}$ is  small and orthogonal to $\psi_l$, as
before. Then for ${\bf B}\psi_l=\lambda_l\psi_l$, we take the
inner product with $\psi_l^\ast$ to get, for $\t \gg 1$,
\begin{equation*}
c_l^{\prime}\approx-|c_l|^{\frac{2k+1}{l+1}} c_l\langle
|\psi_l|^{\frac{2k+1}{l+1}}\psi_l, \psi_l^\ast\rangle_\ast.
\end{equation*}
For convenience, let us assume, as usual, that
\begin{equation*}
\gamma_l = \langle|\psi_l|^{\frac{2k+l}{l+1}}\psi_l ,
\psi_l^\ast\rangle_\ast >0.
%%%%%\quad (\not = 0).
\end{equation*}
Hence, assuming also that $c_l(\t)>0$, we have, for $\tau \gg 1$,
\begin{equation*}
c_l^{-\frac{2k+2+l}{l+1}}c_l^\prime \approx -\gamma_l
 \quad \Longrightarrow
 \quad
-\mbox{$\frac{l+1}{2k+1}$}c_l^{-\frac{2k+1}{l+1}}\approx-\gamma_l\tau
\quad \Longrightarrow
 \quad c_l(\tau)\approx\big(\mbox{$\frac{2k+1}{l+1}$}\gamma_l\tau\big)^{-\frac{l+1}{2k+1}}.
  %%%\mbox{so}
\end{equation*}
%%therefore
%%\begin{equation*}
%%-\mbox{$\frac{l+1}{2k+1}$}c_l^{-\frac{2k+1}{l+1}}\approx-\gamma_l\tau,
%%\end{equation*}
%%so that
%%\begin{equation*}
%%c_l(\tau)\approx\big(\mbox{$\frac{2k+1}{l+1}$}\gamma_l\tau\big)^{-\frac{l+1}{2k+1}}.
%%\quad \text{for} \quad \tau\gg 1.
%%\end{equation*}

Hence,
 the  stable subspace behaviour of $u(x,t)$ as $t\to\infty$,
for all critical exponents $p=p_l$, is given by
\begin{equation*}
u(x,t)\approx
\mbox{$(1+t)^{-\frac{1+l}{2k+1}}$}\big[\mbox{$\frac{2k+1}{l+1}\,
\gamma_l\ln(1+t)$}\big]^{-\frac{l+1}{2k+1}}\psi_l\big(x/(1+t)^{\frac{1}{2k+1}}\big).
\end{equation*}
So, there exists a countable set of stable subspace behaviours
governed by eigenfunctions, corresponding to the point spectrum
$\sigma({\bf B}) = \{-\mbox{$\frac{l}{2k+1}$}, \quad l\geq0\}$.
%%\\[5mm]

 %%\smallskip

%%\noindent{{\bf Remark:

%%%%%%%%%%%%%%%%%%%%%%%%%%%%%%%%%%%%%%%%%%%%%%%%%%%%%%%%%%%%%
 \subsection{Why $\gamma_l>0$: numerics}

   For the above analysis to hold, it is assumed that
 %%%%convenient to assume that (though this is not very crucial)
\begin{equation}
 \label{l01}
\gamma_l = \langle |\psi_l|^{\frac{2k+1+l}{l+1}}\psi_l,
\psi_l^\ast\rangle_* >0,
\end{equation}
which is not an easy inequality to prove. Note that, for all $l
\ge 0$, the eigenfunctions $\psi_l(y)$ and $\psi_l^*(y)$ are of
the same parity, so that the indefinite metric $\langle \cdot,
\cdot \rangle_*$ plays no role here and  can be  replaced by the
standard $L^2$-one (for odd $l$'s, with the minus sign).

Let us first look at the case $l=0$. We need to show that
\begin{equation*}
\gamma_0 = \langle |\psi_0|^{2k+1}\psi_0, \psi_0^\ast\rangle_*
>0.
\end{equation*}
However, in this case we have that $\psi_0^\ast\equiv 1$, so in
essence it is enough to prove that
\begin{equation*}
 \mbox{$
\gamma_0 = \int |\psi_0|^{2k+1}\psi_0
>0.
 $}
\end{equation*}
Attempting to prove this rigorously, is very difficult as well.
However, numerically it can be shown, for $k=1$ at least, that this
inequality is true.

Certainly, we know that the first eigenfunction $\psi_0 \equiv
F(y)$ can be found, for the lower-order case $k=1$.  Using the
{\tt MatLab} function {\tt trapz}, which uses a trapezoidal method
of integration,
%%%%we can solve. Using this method,
the integral can
be approximated to
\begin{equation}
\label{trapz1} \int |F|^3 F = 0.0300\hdots > 0.
\end{equation}
It is noted that, since this numerical solution $F$, is not the
fundamental kernel satisfying $\int F = 1$, we scale the
calculations such that this is true.  So \eqref{trapz1} holds for
the true rescaled fundamental kernel.

Solving the problem (\ref{l01}) for $k \ge 2$ is more difficult as
the shooting problem to find further eigenfunctions
 is not easy, but not impossible. Similarly, it is possible to
construct solutions for all $\psi_l$, where $l\geq 0$.
%% using
 %%conservation laws.
  Again the problem of shooting is difficult,
but not impossible to justify the inequality (\ref{l01}).

%%This work is ongoing and we hope to prove this more rigourously for
%%all $l$ and $k$.

%For the case $k=1$, this can be shown numerically, since in this
%case $\psi_0(y)\approx \mathrm{Ai}(-y)$.  We use the computer
%program {\em Maple} to integrate $\int \mathrm{Ai}^4(-y)$.  Whilst
%the program cannot evaluate the integral over $\mathbb{R}$, we can
%integrate over a finite range without loss of generality since we
%know that $F\to 0$ as $y\to\pm \infty$.  Also knowing that the
%behaviour of the solution is oscillatory for $x\gg 1$, we can just
%integrate over values relatively close to $0$.

%So using Maple we find that
%\begin{equation}
%\int_{-10}^{100} \mathrm{Ai}^4(-y) \, \mathrm{d}y = 0.2422609017 \,
%>0.
%\end{equation}

%%%%%%%%%%%%%%%%%%%%%%%%%%%%%%%%%%%%%%%%%%%%%%
%%%%%%%%%%%%%%%%%%%%%%%%%%%%%%%%%%%%%%%%%%%%%%
%%%%%%%%%%%%%%%%%%%%%%%%%%%%%%%%%%%%%%%%%%%%%%
%%%%%%%%%%%%%%%%%%%%%%%%%%%%%%%%%%%%%%%%%%%%%%
%%%%%%%%%%%%%%%%%%%%%%%%%%%%%%%%%%%%%%%%%%%%%%

\subsection{Bifurcation points and $p$-diagrams}
 \label{S7.4}

\noindent Once again, we look at our semilinear equation given in
\eqref{conv}, which we can write now as
\begin{equation*}
%\begin{split}
\mbox{${\bf B}f +\big(\frac{1}{p-1}-\frac{1}{2k+1}\big)f -
|f|^{p-1}f$} = 0 \Longleftrightarrow\mbox{${\bf B}_1f-|f|^{p-1}f$} =
0,
%\end{split}
\end{equation*}
where ${\bf B}_1= {\bf B}+ d_1I$ is defined as before. Critical
exponents $\{p_l\}$ occur when
\begin{equation*}
\mbox{$d_1\equiv\frac{1}{p_l-1} - \frac{1}{2k+1} $}=-\lambda_l,
 \quad \mbox{i.e., when} \quad
%%\end{equation*}
%%i.e., when
%%\begin{equation*}
\mbox{$\frac{1}{p_l-1} - \frac{1}{2k+1}$}=\mbox{$\frac{l}{2k+1}$}.
\end{equation*}
Therefore our critical exponents, $p_l$, are given by
\begin{equation*}
\mbox{$p_l = 1 + \frac{2k+1}{l+1}$}, \quad l=0,1,2,...\, .
\end{equation*}
We can see from this that $p_l \rightarrow 1^+$ as $l \rightarrow
+\infty$.

We look at  $p$ near these critical values, so that $p\approx
p_l$. We set $\e = p_l - p$, and then
\begin{equation}
\label{reducedsembif} ({\bf B} - \lambda_lI)f + \e a_0f =
|f|^{p-1}f + O(\e^2) \whereA
%%\end{equation}
%%where $a_0$ is some constant. We can find $a_0$, by substituting
%%in $\e$ and so we have
%%\begin{equation*}
a_0 =
\mbox{$\frac{1}{(p_l-1)^2}$}=\mbox{$\big(\frac{l+1}{2k+1}\big)^2$}.
\end{equation}

By classic bifurcation theory \cite{VainbergTr}, our solution $f$
can be given by
\begin{equation*}
f = C\psi _l + w^\bot,
\end{equation*}
where $w^\bot$ is orthogonal to $\psi _l$,
\begin{equation*}
\langle w^\bot,\psi _l^{\ast}\rangle_\ast = 0.
\end{equation*}
 %%So $f \approx C\psi _l$ and
Thus, taking the inner product of \eqref{reducedsembif} with
$\psi_l^\ast$, we have that
\begin{equation*}
a_0\e C\langle\psi _l, \psi _l^{\ast}\rangle_\ast =
|C|^{p-1}C\langle|\psi_l|^{p-1}\psi_l, \psi_l^{\ast}\rangle_\ast .
\end{equation*}
Since $\langle\psi _l, \psi _l^{\ast}\rangle_* = 1$, we find
\begin{equation}
 \label{C01}
\mbox{$|C|^{p-1} = \frac{a_0\e}{\langle|\psi_l|^{p-1}\psi_l,
\psi_l^{\ast}\rangle_\ast} =\frac{1}{\kappa
_l}\big(\frac{l+1}{2k+1}\big)^2\e$}, \quad
%%\end{equation*}
 \mbox{where}\quad
%%\begin{equation*}
\kappa_l=\langle|\psi_l|^{p-1}\psi_l, \psi_l^\ast\rangle_\ast.
\end{equation}
For $p=1$, we have that
%%\begin{equation*}
$\kappa _l=\langle\psi _l, \psi _l^{\ast}\rangle_\ast = 1$,
 %%\end{equation*}
and so by continuity with respect to $p$, we must have that
\begin{equation*}
\kappa _l>0 \quad \text{for all} \quad p\approx 1^+.
\end{equation*}
Therefore, (\ref{C01}) indicates  a countable number of {\em
subcritical pitchfork} bifurcations.

\begin{figure}[htbp]
\begin{center}
\includegraphics[scale=0.8]{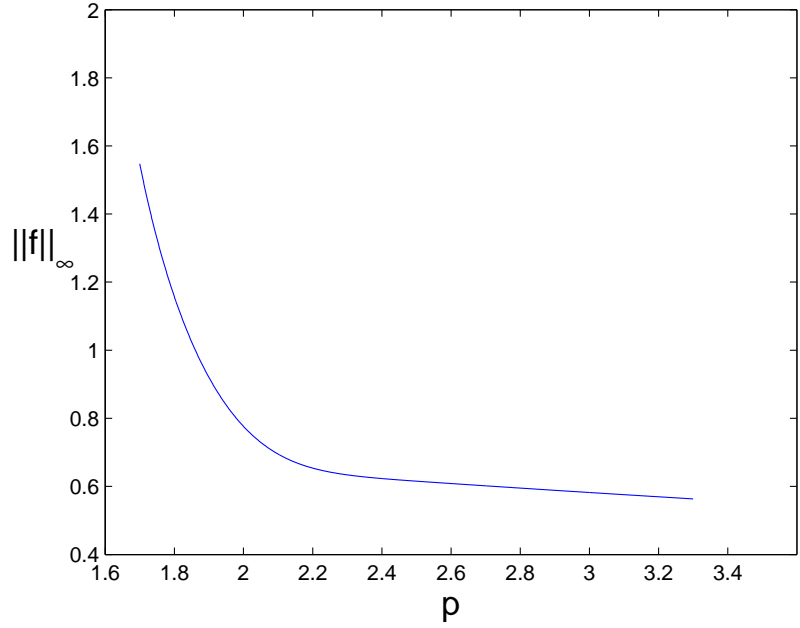}
\caption{\small Bifurcation branch for $l=0$ and $k=1$, for the
semilinear ODE \eqref{conv}.} \label{FFF1}
\end{center}
\end{figure}

Figure \ref{FFF1} shows a numerical calculation for the first
branch of the bifurcation diagram, where we take $l=0$ and $k=1$.
Hence, in this case, the critical point is $p_0=4$. During each
iteration of the numerical program, the calculation uses the
previous results to calculate the next step, thus improving the
accuracy. The step size used here is $0.001$, in the range of
$p=1.7$ to $p=3.3$.  Extending the range of values of $p$ closer
to the critical $p_0=4$ proves to be difficult.

\begin{figure}[htbp]
\begin{center}
\psfrag{0}{$0$}
\psfrag{p0}{$p_0$}
\psfrag{p1}{$p_1$}
\psfrag{p2}{$p_2$}
\psfrag{p3}{$p_3$}
\psfrag{pl}{$p_l$}
\psfrag{pl+1}{$p_{l+1}$}
\psfrag{1}{$1$}
\psfrag{p}{$p$}
\psfrag{finf}{$\|f\|_{\iy}$}
\psfrag{...}{$...$}
\includegraphics[scale=0.4]{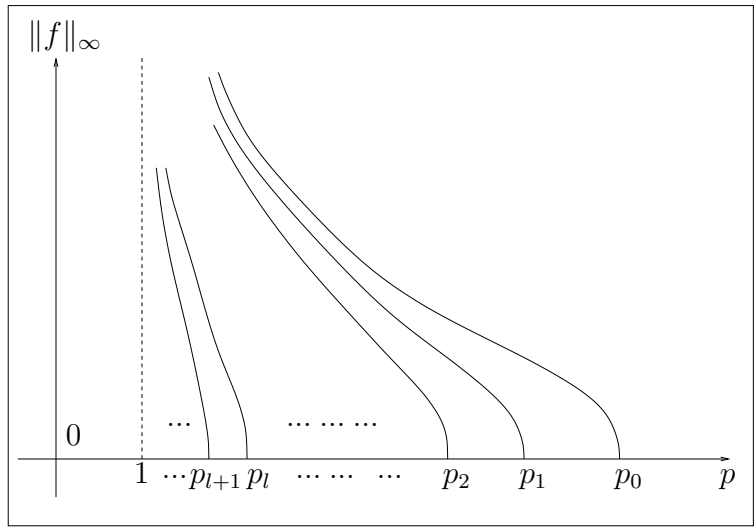}     %%%{bifdiag2anew.eps}
\caption{\small A countable family of expected $p$-bifurcation
branches for the ODE \eqref{conv}.}
  \label{FFF2}
\end{center}
\end{figure}

Figure \ref{FFF2} shows how the $p$-bifurcation diagram for the
semilinear equation is expected to look like, if the branch is
extended, given the numerics in Figure \ref{FFF1}.

%%%%%%%%%%%%%%%%%%%%%%%%%%%%%%%%%%%%%%%%%%%%%%%%%%%%%%%%%%%%%%%%%%%%%%%
%%%%%%%%%%%%%%%%%%%%%%%%%%%%%%%%%%%%%%%%%%%%%%%%%%%%%%%%%%%%%%%%%%%%%%%

%%%%%%%%%%%%%%%%%%%%%%%%%%%%%%%%%%%%%%%%%%%%%%%%%%%%%%%%%%%%%%%%%%%%%%%

%%%%%%%%%%%%%%%%%%%%%%%%%%%%%%%%%%%%%%%%%%%%%%%%%%%%%%%%%%%%%%%%%%%%%%%

However the numerical analysis done in Section \ref{NumerSem} has
failed to provide firm evidence that Figure \ref{FFF2} shows how
the $p$-branches behave, since we were not able to get reliable
numerics close to $p_0=4$. As before, this is related to the
extremal oscillatory behaviour of similarity profiles at the
right-hand side, which does not allow us using standard numerical
codes of continuation with respect to the parameter $p$, i.e.,
numerically construct the so-called $p$-branches of solutions. In
Figure \ref{FFF2}, we use the analytical evidence of such
bifurcations from zero at $p=p_l$, which also requires extra
difficult mathematical justification.

%%%%%%%%%%%%%%%%%%%%%%%%%%%%%%%%%%%%%%%%%%%%%%%%%%%%%%%%%%%%%%%%%%%
%\addcontentsline{toc}{chapter}{Bibliography}
\bibliographystyle{amsplain}
\bibliography{biblioFG}

\providecommand{\bysame}{\leavevmode\hbox to3em{\hrulefill}\thinspace}
\providecommand{\MR}{\relax\ifhmode\unskip\space\fi MR }
\providecommand{\MRhref}[2]{%
  \href{http://www.ams.org/mathscinet-getitem?mr=#1}{#2}
}
\providecommand{\href}[2]{#2}
\begin{thebibliography}{10}

\bibitem{linopindefmet}
T.Ya. Azizov and I.S. Iokhvidov, \emph{\rm {L}inear {O}perators in {S}paces
  with {I}ndefinite {M}etric}, Wiley-Interscience, New York/Toronto, 1989.

\bibitem{AsympPert}
C.M. Bender and S.A. Orszag, \emph{\rm {A}dvanced {M}athematical {M}ethods for
  {S}cientists and {E}ngineers; {A}symptotic {M}ethods and {P}erturbation
  {T}heory}, Springer, New York, 1999.

\bibitem{BS}
M.S. Birman and M.Z. Solomjak, \emph{\rm {S}pectral {T}heory of
  {S}elf-{A}djoint {O}perators in {H}ilbert {S}pace}, D. Reidel,
  Dordrecht/Tokyo, 1987.

\bibitem{Blas}
H.~Blasius, \emph{Grenzschichten in {F}l\"ussigkeiten mit kleiner {R}einbung},
  Z. Math. Phys. \textbf{56} (1908), 1--37.

\bibitem{Bouss}
J.~Boussinesq, \emph{Th\'{e}ories des ondes et des remous qui se propagent le
  long d'un canal rectangulaire horizontal, en communiquant au liquide contenu
  dans ce canal des vitesses sensiblement pareilles de la surface au fond}, J.
  Math. Pures Appl. \textbf{17} (1872), 55--108.

\bibitem{VssHeatAbs}
H.~Brezis, L.A. Peletier, and D.~Terman, \emph{A very singular solution of the
  heat equation with absorption}, Arch. Rational Mech. Anal. \textbf{95}
  (1986), 185--209.

\bibitem{Cai97}
H.~Cai, \emph{Dispersive smoothing effects for {K}d{V} type equations}, J.
  Differ. Equat. \textbf{136} (1997), 191--221.

\bibitem{Carl39}
T.~Carleman, \emph{Sur un probl\`eme d'unicit\'e pur les syst\`emes
  d'\'equations aux d\'eriv\'ees partielles \`a deux variables
  ind\'ependantes}, Ark. Mat. Astr. Fys. \textbf{26} (1939), 1--9.

\bibitem{SubcritGKdV}
R.~C\^{o}te, \emph{Construction of solutions to the subcritical g{K}d{V}
  equations with a given asymptotical behaviour}, J. Func. Anal. \textbf{241}
  (2006), 143--211.

\bibitem{Cr90}
W.~Craig and J.~Goodman, \emph{Linear dispersive equations of {A}iry type}, J.
  Differ. Equat. \textbf{87} (1990), 38--61.

\bibitem{Cr92}
W.~Craig, T.~Kappeler, and W.~Strauss, \emph{Gain of regularity for equations
  of {K}d{V} type}, Ann. Inst. H. Poincare \textbf{9} (1992), 147--186.

\bibitem{UniqPropDispEq}
L.L. Dawson, \emph{Uniqueness properties of higher order dispersive equations},
  J. Differ. Equat. \textbf{236} (2007), 199--236.

\bibitem{SemiPDE}
Yu.V. Egorov, V.A. Galaktionov, V.A. Kondratiev, and S.I. Pohozaev,
  \emph{Global solutions of higher- order semilinear parabolic equations in the
  supercritical range}, Adv. Differ. Equat. \textbf{9} (2004), 1009--1038.

\bibitem{Fam02}
A.V. Faminskii, \emph{On the mixed problem for quasilinear equations of the
  third order}, J.~Math. Sci. \textbf{110} (2002), 2476--2507.

\bibitem{FerGalII}
R.S. Fernandes and V.A. Galaktionov, \emph{{E}igenfunctions and very singular
  similarity solutions of odd-order nonlinear dispersion {P}{D}{E}s},
  submitted.

\bibitem{Dos05}
D.~Dos~Santos Ferreira, \emph{Sharp ${L}^p$ {C}arleman estimates and unique
  continuation}, Duke Math. J. \textbf{129} (2005), 503--550.

\bibitem{GalGeom}
V.A. Galaktionov, \emph{\rm {G}eometric {S}turmian {T}heory of {N}onlinear
  {P}arabolic {E}quations and {A}pplications}, Chapman and Hall/CRC, Florida,
  2004.

\bibitem{2mSturm}
\bysame, \emph{Sturmian nodal set analysis for higher-order parabolic equations
  and applications}, Adv. Differ. Equat. \textbf{12} (2007), 669--720.

\bibitem{GalANS}
\bysame, \emph{Non-radial very singular solutions of absorption-diffusion
  equations with non-homogeneous potentials}, Adv. Nonl. Stud. \textbf{8}
  (2008), 429--454.

\bibitem{GPndeII}
\bysame, \emph{Nonlinear dispersion equations: smooth deformations, compactons,
  and extensions to higher orders}, Comput. Math. Math. Phys. \textbf{48}
  (2008), 1823--1856 (arXiv:0902.0275).

\bibitem{GalKamLSE}
V.A. Galaktionov and I.V. Kamotski, \emph{{R}efined scattering and {H}ermitiam
  spectral theory for linear {S}chr\"odinger eqautions with applications}, in
  preparation.

\bibitem{AsymEigNonlinPara}
V.A. Galaktionov, S.P. Kurdyumov, and A.A. Samarski\u{\i}, \emph{On asymptotic
  ``eigenfunctions" of the {C}auchy problem for a nonlinear parabolic
  equation}, Math. USSR Sbornik \textbf{54} (1986), 421--455.

\bibitem{MajOrdOp}
V.A. Galaktionov and S.I. Pohozaev, \emph{Existence and blow-up for
  higher-order semilinear parabolic equations: majorizing order preserving
  operators}, Indiana Univ. Math. J. \textbf{51} (2002), 243--280.

\bibitem{GPnde}
\bysame, \emph{Third-order nonlinear dispersive equations: shocks, rarefaction,
  and blow-up waves}, Comput. Math. Math. Phys. \textbf{48} (2008), 1784--1810
  (arXiv:0902.0253).

\bibitem{GSVR}
V.A. Galaktionov and S.R. Svirshchevskii, \emph{\rm {E}xact {S}olutions and
  {I}nvariant {S}ubspaces of {N}onlinear {P}artial {D}ifferential {E}quations
  in {M}echanics and {P}hysics}, Chapman and Hall/CRC, Florida, 2007.

\bibitem{StabTechPdeVSS}
V.A. Galaktionov and J.L. V\'{a}zquez, \emph{\rm {A} {S}tability {T}echnique
  for {E}volution {P}artial {D}ifferential {E}quations}, Birkh\"{a}user,
  Berlin/Boston, 2004.

\bibitem{VSSSParaPDEs}
V.A. Galaktionov and J.F. Williams, \emph{On very singular similarity solutions
  of a higher-order semilinear parabolic equation}, Nonlinearity \textbf{17}
  (2004), 1075--1099.

\bibitem{Hos99}
T.~Hoshiro, \emph{Mouree's method and smoothing properties of dispersive
  equations}, Comm. Math. Phys. \textbf{202} (1999), 255--265.

\bibitem{UnitOpIndef}
I.S. Iokhvidov, \emph{Unitary operators in a space with an indefinite metric},
  Zap., N.I.I. Mat. i Mekh. Khar'kov Gos. Univ. Mat. Obsch. (1949), 79--86.

\bibitem{Ion04}
A.D. Ionescu and C.E. Kenig, \emph{${L}^p$ {C}arleman inequalities and
  uniqueness of solutions of nonlinear {S}chr\"odinger equations}, Acta. Math.
  \textbf{193} (2004), 193--239.

\bibitem{Ion06}
\bysame, \emph{Uniqueness properties of solutions of {S}chr\"odinger
  equations}, J. Funct. Anal. \textbf{232} (2006), 90--236.

\bibitem{SSHeatAbs}
S.~Kamin and L.A. Peletier, \emph{Singular solutions of the heat equation with
  absorption}, Proc. Math. Soc. \textbf{95} (1985), 205--210.

\bibitem{ExisUniqVSSPorMed}
S.~Kamin and L.~Veron, \emph{Existence and uniqueness of the very singular
  solution of the porous media equation with absorption}, J. Anal. Math.
  \textbf{51} (1988), 245--258.

\bibitem{QuasiEvoPDE}
T.~Kato, \emph{Quasi-linear equations of evolution, with applications to
  partial differential equations}, {S}pectral {T}heory and {D}ifferential
  {E}quations ({P}roc. {S}ympos., {D}undee, 1974), pp.~27--50. Lecture Notes in
  Math., Vol. 448, Springer, Berlin, 1975.

\bibitem{OnKdV}
\bysame, \emph{On the {K}ortweg-de {V}ries equation}, Manuscripta Math.
  \textbf{28} (1979), 89--99.

\bibitem{BilinEst}
C.E. Kenig, G.~Ponce, and L.~Vega, \emph{A bilinear estimate with applications
  to the {K}d{V} equation}, Proc. Amer. Math. Soc. \textbf{9} (1996), 573--603.

\bibitem{FuncAnaly}
A.N. Kolmogorov and S.V. Fomin, \emph{\rm {F}unctional {A}nalysis: {V}olume 1},
  4th ed., Graylock, Rochester, 1957.

\bibitem{KolmF}
\bysame, \emph{\rm {E}lements of {T}heory of {F}unctions and {F}unctional
  {A}nalysis}, 4th ed., Nauka, Glav. Red. Fiz--Mat. Lit., Moscow, 1976.

\bibitem{HelCurv}
M.G. Krein, \emph{Helical curves in an infinite-dimensional {L}obachevskiy
  space and {L}orents transformation}, Izvestiya Akad. Nauk USSR, Ser. Matem.
  \textbf{3} (1948), 3.

\bibitem{Lar06}
N.A. Larkin, \emph{Modified {K}d{V} equation with a source term in a bounded
  domain}, Math. Meth. Appl. Sci. \textbf{29} (2006), 751--765.

\bibitem{Led08}
Yu.S. Ledyaev, \emph{Private communication}, 2008.

\bibitem{Lev01}
J.L. Levandosky, \emph{Smoothing properties of nonlinear dispersive equations
  in two spatial dimensions}, J. Differ. Equat. \textbf{175} (2001), 275--372.

\bibitem{LustSob}
L.~Ljusternik and V.~Sobolev, \emph{\rm {E}lements of {F}unctional {A}nalysis},
  Ungar Publ. Comp., New York, 1961.

\bibitem{Maz}
V.~Maz'ya, \emph{\rm {S}obolev {S}paces}, Springer-Verlag, Berlin, 1985.

\bibitem{Miz06}
R.~Mizuhara, \emph{The initial value problem for third and fourth order
  dispersive equations in one space dimension}, Funk. Ekvacioj \textbf{49}
  (2006), 1--38.

\bibitem{NaimarkI}
M.A. Naimark, \emph{\rm {L}inear {D}ifferential {O}perators, {P}art {I}}, Ungar
  Publ. Comp., New York, 1967.

\bibitem{HermOpIndefMet}
L.S. Pontryagin, \emph{Hermitian operator in spaces with indefinite metric},
  Izvestiya Akad. Nauk USSR, Ser. Matem. \textbf{8} (1944), 243--80.

\bibitem{Prandtl}
L.~Prandtl, \emph{\rm \"{U}ber {F}l\"ussigkeitsbewegung bei sehr kleiner
  {R}eibung, {I}n: {P}roc. 3rd {I}nt. {C}ongr., {H}eidelberg, 1904}, Teubner,
  Leipzig, 1905.

\bibitem{quasilin}
A.A. Samarskii, V.A. Galaktionov, S.P. Kurdyumov, and A.P. Mikhailov, \emph{\rm
  {B}low-up in {Q}uasilinear {P}arabolic {E}quations}, Walter de Gruyter,
  Berlin, 1995.

\bibitem{PDEPhysSommer}
A.~Sommerfeld, \emph{\rm {P}artial {D}ifferential {E}quations in {P}hysics},
  Academic, New York, 1949.

\bibitem{DispEnerWav}
W.A. Strauss, \emph{Dispersion of low-energy waves for two conservative
  equations}, Arch. Ration. Mech. Anal. \textbf{55} (1974), 86--92.

\bibitem{St}
C.~Sturm, \emph{M\'{e}moire sur une classe d'\'{e}quations \`{a}
  diff\'{e}rences partielles}, J. Math. Pures Appl. \textbf{1} (1836),
  373--444.

\bibitem{Tak06}
H.~Takuwa, \emph{Microlocal analytic smoothing effects for operators of real
  principal type}, Osaka J.~Math. \textbf{43} (2006), 13--62.

\bibitem{Tao00}
T.~Tao, \emph{Multilinear weighted convolution of $l^2$ functions, and
  applications to nonlinear dispersive equations}, Amer. J. Math. \textbf{6}
  (2000), 839--908.

\bibitem{Tao08}
X.~Tao and S.~Zhang, \emph{Weighted doubling properties and unique continuation
  theorems for the degenerate {S}chr\"odinger equations with singular
  potentials}, J. Math. Anal. Appl. \textbf{339} (2008), 70--84.

\bibitem{VainbergTr}
M.A. Vainberg and V.A. Trenogin, \emph{\rm {T}heory of {B}ranching of
  {S}olutions of {N}on-{L}inear {E}quations}, Noordhoff Int. Publ., Leiden,
  1974.

\bibitem{SommerVibrFluid}
J.T. Xing, \emph{An investigation into natural vibrations of fluid-structure
  interaction systems subject to {S}ommerfeld radiation condition}, Acta Mech.
  Sin. \textbf{24} (2008), 69--82.

\end{thebibliography}

%\nocite{SolitAndNonlinEq, quasilin, DSandErgodic, LinopAndLinsys,
%IntroToNPDEs, FuncAnaly, AppNonSem, AppNonAnaly, SemiPDE,
%VSSSParaPDEs, NonFuncAnaly, GeomMethNonAn}

%%%%%%%%%%%\begin{appendix}
%%%%%%%%%%%%%%%%%%%%%%%%%%%%%%%%%%%%%%%%%%%%%%%%%%%%%%%%%%%%%
%%%%%%%%%%%%%%%%%%%%%%%%%%%%%%%%%%%%%%%%%%%%%%%%%%%%%
\begin{appendix}
\section*{Appendix A. ``Radiation conditions" in spectral $\{\BB,\BB^*\}$-theory}
 %%\label{Sect4.N}
 \label{SectRadcond}
 \setcounter{section}{1}
\setcounter{equation}{0}

\begin{small}

%%%%%%%%%%%%%%%%%%%%%%%%%%%%%%%%%%%%%%%%%%%%%%%%%
%%\section{``Radiation conditions" in spectral theory}
%%\label{SectRadcond}

 Here, we  clarify the ``radiation-type conditions"  posed at
infinity in the domains of operators $\BB$ and $\BB^*$, which
allow the operator pair $\{{\bf B},\, {\bf B}^*\}$ to have purely
discrete spectra $\{- \frac l{2k+1}, \, l \ge 0\}$ already
detected by eigenfunction expansion of the corresponding
semigroups. This serves as a description of new properties of the
spaces of closure $\tilde L^2_\rho(\re)$ and $\tilde
L^2_{\rho^*}(\re)$.

%%%%%%%%%%%%%%%%%%%%%%%%%%%%%%%%%%%%%%%%%%%%%%%%%%%%%%%%
\subsection{Domain $\tilde {\mathcal D}(\BB)$ of  the linear operator ${\bf B}$}

To this end,  in order for our linear ODE to be ``well-posed"
(i.e., with a proper number of boundary conditions at infinity),
both in a mathematical sense, as well as a physical sense, we look
for conditions, which must be satisfied for the  eigenvalue
equation:
\begin{equation*}
{\bf B}\psi_l(y) = \lambda_l\psi_l(y) \quad \mbox{in} \quad \re.
\end{equation*}
%%is satisfied.
 This can be rewritten as
\begin{equation}
\label{eigvalradcond} \mbox{$(-1)^{k+1}\psi_l^{(2k+1)}(y) +
\frac{1}{2k+1}\,\psi_l(y) + \frac{1}{2k+1}\,y\psi_l^\prime(y)  =
\lambda_l\psi_l(y)$}.
\end{equation}
Since the order of this ODE is $2k+1$, it is natural that there
must also be $2k+1$ boundary conditions placed, as classic theory
of ordinary differential operators suggests; see Naimark's
monograph \cite{NaimarkI}.

\smallskip

\noi\underline{\sc Eliminating exponentially growing bundles}.
First consider the problem, as $y\to+\infty$. Attempting to
balance leading order terms in \eqref{eigvalradcond}, leads to
\begin{equation}
\label{eigvalode} \mbox{$(-1)^{k+1}\psi_l^{(2k+1)}(y) +
\frac{1}{2k+1}\,y\psi_l^\prime (y) \sim 0$}.
\end{equation}
As $y\to+\infty$, we have that
 \begin{equation}
  \label{exp56}
\psi_l(y)\sim\mathrm{e}^{b y^{\frac{2k+1}{2k}}}, \quad b \in
{\mathbb C}, \quad b \not = 0,
 \end{equation}
 hence
substituting this into the above equation, yields
\begin{equation*}
 %%\begin{split}
\mbox{$(-1)^{k+1}\big(\frac{2k+1}{2k}\,b\big)^{2k+1} +
\frac{1}{2k+1}\,\big(\frac{2k+1}{2k}\,b\big)$} \sim 0 \quad
\Longrightarrow \quad b^{2k} \sim
\mbox{$(-1)^{k}\big(\frac{2k}{2k+1}\big)^{2k}\frac{1}{2k+1}$}.
%%\end{split}
\end{equation*}
It can be seen that we end up with two cases, dependent on the
value of the parameter $k$, which will determine the sign of
$b^{2k}$ and therefore its roots. Hence, for now, we ignore the
term $\big(\frac{2k}{2k+1}\big)\frac{1}{2k+1}$ and just look at
the value of $(-1)^k$, assuming that
 $
 \mbox{$
 b=\big(\frac{2k}{2k+1}\big)\frac{1}{2k+1}\, \hat b.
 $}
 $

{\em When $k$ is even}: for even values of $k$, it is noted that
 %%\begin{equation*}
  $\hat b^{2k} = 1$.
 %%\end{equation*}
Hence there must $2k$ roots for $ \hat b$, which are given by
 %%\begin{equation*}
 $\mbox{$\hat b_m = \mathrm{e}^{\frac{2\pi m \mathrm{i}}{2k}}$}$
 %%\end{equation*}
for $m= 0,1,\hdots,2k-1$.

{\em When $k$ is odd}:  for odd values of $k$, there is now a
negative sign, such that
%%\begin{equation*}
$\hat b^{2k} = -1$.
 %%\end{equation*}
Similarly, as before, we derive $2k$ roots, where now
%%\begin{equation*}
$\mbox{$\hat b_m = \mathrm{e}^{\frac{(\pi + 2\pi
m)\mathrm{i}}{2k}}$}$
%% \quad \text{for} \quad
for $ m = 0,1,\hdots,2k-1$.
%%\end{equation*}

For the problem to be well posed on the space $L^2_\rho$, it is
important to look for roots such that there is exponential decay,
rather than growth.  In other words, it must satisfy the condition
$\rho = \mathrm{e}^{-ay^\alpha}$, as $y\to+\infty$.  Hence we need
to eliminate any roots such that $\mathrm{Re}\, \hat b_m >0$. It
is first noted that equality, $\mathrm{Re}\, \hat b_m =0$, occurs
when
\begin{equation*}
 \mbox{$\frac{\pi}{2}$} =
\begin{cases}
\mbox{$\frac{\pi m }{k}$} \quad  \text{for even $k$},\\
\mbox{$\frac{(\pi + 2\pi m)}{2k}$} \quad \text{for odd $k$},
\end{cases}
\end{equation*}
with the same applying for $\frac{3\pi}{2}$.  For $\mathrm{Re}\,\,
\hat b_m =0$, we must have that
\begin{equation*}
m=
\begin{cases}
\mbox{$\frac{k}{2}$} \quad  \text{for even $k$},\\
\mbox{$\frac{k-1}{2}$} \quad \text{for odd $k$}
\end{cases}
 \quad
 \mbox{and} \quad
m=
\begin{cases}
\mbox{$\frac{3k}{2}$} \quad  \text{for even $k$},\\
\mbox{$\frac{3k-1}{2}$} \quad \text{for odd $k$}.
\end{cases}
\end{equation*}
%%and
%%\begin{equation*}
%%m=
%%\begin{cases}
%%\mbox{$\frac{3k}{2}$} \quad  \text{for even $k$},\\
%%\mbox{$\frac{3k-1}{2}$} \quad \text{for odd $k$}.
%%\end{cases}
%%\end{equation*}
Hence, in order to eliminate roots which give rise to exponential
growth, the conditions placed must be such that we do not include
roots such that
\begin{equation*}
\label{yinftycondns}
\begin{cases}
\mbox{$m < \frac{k}{2} \quad \text{and} \quad m > \frac{3k}{2}$}
\quad  \text{for even $k$},\\ \mbox{$m < \frac{k-1}{2} \quad
\text{and} \quad m > \frac{3k-1}{2}$} \quad \text{for odd $k$}.
\end{cases}
\end{equation*}
Therefore, by taking the weight
\begin{equation}
 \label{exp78}
\mbox{$\rho(y) = \mathrm{e}^{-ay^{\frac{2k+1}{2k}}}$}, \quad
\mbox{with any sufficiently small $a>0$}
\end{equation}
we eliminate all exponentially growing oscillatory bundles. A
sharp bound on admissible $a>0$ will be derived.
%%\eqref{yinftycondns}.
Thus,  we have $k-1$ conditions placed
%% here, satisfying
%%\eqref{eigvalode}
  as $y\to+\infty$.

Similarly, we can do the same analysis for $y\to-\infty$. As
$y\to-\infty$ in \eqref{eigvalode}, we have that
  \begin{equation}
   \label{exp77}
   \psi(y)\sim\mathrm{e}^{b
(-y)^{\frac{2k+1}{2k}}}
 \quad \Longrightarrow \quad
\mbox{$b^{2k} \sim
(-1)^{k+1}\big(\frac{2k}{2k+1}\big)\frac{1}{2k+1}$}.
 \end{equation}
  %%and substituting into equation
%%\eqref{eigvalode}, yields
%%\begin{equation*}
%%\begin{split}
%%\mbox{$b^{2k} \sim
%%(-1)^{k+1}\big(\frac{2k}{2k+1}\big)\frac{1}{2k+1}$}.
%%\end{split}
 %%\end{equation*}
As expected,  this only differs from the $y\to+\infty$ case by the
opposite sign.  This leads to
\begin{equation*}
\hat{b}_m =
\begin{cases}
\mbox{$\mathrm{e}^{\frac{(\pi + 2\pi m)\mathrm{i}}{2k}}$}, \quad
\text{for even $k$},\\ \mbox{$\mathrm{e}^{\frac{\pi m
\mathrm{i}}{k}}$}, \quad  \text{for odd $k$}, \,\,\, \mbox{for $m
= 0,1,\hdots,2k-1$}.
\end{cases}
\end{equation*}
%%%for $m = 0,1,\hdots,2k-1$.

As before, we do not want roots such that $\mathrm{Re}\,\, b_m>0$,
where in this case
\begin{equation*}
\begin{cases}
\mbox{$m <  \frac{k-1}{2} \quad \text{and} \quad m >
 \frac{3k-1}{2} $} \quad  \text{for even $k$},\\
\mbox{$m <  \frac{k}{2}  \quad \text{and} \quad m > \frac{3k}{2}
$} \quad \text{for odd $k$}.
\end{cases}
\end{equation*}
Hence, for the weight \ef{H03},
%%\begin{equation}
%% \label{w910}
%%\mbox{$ \rho(y) = \mathrm{e}^{a|y|^{\frac{2k+1}{2k}}}$},
%% \quad
%%\mbox{with sufficiently small $a>0$},
%%\end{equation}
we eliminate these conditions. This leads to a further $k$
conditions, which are placed at $y\to-\infty$. All these
correspond to eliminating exponentially growing asymptotic bundles
that is not enough as we explain below.

\smallskip

\noi\underline{\sc Domain $\tilde {\mathcal D}(\BB)$: radiation
condition}.  As stated
 before, we look for $2k+1$ conditions to be posed onto the
problem.  Hence there are two more conditions needed. These
conditions at $y\to\pm\infty$ are known as {\em radiation
conditions}. In classic problems of quantum mechanics, acoustics,
and  physics, the general idea behind radiation conditions is that
energy sources must exactly be that and not sinks of energy. Hence
all energy must be radiated from a point and scatter to infinity.
We refer to the book by Sommerfeld \cite{PDEPhysSommer},  who (in
1912) first proposed radiation conditions for the {\em Helmholtz
equation},
 and to
 %%We also refer to the paper by Xing
 \cite{SommerVibrFluid} for recent applications and references.
 %%%%, which applies the radiation condition.

In our problem, the radiation conditions are rather tricky and
have almost nothing to do with the classic
   ones. We recall
that we identify those just for convenience (to verify the domain
of ${\bf B}$ and ${\bf B}^*$), since the eigenfunction expansions
of the semigroups, as the main tool of our asymptotic analysis,
automatically include the necessary two conditions at infinity.
 %%as we show below.

\ssk

\noi{\bf Remark.} As is shown in \cite{GalKamLSE}, for the linear
rescaled Schr\"odinger operator in $L^2_\rho(\ren)$,
$\rho={\mathrm e}^{-a|y|^{\a}}$, $\a= \frac{2m}{2m-1}$ (a standard
setting),
 $$
  \mbox{$
 \BB=-\ii (-\D)^m + \frac {\rm 1}{2m}\, y \cdot \n + \frac {N}{2m}\, I,
  \quad \BB^*=-\ii (-\D)^m - \frac {\rm 1}{2m}\, y \cdot \n,
  $}
 $$
the absence of a ``radiation condition", accepting ``oscillatory
functions" at infinity only, will consequence in the fact that, in
addition to the standard required spectrum $\s(\BB)=\{-
\frac{|\b|}{2m}, \, |\b| \ge 0\}$, there appears another symmetric
real one $\{+\frac{|\b|}{2m}, \, |\b| \ge 0\}$ with generalized
Hermite polynomial eigenfunctions constructed as in
(\ref{poleigfunc}). Fortunately, this is not the case for our
operators (\ref{linODE}), since  the polynomials
%% are not in
$ \not \in L^2_\rho$ for the weight in (\ref{H03}).

\ssk

%%%The necessity of radiation conditions at infinity is clearly seen
%%%not only form the above ``algebraic discrepancy", but also  from
%%%the fact that the operator $\BB$ in \eqref{linODE} in

The origin of the radiation condition for {\bf B} is as follows:
Let us now balance all lower-order terms in the eigenvalue problem
\eqref{eigvalradcond}, so that
\begin{equation*}
\mbox{$\frac{1}{2k+1}\,\psi_l(y) +
\frac{1}{2k+1}\,y\psi_l^\prime(y) $}\sim
\mbox{$\lambda_l\psi_l(y)$}.
\end{equation*}
By integration we can easily see that
%\begin{equation*}
%\begin{split}
%\int\mbox{$\frac{1}{\psi_l(y)}$}\,\mathrm{d}\psi & \sim\int
%\mbox{$\frac{1}{y}\,\big[(2k+1)\lambda_l - 1\big]\,\mathrm{d}y$} \\[1mm]
%\Longrightarrow \quad \psi_l(y) &\sim A y^{(2k+1)\lambda_l - 1},
%\end{split}
%\end{equation*}
\begin{equation}
\label{PsiRadRatPoly}
 \psi_l(y) \sim A y^{(2k+1)\lambda_l - 1},
\end{equation}
for some constant $A$. Of course, this corresponds to the obvious
root $b=0$ in the exponential expansions \eqref{exp56} and
\eqref{exp77}.

We note that \eqref{PsiRadRatPoly}, is a ``rational" function,
unlike the exponentially oscillatory bundles in \eqref{exp56} and
\eqref{exp77}. For $y\to-\infty$, we know that, for any $\l \in
{\mathbb C}$, rational solutions such as (\ref{PsiRadRatPoly}) do
not belong to the space $L^2_\rho$ with the exponentially growing
weight \eqref{H03}.

Thus, overall, we conclude as follows:
 \begin{equation}
  \label{pr33}
  \mbox{at $y=-\infty$, the proper weight (\ref{H03}) generates
  $k+1$ conditions.}
   \end{equation}
So, this is a usual and a standard situation, so that the singular
point $y=-\infty$ does not require any radiation-type condition.
This is not the case for the ``oscillatory" end-point $y=+\iy$.

For $y\to+\infty$, consider all complex ``eigenvalues",
$\lambda_l\in\mathbb{C}$, such that $\lambda_l = P+\mathrm{i}Q$,
for some $P,Q\in\mathbb{R}$.  From \eqref{PsiRadRatPoly},
$\psi_l(y)$ may now be given by
\begin{equation}
\label{complbundle}
 \psi_l(y) \sim y^{\hat{P} +\mathrm{i}\hat{Q}},
\end{equation}
for $\hat{P},\hat{Q}\in\mathbb{R}$.  Hence we see that
\begin{equation}
 %%\begin{split}
\label{RadWeakOscill} \psi_l(y) \sim
y^{\hat{a}}\mathrm{e}^{\mathrm{i}\hat{b}\ln y }\\[1mm ] \sim
y^{\hat{a}}\big(\cos(\hat b\ln y)+\mathrm{i}\sin(\hat b\ln
y)\big),
 %%%\end{split}
\end{equation}
as $y\to+\infty$.  However, we know  from the asymptotic analysis
that the behaviour of proper eigenfunctions is different and given
by a distinct type of higher  oscillatory functions:
\begin{equation*}
\psi_l(y) \sim  y^{-\frac{2k-1}{4k}}\,\cos(y^{\frac{2k+1}{2k}}),
\end{equation*}
which obviously gives a stronger oscillatory behaviour than a pure
$\cos(\hat{b} \ln y)$ in \eqref{RadWeakOscill}. However,
\eqref{RadWeakOscill} admits  weaker oscillatory behaviour and so
we must place a condition to eliminate this behaviour. We recall
the all the proper eigenfunctions being given by the generating
formula
 $
 %%\mbox{$
  \psi_\beta(y) = \frac 1{\sqrt{\beta!}} D^\beta F(y),
  %%$}
  $
{\em do not contain the bundle} \eqref{PsiRadRatPoly}, since the
fundamental rescaled kernel $F(y)$ does not by the known
divergence of the operator {\bf B} (the equation for $F$ has been
integrated once with the zero constant of integration that
eliminated any trace of \eqref{PsiRadRatPoly}); see %%computations
below.

Thus, the generalized radiation condition, that is necessary for
the proper discreteness of the spectrum of ${\bf B}$  in
$L^2_\rho$, can formally be formulated as follows:
%%\\[5mm]
%%\quad
  \be
  \label{RD32}
\fbox{$
  \mbox{for the eigenvalue equation (\ref{eigvalradcond}),
the bundle (\ref{complbundle}) at $y=+\infty$ is absent.}
 $}
 \ee
 %%%\\[5mm]
 Then, as usual, the domain
 By the construction of the space of closure $\tilde L^2_\rho$ in
 Section \ref{S4.3}, the condition \ef{RD32} is assumed to be
 included, though we do not prove this here; cf \cite{GalKamLSE}
 for linear rescaled Schr\"odinger operators.
%%$\tilde {\mathcal D}(\BB)$ is obtained by closure of the
%%eigenfunction subset $\{\psi_\b\}$ (not an easy procedure in
%%general).

Actually, it is easy to see that all our eigenfunctions
$\{\psi_l(y)\}$ satisfy the above radiation condition. Indeed,
$F(y)\equiv\psi_0(y)$ does satisfy this, by integrating once,
where we have that
\begin{equation*}
 %%\begin{split}
\mbox{$(-1)^{k+1}F^{(2k+1)} + \frac{1}{2k+1}\,(Fy)^\prime$} = 0
\quad  \Longrightarrow \quad \mbox{$(-1)^{k+1}F^{(2k)} +
\frac{1}{2k+1}\,Fy$} = C,
 %%\end{split}
\end{equation*}
for some constant of integration $C$. The last term precisely
shows that such a rational behaviour is absent, since we have that
\begin{equation*}
\mbox{$F(y) \sim \frac{(2k+1)C}{y}$},
\end{equation*}
which implies that $C\equiv0$.  Then each eigenfunction
 %%\begin{equation*}
$\psi_l(y) = \frac{1}{\sqrt{l !}}\,F^{(l)}(y)$
%% \quad
%%\text{for all} \quad
 for all $l\geq 0$,
 %%%%\end{equation*}
also satisfies this condition.

%%%%%%%%%%%%%%%%%%%%%%%%%%%%%%%%%%%%%%%%%%%%%%%%%%%%%%%%%%%
\subsection{Radiation condition in the domain $\tilde {\mathcal D}(\BB^*)$
 of  adjoint operator ${\bf B}^\ast$}

We apply the same analysis to the adjoint operator ${\bf B}^\ast$,
where the eigenvalue problem is now
\begin{equation}
\label{radadjeig} \mbox{$(-1)^{k+1}{\psi_l^{\ast}}^{(2k+1)}(y)
-\frac{1}{2k+1}\,y{\psi^{\ast}}^\prime(y) =
\lambda_l\psi^\ast(y)$}.
\end{equation}
Similarly, as before, we look at exponential bundles with $b \not
= 0$:
\begin{equation*}
\psi^\ast(y) \sim \mathrm{e}^{b y^{\frac{2k+1}{2k}}}
 \,\,\mbox{as} \,\, y \to +\iy, \quad \mbox{and}
 \quad \psi^\ast(y) \sim \mathrm{e}^{b |y|^{\frac{2k+1}{2k}}}
 \,\,\mbox{as} \,\, y \to -\iy.
\end{equation*}
%%as $y\to+\infty$ and
%%\begin{equation*}
%%\psi^\ast(y) \sim \mathrm{e}^{b y^{\frac{2k+1}{2k}}},
%%\end{equation*}
%%%as $y\to-\infty$.

Balancing terms
\begin{equation}
\mbox{$(-1)^{k+1}{\psi_l^{\ast}}^{(2k+1)}(y) =
\frac{1}{2k+1}\,y{\psi^{\ast}}^\prime(y)+... $},
\end{equation}
it can easily be seen that the case for the adjoint operator ${\bf
B}^\ast$, only differs from the linear operator ${\bf B}$, with
respect to a change of sign.  Hence it is seen that, as $y\to +
\infty$ and $y \to - \iy$ respectively, for $m = 0,1,\hdots,2k-1$,
\begin{equation*}
\hat{b}_m =
\begin{cases}
\mbox{$\mathrm{e}^{\frac{(\pi + 2\pi m)\mathrm{i}}{2k}}$}, \quad
\text{for even $k$},\\ \mbox{$\mathrm{e}^{\frac{\pi m
\mathrm{i}}{k}}$}, \quad  \text{for odd $k$},
\end{cases}
 \quad \mbox{and}
 \quad
\hat{b}_m =
\begin{cases}
\mbox{$\mathrm{e}^{\frac{\pi m \mathrm{i}}{k}}$}, \quad  \text{for
even $k$},\\ \mbox{$\mathrm{e}^{\frac{(\pi + 2\pi
m)\mathrm{i}}{2k}}$}, \quad \text{for odd $k$}.
\end{cases}
\end{equation*}
%%for $m = 0,1,\hdots,2k-1$.  When $y\to-\infty$, we have that
%%\begin{equation*}
%%a_m =
%%\begin{cases}
%%\mbox{$\mathrm{e}^{\frac{\pi m \mathrm{i}}{k}}$}, \quad  \text{for even $k$},\\
%%\mbox{$\mathrm{e}^{\frac{(\pi + 2\pi m)\mathrm{i}}{2k}}$}, \quad
%%\text{for odd $k$},
%%\end{cases}
%%\end{equation*}
%%for $m = 0,1,\hdots,2k-1$.

It is noted that we look at the problem in the space
$L^2_{\rho^\ast}$, with weight
\begin{equation*}
\mbox{$\rho^\ast(y) = \mathrm{e}^{-a|y|^{\frac{2k+1}{2k}}}$},
\quad a>0 \,\,\, \mbox{is small enough}.
\end{equation*}
Hence, for this weight, for $y\to+\infty$, we eliminate the roots
such that
\begin{equation*}
\begin{cases}
\mbox{$m <  \frac{k-1}{2} \quad \text{and} \quad m >
\frac{3k-1}{2}$} \quad  \text{for even $k$},\\ \mbox{$m <
\frac{k-1}{2} \quad \text{and} \quad m > \frac{3k-1}{2}$} \quad
\text{for odd $k$},
\end{cases}
\end{equation*}
and for $y\to-\infty$
\begin{equation*}
\begin{cases}
\mbox{$m <  \frac{k}{2} \quad \text{and} \quad m >  \frac{3k}{2}
$} \quad  \text{for even $k$},\\ \mbox{$m <  \frac{k-1}{2}  \quad
\text{and} \quad m > \frac{3k-1}{2} $} \quad \text{for odd $k$}.
\end{cases}
\end{equation*}
These give $2k-1$ conditions in total.

Whilst before balancing the rest of the terms will lead to our
radiation conditions, this is not possible to do this in similar
lines (as for {\bf B}) in the case for ${\bf B}^\ast$.  Note that,
on integration, there holds:
\begin{equation*}
\mbox{$-\frac{1}{2k+1}\,y{\psi^{\ast}}^\prime(y) \sim
\lambda_l\psi^\ast(y)$} \quad \Longrightarrow \quad \psi^\ast_l
(y) \sim A_0y^{-(2k+1)\lambda_l}
\end{equation*}
%%after integration, it can be shown that
%%\begin{equation*}
%%\psi^\ast_l (y) \sim A_0y^{-(2k+1)\lambda_l},
%%\end{equation*}
for some constant $A_0$.  However, this behaviour is perfectly
acceptable in $L^2_{\rho^*}$, and  it is satisfied by the
polynomial adjoint eigenfunctions $\{\psi^\ast_l (y)\}$, which are
our generalized Hermite polynomials.

Instead we look for the remaining two conditions for the problem
as $y\to -\infty$.  In light of the definition of ${\bf B}^\ast$,
we assume to have another type of radiation condition, which for
convenience we denote as the {\em ``Adjoint Radiation Condition"},
such that we
 $$
\fbox{$ \mbox{exclude two bundles  with ``maximal" oscillatory
components at $y=-\infty$.}
 $}
 $$
Again, this extra radiation condition makes the total number of
conditions to be equal to $2k+1$, which is the differential order
of the operator ${\bf B}^*$, so that the eigenvalue problem
becomes
 algebraically consistent \cite{NaimarkI}. In other words, we have
 an algebraic inhomogeneous system of $2k+1$ equations with analytic coefficients
  with $2k+1$ unknowns. Such systems do not have more than a countable set
  of solutions, which are eigenvalues of ${\bf B}^*$, which is
  defined in such a way.

Hence, we now restrict those roots such that $\mathrm{Re}\,\,a_m =
0$. Hence the conditions, as $y\to+\infty$ are now given by the
following distribution of the acceptable coefficients $\{b_m\}$:
\begin{equation*}
\begin{cases}
\mbox{$m \leq  \frac{k}{2} \quad \text{and} \quad m \geq
\frac{3k}{2} $} \quad  \text{for even $k$},\\ \mbox{$m \leq
\frac{k-1}{2}  \quad \text{and} \quad m \geq \frac{3k-1}{2} $}
\quad \text{for odd $k$}.
\end{cases}
\end{equation*}

%%%%%%%%%%%%%%%%%%%%%%%%%%%%%%%%%%%%%%%%%%%%%%%%%%%%%%%%%%
%%%%%%%%%%%%%%%%%%%%%%%%%%%%%%%%%%%%%%%%%%%%%%%%%%%%%%%%%%
%%%%%%%%%%%%%%%%%%%%%%%%%%%%%%%%%%%%%%%%%%%%%%%%%%%%%%%%%%
%%%%%%%%%%%%%%%%%%%%%%%%%%%%%%%%%%%%%%%%%%%%%%%%%%%%%%%%%%
%%%%%%%%%%%%%%%%%%%%%%%%%%%%%%%%%%%%%%%%%%%%%%%%%%%%%%%%%%
%%%%%%%%%%%%%%%%%%%%%%%%%%%%%%%%%%%%%%%%%%%%%%%%%%%%%%%%%%
%%%%%%%%%%%%%%%%%%%%%%%%%%%%%%%%%%%%%%%%%%%%%%%%%%%%%%%%%%
%%%%%%%%%%%%%%%%%%%%%%%%%%%%%%%%%%%%%%%%%%%%%%%%%%%%%%%%%%

\subsection{Calculations for the weights, $\rho(y)$ and $\rho^\ast(y)$}

We now determine the sharp  distance between the principle root
$b_{m_c}$ such that $\mathrm{Re}\,\,b_{m_c} = 0$ and the previous
root $b_{m_c-1}$, where $\mathrm{Re}\,\,b_{m_c-1} > 0$. In doing
so we may find our weight $\rho(y)$ such that it cuts off all
unwanted roots, for which $\mathrm{Re}\,\,b > 0$.

First consider the case as $y\to+\infty$.  For
$\mathrm{Re}\,\,b_{m_c} = 0$, it is known that there is a root
here and this is given by
\begin{equation*}
b_{m_c} = d_k\mathrm{i}, \quad \mbox{where} \quad d_k =
\mbox{$2k\big(\frac{1}{2k+1}\big)^{\frac{2k+1}{2k}}$} \quad
\mbox{and} \quad m_c =
\begin{cases}
\mbox{$\frac{k}{2}$} \quad \text{for even $k$},\\
\mbox{$\frac{k-1}{2}$} \quad \text{for odd $k$}.
\end{cases}
\end{equation*}
%%where $m_c$ is given to be
%%\begin{equation*}
%%m_c =
%%\begin{cases}
%%\mbox{$\frac{k}{2}$} \quad \text{for even $k$},\\
%%\mbox{$\frac{k-1}{2}$} \quad \text{for odd $k$}.
%%\end{cases}
%%\end{equation*}
%%Here
%%\begin{equation*}
%%d_k = \mbox{$2k\big(\frac{1}{2k+1}\big)^{\frac{2k+1}{2k}}$},
%%\end{equation*}
%%as before.

%%??????????????????????????????

%%To RF: there are discrepancies with Thesis; please check all
%%carefully!!!!

%%???????????????????????????????????????

  Hence,
the root $m_c -1$ is given by
\begin{equation*}
m_c-1 =
\begin{cases}
\mbox{$\frac{k-2}{2}$} \quad \text{for even $k$},\\
\mbox{$\frac{k-3}{2}$} \quad \text{for odd $k$}.
\end{cases}
\end{equation*}
This yields
\begin{equation*}
b_{m_c-1} = \mbox{$\big(\cos{\frac{(k-2)\pi}{2k}} +
\mathrm{i}\sin{\frac{(k-2)\pi}{2k}}\big)d_k$},
\end{equation*}
for all $k$. Hence the distance between the two roots (in the real
axis) is
\begin{equation*}
\mbox{$d = d_k\, \cos{\frac{(k-2)\pi}{2k}}>0$}.
\end{equation*}

For $y\to-\infty$, we do not have any roots
$\mathrm{Re}\,\,b_{m_c} = 0$, for any $k$.  However we look at the
distance between the real axis and the next root such that
$\mathrm{Re}\,\,b > 0$.

In this case we find the distance between the real axis and roots
such that
\begin{equation*}
m_c-1= \lfloor \mbox{$\frac{k-1}{2}$} \rfloor,
\end{equation*}
for all $k$.  Hence we have that
\begin{equation*}
b_{m_c-1} =
\begin{cases}
\mbox{$d_k\big[\cos(\lfloor \frac{k-1}{2}\rfloor \frac{\pi}{k}
+\frac{\pi}{2k}) + \mathrm{i}\sin(\lfloor \frac{k-1}{2}\rfloor
\frac{\pi}{k} +\frac{\pi}{2k})\big]$} \quad \text{for even $k$},\\
\mbox{$d_k\big[\cos( \frac{k-1}{2} \frac{\pi}{k} ) +
\mathrm{i}\sin( \frac{k-1}{2} \frac{\pi}{k})\big]$} \quad
\text{for odd $k$}.
\end{cases}
\end{equation*}
This gives the distance between this root and the real axis
\begin{equation*}
d =
\begin{cases}
\mbox{$d_k\,\cos(\lfloor \frac{k-1}{2}\rfloor \frac{\pi}{k}
+\frac{\pi}{2k})$} \quad \text{for even $k$},\\
\mbox{$d_k\,\cos(\frac{k-1}{2} \frac{\pi}{k} )$} \quad \text{for
odd $k$}.\\
\end{cases}
\end{equation*}
This characterises our weighted space $L^2_{\rho}(\mathbb{R})$,
with the exponential weight
\begin{equation*}
\rho(y) = \begin{cases} \mathrm{e}^{a|y|^{\frac{2k+1}{2k}}}
&\text{for $y\leq-1$},\\ \mathrm{e}^{-ay^{\frac{2k+1}{2k}}}
&\text{for $y\geq 1$},
\end{cases}
\end{equation*}
where $a\in(0,2d)$.

Similarly we can find the weight $\rho^\ast(y)$.  We note that the
calculations are exactly the same as in the case for $\rho(y)$,
except a difference in sign when calculating the roots.  This
leads to
\begin{equation*}
d =
\begin{cases}
\mbox{$d_k\,\cos(\lfloor \frac{k-1}{2}\rfloor \frac{\pi}{k}
+\frac{\pi}{2k})$} \quad \text{for even $k$},\\ \mbox{$d_k\,\cos(
\frac{k-1}{2} \frac{\pi}{k} )$} \quad \text{for odd $k$}.\\
\end{cases}
\end{equation*}
as $y\to+\infty$ and
\begin{equation*}
d =
\begin{cases} \mbox{$d_k\, \cos{\big(\frac{k-2}{2}\frac{\pi}{k}\big)} \quad
\text{for even $k$} $},\\
\mbox{$d_k\,\cos{\big(\lfloor\frac{k-2}{2}\rfloor \frac{\pi}{k}
+\frac{\pi}{2k}\big)} \quad \text{for odd $k$} $},
\end{cases}
\end{equation*}
as $y\to-\infty$.  Here the weight $\rho^\ast(y)$ is defined as
 in \ef{rho*}.
%%\begin{equation*}
%%\rho^\ast(y) = {\mathrm e}^{-a|y|^{\frac{2k+1}{2k}}} \quad
%%\text{for all} \quad |y| \ge 1,
%%\end{equation*}
%%for $a\in(0,2d)$.

%%%%%%%%%%%%%%%%%%%%%%%%%%%%%%%%%%%%%%%%%%%%%%%%
%%%%%%%%%%%%%%%%%%%%%%%%%%%%%%%%%%%%%%%%%%%%%%%%
%%%%%%%%%%%%%%%%%%%%%%%%%%%%%%%%%%%%%%%%%%%%%%%%
%%%%%%%%%%%%%%%%%%%%%%%%%%%%%%%%%%%%%%%%%%%%%%%%
%%%%%%%%%%%%%%%%%%%%%%%%%%%%%%%%%%%%%%%%%%%%%%%%
%%%%%%%%%%%%%%%%%%%%%%%%%%%%%%%%%%%%%%%%%%%%%%%%
%%%%%%%%%%%%%%%%%%%%%%%%%%%%%%%%%%%%%%%%%%%%%%%%
%%%%%%%%%%%%%%%%%%%%%%%%%%%%%%%%%%%%%%%%%%%%%%%%
%%%%%%%%%%%%%%%%%%%%%%%%%%%%%%%%%%%%%%%%%%%%%%%%

\end{small}
\end{appendix}

%%%%%%%%%%%\begin{appendix}
%%%%%%%%%%%%%%%%%%%%%%%%%%%%%%%%%%%%%%%%%%%%%%%%%%%%%%%%%%%%%
%%%%%%%%%%%%%%%%%%%%%%%%%%%%%%%%%%%%%%%%%%%%%%%%%%%%%
\begin{appendix}
\section*{Appendix B. Estimates on the fundamental kernel and majorizing
operator}
 %%\label{Sect4.N}
 \label{S3M}
 \setcounter{section}{2}
\setcounter{equation}{0}

\begin{small}

%%%%\section{}
 %%% and comparison}

 Here, we develop some ideas concerning majorizing-comparison
 issues for the odd-order operators involved, which literally do
 not have any connections with the standard Maximum Principle.

%%%%%%%%%%%%%%%%%%%%%%%%%%%%%%%%%%%%%%%%%%%%%%%%%%%%%%%%%%%%%%%%%%%%
%%%%%%%%%%%%%%%%%%%%%%%%%%%%%%%%%%%%%%%%%%%%%%%%%%%%%%
\subsection{Estimates of the rescaled kernel}

Recalling our estimate given by \eqref{initFest}, we now look to
estimate our rescaled fundamental kernel $F(y)$ by
\begin{equation}
\label{estineq}
 \mbox{$
 |F(y)|\leq \bar{D}\bar{F}(y), \quad \text{where}
\quad\bar{F}(y)>0 \quad \text{and} \quad\int\bar{F}(y)=1.
 $}
\end{equation}
Here $\bar{D}$ is a normalisation constant, obviously satisfying
$\bar{D}>1$. There exists infinitely many functions which satisfy
\eqref{estineq}, but since our kernel $F$ is changing sign, it is
not possible to find an optimal analytic function
\begin{equation*}
\bar{F}_{\rm opt}(y) = \mbox{$\omega_1|F(y)|$},
\end{equation*}
where $\omega_1>0$ is a normalisation constant, such that
\begin{equation*}
 \mbox{$
\int \bar{F}_{\rm opt} = 1
 \quad \Longrightarrow \quad
 \mbox{$
\omega_1 = \big(\int |F|\big)^{-1} >1.
 $}
 $}
\end{equation*}
%%In this case
%%\begin{equation*}
%% \mbox{$
%%\omega_1 = \Big(\int |F|\Big)^{-1} >1,
%% $}
%%\end{equation*}
%%since $\int F = 1$.
 However, we can find an analytical
approximation of $\bar{F}_{\rm opt}$ such that \eqref{estineq} is
satisfied, but is non-optimal.  One such function is given by
\begin{equation*}
\bar{F}_{\ast}(y)=\omega_1\mbox{$(1+y^2)^{-\frac{(2k-1)}{8k}}\big(\frac{1}{1+\mathrm{e}^{-y}}+\frac{1}{1+\mathrm{e}^y}\,
\mathrm{e}^{-a(1+y^2)^{\frac{\alpha}{2}}}\big)$}.
\end{equation*}
A sketch of the function $\bar{F}_{\ast}$ is shown by Figure
\ref{estFbar} and comparison with the numerics in Section
\ref{NumConstr} shows how the function may be an upper bound for
$F(y)$.
\begin{figure}[htbp]
\label{estFbar}
\begin{center}
\fbox{\includegraphics[scale=0.6]{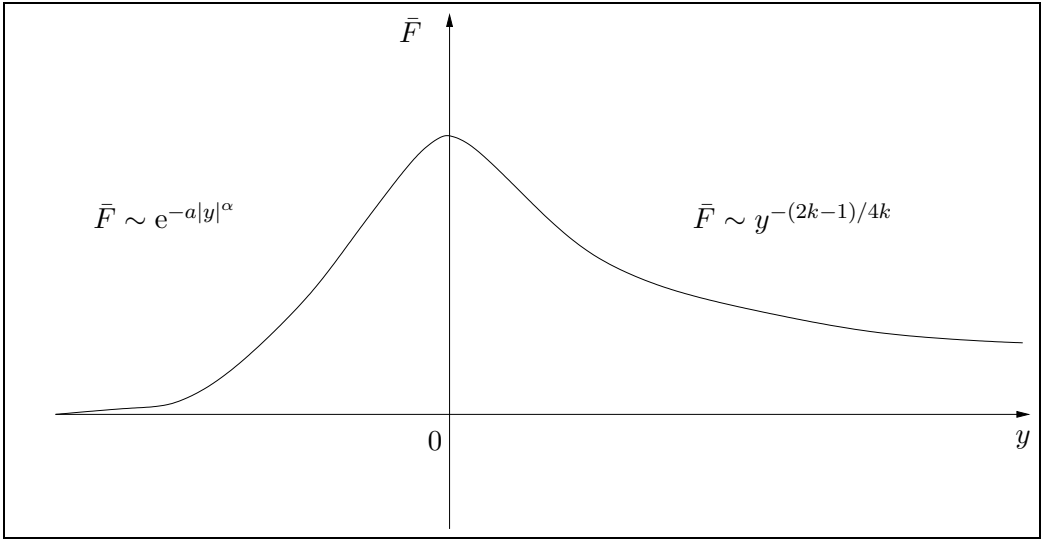}} \caption{\small
Majorizing kernel $\bar{F}(y)$.}
\end{center}
\end{figure}

By using the asymptotic analysis set out in Section
\ref{asymptotics}, we can also find that all higher order
derivatives are estimated by
\begin{equation}
\label{EstDeriv} |D^{\beta}F(y)|\leq
\begin{cases}\bar{c}^{\beta}\beta^{\frac{(\alpha-1)}{\alpha}\beta}\mathrm{e}^{-a|y|^\alpha}
\quad \text{for} \quad y\leq -1,\\
 \bar{c}^{\beta}y^{\frac{\beta}{2k}} \quad \text{for}\quad y\geq 1,
\end{cases}
\end{equation}
where $\bar{c}$ is dependent on $k$ only.

%%%%%%%%%%%%%%%%%%%%%%%%%%%%%%%%%%%%%%%%%%%%%
\subsection{Majorizing kernel and spectral properties}

We now introduce the positive majorizing ``fundamental solution"
\begin{equation*}
\bar{b}(x,t) = t^{-\frac{1}{2k+1}}\bar{F}(y), \quad
y=x/t^{\frac{1}{2k+1}},
\end{equation*}
with the majorizing $\bar F$ relative to the kernel $F$ of
$b(x,t)$ in \ef{kern}. We then end up with  a formal integral
evolution equation written in the standard form
\begin{equation}
\label{EstOpEq} \bar{u}_t={\bf\bar{A}}(t)\bar{u},
\end{equation}
for some linear  operator ${\bf\bar{A}}(t)$, with the
``fundamental solution"  $\bar{b}(x,t)$. This formal
non-autonomous (in time $t$) evolution equation is understood in
the sense that the Cauchy problem for \eqref{EstOpEq} with initial
data
\begin{equation}
\label{initcondmaj} \bar{u}(x,0)=\bar{u}_0(x)\geq 0 \quad
\text{in} \quad \mathbb{R},
\end{equation}
is given by the convolution
\begin{equation}
\label{EstConv}
 \mbox{$
 \bar{u}(x,t)=\bar{{\bf M}}(t)\bar{u}_0(x)\equiv
\bar{b}(t)\ast\bar{u}_0=t^{-\frac{1}{2k+1}}\int
\bar{F}\big((x-z)t^{-\frac{1}{2k+1}}\big)\bar{u}_0(z) \mathrm{d}z.
 $}
\end{equation}
It can be seen from this, that for general kernels $\bar{F}$,
majorizing semigroups do not exist, so equation \eqref{EstOpEq}
does not admit translation in time. This defines the corresponding
majorizing integral equation. For higher-order parabolic
(poly-harmonic) equation, the idea of majorizing integral
operators was introduced and applied in blow-up studied in
\cite{MajOrdOp}.

As before, let us now introduce rescaled variables
\begin{equation*}
\bar{u}(x,t)=t^{-\frac{1}{2k+1}}\bar{w}(y,\tau), \quad
y=xt^{-\frac{1}{2k+1}}, \quad \tau = \ln t.
\end{equation*}
Hence from convolution in \eqref{EstConv}
\begin{equation*}
 \mbox{$
\bar{w}(y,\tau)\equiv \int
\bar{F}(y-z\mathrm{e}^{-\frac{\tau}{2k+1}})\bar{u}_0(z) \,
\mathrm{d}z.
 $}
\end{equation*}
Using Taylor's power series we formally obtain
%%% that
\begin{equation}
%%\begin{split}
\label{estSer}
 \mbox{$
\mbox{$\bar{F}\big(y-z\mathrm{e} ^{-\frac{\tau }{2k+1}}\big)$} =
\sum\limits _{(\beta )} \mbox{$\mathrm{e}^{-\frac{\beta \tau
}{2k+1}}\frac{(-1)^{\beta }}{\beta !}D_y ^{\beta
}\bar{F}(y)z^{\beta }$}
  \equiv \sum\limits _{(\beta )} \mbox{$\mathrm{e}
^{-\frac{\beta \tau }{2k+1}}\frac{1}{\sqrt{\beta !}}\, \bar{\psi}
_{\beta }(y)z^{\beta }$},
 $}
 %%\end{split}
\end{equation}
%%where
\begin{equation*}
 \mbox{where} \quad
\mbox{$\bar{\psi} _{\beta }(y) = \frac{(-1)}{\sqrt{\beta
!}}^{\beta }D_y ^{\beta }\bar{F}(y)$}.
\end{equation*}
The convergence of \eqref{estSer} on bounded intervals is
guaranteed by the estimates of $D_y ^{\beta }\bar{F}(y)$ given in
\eqref{EstDeriv}. A proper convergence in $L^2_\rho$ can be also
included as above.

The solution can then be represented by
\begin{equation*}
 \mbox{$
\bar{w}(y,\tau ) = \sum\limits _{(\beta )} \mathrm{e}
^{-\frac{\beta \tau }{2k+1}}\bar{M}_{\beta }(u_0)\bar{\psi}
_{\beta }(y), \,\,
%%\mbox{where}
 \bar{\lambda} _{\beta } = -\frac{\beta }{2k+1},
 \,\, \mbox{where}
 \,\,\,
\bar{M}_{\beta }(\bar{u}_0) = \mbox{$\frac{1}{\sqrt{\beta !}}$}
\int\limits_{\mathbb{R}}z^{\beta }\bar{u_0}(z) \, \mathrm{d}z
 $}
 \quad
\end{equation*}
%%where we define $\bar{\lambda} _{\beta } = -\frac{|\beta |}{2k+1}$
%%and
%%\begin{equation*}
%%\bar{M}_{\beta }(\bar{u}_0) = \mbox{$\frac{1}{\sqrt{\beta !}}$}
%%\int_{\mathbb{R}}z^{\beta }\bar{u_0}(z) \, \mathrm{d}z,
%%\end{equation*}
are the corresponding moments of the initial data.
\begin{proposition}
There exists some formal operator ${\bf \bar{B}}$, such that
\begin{equation*}
\bar{w}_\tau = {\bf \bar{B}}\bar{w}
\end{equation*}
and this induces the majorizing semigroup $\{\mathrm{e}^{{\bf
\bar{B}}\tau}\}$, formulated by the rescaled variables.
\end{proposition}

It follows from \eqref{wsol} that ${\bf \bar{B}}$ has point
spectrum given by
 %%\begin{equation*}
$\sigma_p({\bf \bar{B}})=\{\bar{\lambda}_\beta\}$
 %%%\end{equation*}
and corresponding eigenfunctions are thus given by $\bar \Phi=
\{\bar{\psi}_\beta\}$. With a standard definition of the space of
$\bar \Phi$-closure $\tilde L^2_\rho$ (see Section \ref{S4.3}),
this becomes the only spectrum of $\bar \BB$.

In order to find a full expansion form of the majorizing
semigroup, we once again perform another rescaling given by
\begin{equation*}
\bar{u} = (1+t)^{-\frac{1}{2k+1}}\bar{w}, \quad y =
x(1+t)^{-\frac{1}{2k+1}}, \quad \tau =
\ln{(1+t)}:\mathbb{R}_+\rightarrow\mathbb{R}_+.
\end{equation*}
Then rescaling the convolution we obtain
\begin{equation*}
%%\begin{split}
 \mbox{$
\bar{w}(y,\tau ) = \mathrm{e} ^{{\bf \bar{B}}\tau}\bar{u}_0
\equiv(1-\mathrm{e} ^{-\tau
})^{-\frac{1}{2k+1}}\int\limits_{\mathbb{R}}
\mbox{$\bar{F}\big((y-z\mathrm{e} ^{-\frac{\tau
}{2k+1}})(1-\mathrm{e}
^{-\tau})^{-\frac{1}{2k+1}}\big)\bar{w}_0(z) \, \mathrm{d} z$}.
 %%\end{split}
 $}
\end{equation*}
Using Taylor expansions, we find that the solution is given by
\begin{equation*}
\begin{split}
 \mbox{$
\bar{w}(y,\tau ) = (1-\mathrm{e}
^{-\tau})^{-\frac{1}{2k+1}}\sum\limits _{(\mu ,\, \nu) }
 $}
&\mbox{$\frac{(-1)}{\mu !\nu !}^{\mu }D^{\nu
}\bar{F}(0)\frac{1}{(\nu - \mu )!}\,y^{\nu - \mu }(\mathrm{e}
^{\tau }-1)^{-\frac{\mu }{2k+1}}$}\\&\mbox{$\times (1-\mathrm{e}
^{-\tau })^{-\frac{\nu }{2k+1}}$}
 \mbox{$
\int\limits_{\mathbb{R}}
 $}
 z^{\mu }\bar{w}_0(z) \, \mathrm{d} z.
\end{split}
\end{equation*}
This again shows the discrete spectrum, as well as some traces of
generalized Hermite polynomials $\{\bar \psi_\b^*\}$.

The adjoint operator $\bar {\bf B}^\ast$ with such polynomial
eigenfunctions $\{\bar{\psi}_\beta^\ast\}$ occurs if we use the
blow-up scaling
\begin{equation*}
\bar{u}(x,t)=\bar{w}(y,\tau), \quad y=x(1-t)^{-\frac{1}{2k+1}},
\quad \tau=-\ln(1-t).
\end{equation*}
We thus obtain
\begin{equation*}
 \mbox{$
\bar{w}(y,\tau )=(1-\mathrm{e} ^{-\tau })\int\limits_{\mathbb{R}}
\mbox{$\bar{F}\big((y\mathrm{e} ^{-\frac{\tau
}{2k+1}}-z)(1-\mathrm{e}
^{-\tau})^{-\frac{1}{2k+1}}\big)\bar{w}_0(z) \, \mathrm{d} z$}.
%%\\[5mm]
 $}
\end{equation*}
Using Taylor expansion yields
\begin{equation*}
\begin{split}
 \mbox{$
\mbox{$\bar{w}(y,\tau)=(1-\mathrm{e} ^{-\tau
})^{-\frac{1}{2k+1}}$}\sum\limits _{(\beta ,\, \nu)}
 $}
&\mbox{$\frac{(-1)^{\beta }}{\beta !\nu !}(1-\mathrm{e}
^{-\tau})^{-\frac{\beta }{2k+1}}(\mathrm{e} ^{\tau
}-1)^{-\frac{\nu }{2k+1}}$}\\&\mbox{$\times D^{\nu
}F(0)\frac{1}{(\nu - \beta )!}y^{\nu - \mu }$}
 \mbox{$
\int\limits_{\mathbb{R}}
 $}
 \mbox{$z^{\beta }\bar{w}_0(z) \, \mathrm{d}z$}.
\end{split}
\end{equation*}
As above, further expanding of exponential terms here leads to the
eigenfunction expansion, which determines finite generalized
Hermite polynomials. This representation is rather technical and
we do not present and analyze it since do not aim using those
polynomials in what follows. Our main application of the
majorizing operators is as follows:

%%%%%%%%%%%%%%%%%%%%%%%%%%%%%%%%%%%%%%%%%%%%%%%555
\subsection{Comparison with majorizing problem}

We look to see how the majorizing kernel relates to our real
solution $u(x,t)$. By looking at the convolution
\eqref{convolution} for the linear PDE, we can easily see that
\begin{equation*}
|u(x,t)|\leq|b(t)|\ast|u_0(x)|.
\end{equation*}
Now looking at our estimate of the majorizing kernel, $
\bar{b}(x,t)$, we have that
\begin{equation*}
|u(x,t)|\leq \bar{D}\bar{b}(x,t)\ast|u_0(x)|\leq \bar{b}(t)\ast
\bar{u}_0(x),
\end{equation*}
where we have to assume the following inequality for initial data:
\begin{equation}
\label{U0Ineq} \bar{D}|u_0(x)|\leq \bar{u}_0(x) \quad \text{in}
\quad \mathbb{R}.
\end{equation}

\begin{proposition}
If \eqref{initcondmaj} and \eqref{U0Ineq} hold then
\begin{equation*}
|u(x,t)|\leq\bar{u}(x,t) \quad \text{in} \quad
\mathbb{R}\times\mathbb{R}_+.
\end{equation*}
\end{proposition}

We introduce the majorizing kernel, since the structural behaviour
of self-similar solutions of the majorizing evolution equation
describes essential features of the solutions to the original PDE.
 Since we have estimated our solution, it can be possible to find properties of the
real solution, which otherwise would be more difficult to achieve.

The linear semigroup for {\bf B} is not order-preserving, since
the kernel $F$ is oscillatory.  Therefore we must have
$\bar{D}>1$, which gives the order deficiency of the linear
operator ${\bf B}$ and of the linear convolution operator
$\bar{M}(t)$, given in \eqref{EstConv}.  The actual defect is, in
fact, characterized by $\bar{D}-1>0$. The linear semigroup for
{\bf B} would only be order-preserving if $F$ did not change sign,
since then we would have that $F=\bar{F}$ and so we would have
$\bar{D}=1$ and so the defect would be zero. As we have shown,
this is not possible for any $k\geq1$.

\end{small}
\end{appendix}

\end{document}